\documentstyle[11pt]{article}
\title{On the Deligne-Simpson problem\footnote{Research partially supported 
by INTAS grant 97-1644}}
\author{Vladimir Petrov Kostov\\ \\ \hspace{7cm}{\sl To prof. V.I.Arnold}} 
\date{}
\newtheorem{tm}{Theorem}
\newtheorem{lm}[tm]{Lemma}
\newtheorem{cor}[tm]{Corollary}
\newtheorem{prop}[tm]{Proposition}
\newtheorem{rem}[tm]{Remark}
\newtheorem{rems}[tm]{Remarks}
\newtheorem{defi}[tm]{Definition}
\newtheorem{ex}[tm]{Example}
\begin{document}
\maketitle

\begin{abstract}
The Deligne-Simpson problem is formulated like this: {\em give necessary and 
sufficient conditions for the choice of the conjugacy classes 
$C_j\subset SL(n,{\bf C})$ or $c_j\subset sl(n,{\bf C})$ so that there exist 
irreducible $(p+1)$-tuples of matrices $M_j\in C_j$ or $A_j\in c_j$ 
satisfying the equality $M_1\ldots M_{p+1}=I$ or $A_1+\ldots +A_{p+1}=0$}. 

We solve the problem for generic eigenvalues with the exception of the 
case of matrices $M_j$ when the greatest common divisor of the 
numbers $\Sigma _{j,l}(\sigma )$ of 
Jordan blocks of a given matrix $M_j$, with a given 
eigenvalue $\sigma$ and of a given size $l$ (taken over all $j$, $\sigma$, 
$l$) is $>1$. 
Generic eigenvalues are defined by explicit algebraic 
inequalities. For such eigenvalues there exist no reducible $(p+1)$-tuples.  

The matrices $M_j$ and $A_j$ are interpreted as monodromy operators of 
regular linear systems and as matrices-residua of fuchsian ones on Riemann's 
sphere. \\  

{\bf Key words:} generic eigenvalues, (poly)multiplicity vector, 
corresponding Jordan normal forms, monodromy operator.

{\bf AMS classification index:} 15A30, 15A24, 20G05
\end{abstract}
\tableofcontents 

\section{Introduction}
\subsection{Formulation of the problem}

The problem which is dealt with in the present paper is connected with 
the theory of {\em fuchsian linear systems} on Riemann's sphere, i.e. 
meromorphic linear systems of differential equations with logarithmic 
poles. Such a 
system can be represented as follows:

\begin{equation}\label{Fuchs}
\dot{X}=\left( \sum _{j=1}^{p+1}\frac{A_j}{t-a_j}\right) X
\end{equation}
where $"^."$ denotes d/d$t$, $t\in {\bf C}P^1$, $A_j\in gl(n,{\bf C})$, 
the $p+1$ poles $a_j$ 
are distinct and we assume without restriction that there is no pole at 
infinity. This last condition implies that the sum of the 
{\em matrices-residua} $A_j$ is 0:

\begin{equation}\label{A_j}
\sum _{j=1}^{p+1}A_j=0
\end{equation}

Fuchsian systems are a particular 
case of {\em regular} systems, i.e. linear systems whose solutions when 
restricted to sectors centered at the poles have a moderate growth 
rate when the argument tends to the pole: $||X||=O(|t-a_j|^{N_j})$ for 
some $N_j\in {\bf R}$.

It is more convenient to consider the dependent variables $X$ as an 
$n\times n$-matrix, i.e. to consider simultaneously $n$ linearly independent 
vector-solutions. This is what we do. 

Fix a base point $a_0$ different from the poles $a_j$. Fix the value 
$B\in GL(n,{\bf C})$ of the solution for $t=a_0$. For each pole $a_j$ define 
a closed contour $\Gamma _j$ containing $a_0$ and freely homotopic to a 
positive loop around $a_j$. 
The contour $\Gamma _j$ consists of a line segment $[a_0,x_j]$ 
where $x_j$ is close to $a_j$, of the circumference $\Theta _j$ 
centered at $a_j$, passing through $x_j$ and circumventing $a_j$ 
counterclockwise (we choose $x_j$ so close to $a_j$ that no other pole of 
the system lies  
inside or on $\Theta _j$), and of the line segment $[x_j,a_0]$. 
We assume 
that for $i\neq j$ one has $\Gamma _i\cap \Gamma _j=\{ a_0\}$ and that when 
one turns around $a_0$ 
clockwise the indices of the contours change from 1 to $p+1$. 

The value at $a_0$ of the 
analytic continuation of the solution along $\Gamma _j$ is 
representable 
in the form $BM_j$. The matrix $M_j\in GL(n,{\bf C})$ is by definition the 
matrix of the 
{\em monodromy operator} corresponding to the class of homotopy equivalence 
of the contour $\Gamma _j$. For a choice of contours $\Gamma _j$ 
like above one has 
 
\begin{equation}\label{M_j}
M_1\ldots M_{p+1}=I
\end{equation}
which is the multiplicative analog of (\ref{A_j}). The monodromy operators 
generate the {\em monodromy group} which is invariant under linear 
transformations of the dependent variables meromorphically depending on the 
time (and, up to conjugacy, the only such invariant). 

\begin{rem}
Note that with this definition the monodromy group is an 
{\em antirepresentation} 
of $\pi _1({\bf C}P^1\backslash \{ a_1,\ldots ,a_{p+1}\})$ into 
$GL(n,{\bf C})$ (because to the product of contours $\Gamma _i\Gamma _j$ there 
corresponds the monodromy operator $M_jM_i$; product in the sense of 
concatenation). To obtain a representation one has to consider the matrices 
$M_j^{-1}$. In the paper we refer to the $(p+1)$-tuples of matrices also 
as to representations.
\end{rem} 

The {\em Deligne-Simpson problem (DSP)} is formulated as follows:

{\em For what $(p+1)$-tuples of conjugacy classes 
$C_j\subset SL(n,{\bf C})$ do there exist irreducible 
$(p+1)$-tuples of matrices $M_j\in C_j$ satisfying (\ref{M_j}) ?} 
(multiplicative version).

{\em For what $(p+1)$-tuples of conjugacy classes $c_j\subset sl(n,{\bf C})$ 
do there exist irreducible 
$(p+1)$-tuples of matrices $A_j\in c_j$ satisfying 
(\ref{A_j}) ?} (additive version).

We give the basic result in Subsection \ref{thebasicresult} 
(see Theorems~\ref{generic}, \ref{genericbis} and \ref{nongenericevs}) 
followed by a plan of the paper, after introducing some definitions in the 
next three subsections.

\begin{rems}
1) "Irreducible" means "not having a common proper invariant 
subspace"; in other words, impossible to bring the $(p+1)$-tuple to a 
block upper-triangular form by simultaneous conjugation. The problem could be 
formulated without the requirement of irreducibility and it would be another 
problem which we do not consider here. However, we consider the problem with 
"irreducible" replaced by "with trivial centralizer", see 
Theorem~\ref{nongenericevs}.

2) In the multiplicative version (i.e. for matrices $M_j$) the problem was 
formulated by P.Deligne, and 
C.Simpson was the first to obtain results towards its solution, see 
\cite{Si1} and \cite{Si2}. Simpson's result is cited in 
Remarks~\ref{aftergeneric}. 

3) The case of nilpotent matrices $A_j$ and of unipotent matrices $M_j$ was 
considered by the author in \cite{Ko1} and \cite{Ko2}. One of the results from 
\cite{Ko2} is used in the present paper, see Theorem~\ref{nilpunip}.

4) We treat the two versions of the problem (additive 
and multiplicative) in parallel. The multiplicative should be considered 
as more important because the monodromy group is a meromorphic invariant 
(up to conjugacy) whereas the $(p+1)$-tuple of matrices-residua is not.

5) One can replace $sl(n,{\bf C})$ or $SL(n,{\bf C})$ by  
$gl(n,{\bf C})$ or $GL(n,{\bf C})$; this is what we do when we solve the 
problem because in the process of solving it one encounters matrices $A_j$ 
and $M_j$ not from $sl(n,{\bf C})$ (resp. not from $SL(n,{\bf C})$).

6) Notice that the Deligne-Simpson problem is formulated in a purely algebraic 
way, without reference to fuchsian or regular systems. Yet they explain the 
interest in solving it.
\end{rems} 

{\bf Acknowledgement.} The author is deeply grateful to P.Deligne and N.Katz 
the exchange of e-mail messages and letters with whom helped him avoid some 
mistakes, and also to J. Brian\c{c}on, A.Galligo, O.Gleizer, O.Laudal, 
Ph.Maisonobe, M.Merle and C.Simpson the discussions with whom were also 
helpful. 

\subsection{Generic eigenvalues\protect\label{gen}}

\begin{defi}
Call {\em Jordan normal form of size $n$} a family 
$J^n=\{ b_{i,l}\}$ ($i\in I_l$, $I_l=\{ 1,\ldots ,s_l\}$, $l\in L$) of 
positive integers $b_{i,l}$ 
whose sum is $n$. The set $L$ is the one of indices of eigenvalues 
(all distinct) and 
$I_l$ is the set of indices of Jordan blocks with the $l$-th eigenvalue; 
$b_{i,l}$ is the 
size of the $i$-th block with this eigenvalue. An $n\times n$-matrix 
$Y$ has the Jordan normal form $J^n$ (notation: $J(Y)=J^n$) if its eigenvalues 
can be indexed by $L$ in such a way that to its 
distinct eigenvalues $\delta _l$ ($l\in L$) there belong Jordan blocks of 
sizes $b_{i,l}$. We assume that for every fixed $l$ one has 
$b_{1,l}\geq b_{2,l}$ $\geq$ $\ldots$ $\geq b_{s_l,l}$.
\end{defi} 

\begin{rem}
The basic result of this paper depends  
actually not on the conjugacy classes but only on the Jordan normal forms of 
the matrices $A_j$ or $M_j$ provided that the eigenvalues remain generic, 
see the definition of generic eigenvalues below.
\end{rem} 

We assume that the following necessary conditions for existence of 
irreducible $(p+1)$-tuples of matrices $A_j$ or $M_j$ (satisfying (\ref{A_j}) 
or (\ref{M_j})) hold:

\begin{equation}\label{tracedet} 
\sum _{j=1}^{p+1}{\rm Tr}c_j=0~~,~~\prod _{j=1}^{p+1}\det C_j=1
\end{equation} 
Denote by $\lambda _{k,j}$, $\sigma _{k,j}$ the eigenvalues of 
$A_j$, $M_j$ (they are {\em not} presumed distinct; a {\em multiplicity} of 
an eigenvalue is by definition the number of eigenvalues equal to it 
including the eigenvalue itself). When $A_j$ are 
the matrices-residua of a fuchsian system with 
monodromy operators $M_j$, then one has 
$\sigma _{k,j}=\exp (2\pi i\lambda _{k,j})$. Equation (\ref{tracedet}) 
admits the following equivalent form:  

\begin{equation}\label{sumEVszero} 
\sum _{k=1}^n\sum _{j=1}^{p+1}\lambda _{k,j}=0~~,~~ 
\prod _{k=1}^n\prod _{j=1}^{p+1}\sigma _{k,j}=1
\end{equation}
(see (\ref{A_j}) and (\ref{M_j})). From now on we presume that the 
eigenvalues satisfy conditions (\ref{sumEVszero}). Call 
{\em non-genericity relation} any equality of the form 

\[ \sum _{j=1}^{p+1}\sum _{k\in \Phi _j}\lambda _{k,j}=0~~
{\rm ,~resp.}~~
\prod _{j=1}^{p+1}\prod _{k\in \Phi _j}\sigma _{k,j}=1~~~~~~~~~~~~~~~~~~~~
(\gamma )\]
where the sets $\Phi _j$ contain one and the same number $s$ of indices 
(with $1<s<n$) for all $j$. Eigenvalues satisfying none of these 
relations are called {\em generic}. 

Eigenvalues $\lambda _{k,j}$ satisfying none of the relations 
$(\gamma )$ modulo ${\bf Z}$ are called {\em strongly generic}. If the 
eigenvalues $\lambda _{k,j}$ are strongly generic, then the eigenvalues 
$\sigma _{k,j}$ are generic.

\begin{rems}
1) Reducible (i.e. block upper-triangular up to conjugacy) 
$(p+1)$-tuples of matrices exist only for non-generic eigenvalues and if the 
$(p+1)$-tuple is in block upper-triangular form, then the 
eigenvalues of its restriction to each diagonal block satisfy some 
relation $(\gamma )$. 

2) If in $(\gamma )$ one replaces each of the sets $\Phi _j$ by its 
complement  
in $\{ 1,\ldots ,n\}$, then one obtains an equivalent non-genericity relation. 
\end{rems}

\subsection{The quantities $r_j$ 
and $d_j$\protect\label{rd}}

For a conjugacy class $C$ in $GL(n,{\bf C})$ or $gl(n,{\bf C})$ denote by 
$d(C)$ its dimension and for a matrix $Y$ from $C$ set 
$r(C):=\min _{\lambda \in {\bf C}}{\rm rank}(Y-\lambda I)$. The integer 
$n-r(C)$ is the maximal number of Jordan blocks of $J(Y)$ with one and the 
same eigenvalue. Set $d_j:=d(C_j)$ (resp. $d(c_j)$), $r_j:=r(C_j)$ 
(resp. $r(c_j)$). The quantities 
$r(C)$ and $d(C)$ depend only on the Jordan normal form $J(Y)=J^n$, 
so we write sometimes $r(J^n)$ and $d(J^n)$.

The following proposition was proved in the multiplicative version by 
C.Simpson in \cite{Si1}. We give a proof for both versions here:

\begin{prop}\label{d_jr_j}
A necessary condition for the existence of irreducible $(p+1)$-tuples 
satisfying respectively (\ref{A_j}) or (\ref{M_j}) is the following couple 
of inequalities to hold:

\[ \begin{array}{rcllllc}d_1+\ldots +d_{p+1}&\geq &2n^2-2&&&&(\alpha _n)\\
{\rm for~all~}j~(r_1+\ldots +\widehat{r_j}+\ldots +r_{p+1})&
\geq &n&&&&(\beta _n)
\end{array}\] 
\end{prop}

Condition $(\beta _n)$ is generalized in Proposition~\ref{necessarynongen}. 
Both propositions are proved in Section~\ref{proofofd_jr_j}.

\begin{rem}
Neither of the inequalities $(\alpha _n)$ and $(\beta _n)$ 
follows from the other one. If 
$n\geq 4$ is even, $p=2$, the matrices $A_j$ are diagonalizable and the
multiplicity vectors of the eigenvalues of $A_1$, $A_2$, $A_3$ equal 
respectively
$(1,\ldots ,1)$, $(n/2,n/2)$, $(n/2,n/2)$, then $(\beta _n)$ holds while
$(\alpha _n)$
does not. If $n$ is even, $n+2=2p$, $p\geq 3$ and the multiplicity vectors of
the eigenvalues of the diagonalizable 
matrices $A_j$ equal $(1,\ldots ,1)$, $(n-1,1)$, 
$\ldots$, $(n-1,1)$, then $(\alpha _n)$ holds while $(\beta _n)$ does not.
\end{rem}

\begin{defi}
In the additive version we say that the DSP is {\em solvable} 
(resp. {\em weakly solvable}) for 
given Jordan normal forms $J_j^n$ and for given eigenvalues if there exists 
an {\em irreducible} $(p+1)$-tuple of matrices $A_j$ satisfying condition 
(\ref{A_j}) with $J(A_j)=J_j^n$ and with the given eigenvalues (resp. if 
there exists such a $(p+1)$-tuple of matrices $A_j$ 
with a {\em trivial centralizer}). In the multiplicative version one 
replaces in the definition the matrices $A_j$ by matrices $M_j$ 
satisfying condition (\ref{M_j}).
\end{defi}

\begin{prop}\label{necessarynongen}
1) For prescribed conjugacy classes $c_j$ of the matrices $A_j$  
(not necessarily diagonalizable or with generic eigenvalues)  
denote by $r_j(b)$ the rank of the matrix $A_j-bI$, $b\in {\bf C}$, 
$A_j\in c_j$. A necessary condition for the solvability of the DSP 
for matrices $A_j\in c_j$ is the following inequality:

\[ \min _{b_i\in {\bf C},b_1+\ldots +b_{p+1}=0}
(r_1(b_1)+\ldots +r_{p+1}(b_{p+1}))\geq 2n\]

2) A necessary condition for the solvability of the DSP for matrices 
$M_j$ from prescribed conjugacy classes $C_j$ 
(not necessarily diagonalizable or with generic eigenvalues) 
is the following inequality:

\[ \min _{b_j\in {\bf C}^*, b_1\ldots b_{p+1}=1}({\rm rk}(b_1M_1-I)+\ldots 
+{\rm rk}(b_{p+1}M_{p+1}-I))\geq 2n\]
\end{prop}

\begin{rems}
1) Remind that the $(p+1)$-tuples $(A_1,\ldots ,A_{p+1})$ and 
$(A_1-b_1I,\ldots ,A_{p+1}-b_{p+1}I)$, $b_1+\ldots +b_{p+1}=0$ are 
simultaneously irreducible. The above minimum is obtained for a 
$(p+1)$-tuple $(b_1,\ldots ,b_{p+1})$ in which at least $p$ of the numbers 
$b_j$ are eigenvalues of the corresponding matrices $A_j$. A similar remark 
is true for the matrices $M_j$ as well. 

2) The proposition generalizes condition $(\beta _n)$. In the case of generic 
eigenvalues it coincides with it, in the case of non-generic ones it 
implies it and is stronger than it.
\end{rems} 

For the formulation of the basic result it will be essential whether the 
quantities $r_j$ satisfy the inequality

\[ (r_1+\ldots +r_{p+1})\geq 2n~~~~~~~~~~~~~~~~(\omega _n)\]

\begin{rems}
1) Evidently, condition $(\omega _n)$ implies condition $(\beta _n)$ 
(because $r_j<n$), but condition $(\omega _n)$ is not 
necessary for the existence of irreducible $(p+1)$-tuples (e.g. 
for $p=2$, $n=2$ there exist irreducible triples of matrices each with two 
distinct eigenvalues, i.e. $r_1=r_2=r_3=1$ and $(\omega _2)$ does not hold). 

2) Condition $(\omega _n)$ arises when the Deligne-Simpson problem is 
considered for nilpotent matrices $A_j$ and for unipotent matrices $M_j$. 
For such matrices the eigenvalues are ``mostly non-generic'', i.e. they 
satisfy all possible non-genericity relations. Moreover, for 
such matrices condition $(\omega _n)$ is necessary for the existence of 
irreducible $(p+1)$-tuples -- it coincides with the necessary conditions 
from Proposition~\ref{necessarynongen}. It turns 
out that condition $(\omega _n)$ is almost sufficient as well 
in the following sense.
\end{rems}

\begin{defi}\label{nicerepr}
A $(p+1)$-tuple of matrices $A_j$ whose sum is 0 or of 
matrices $M_j$ whose product is $I$ is said to define a {\em nice} 
representation (or for short, to be {\em nice}) 
if its centralizer is trivial and either the 
$(p+1)$-tuple is irreducible or it is reducible and one can conjugate it to 
a block upper-triangular form in which the diagonal 
blocks are all of sizes $g_i>1$ and define irreducible representations. 
Thus the matrix algebra ${\cal A}$ defined 
by the matrices $A_j$ or $M_j-I$ contains a non-degenerate matrix 
(i.e. with non-zero 
determinant) -- by the Burnside theorem the restrictions to the diagonal 
blocks of ${\cal A}$ equal $gl(g_i,{\bf C})$.
\end{defi} 

\begin{tm}\label{nilpunip}
If for the nilpotent (resp. unipotent) conjugacy classes $c_j$ (resp. $C_j$) 
condition $(\omega _n)$ holds and the following four particular cases are 
avoided, then there exist nice $(p+1)$-tuples of matrices $A_j\in c_j$ 
whose sum 
is 0 (resp. of matrices $M_j\in C_j$ whose product is $I$).  
In the four particular cases each conjugacy class has Jordan blocks of one and 
the same size (denoted by $l_j$). The cases are:

1) $n=2k$, $k>1$, $p=3$, $l_1=l_2=l_3=l_4=2$;

2) $n=3k$, $k>1$, $p=2$, $l_1=l_2=l_3=3$;

3) $n=4k$, $k>1$, $p=2$, $l_1=l_2=4$, $l_3=2$;

4) $n=6k$, $k>1$, $p=2$, $l_1=6$, $l_2=3$, $l_3=2$.
\end{tm}

The above theorem is part of Theorem 34 from \cite{Ko2}. 
If another four particular cases are avoided, then there exist $(p+1)$-tuples 
defining irreducible representations, see \cite{Ko2}.

\subsection{(Poly)multiplicity vectors\protect\label{poly}}

\begin{defi}
A {\em polymultiplicity vector (PMV)} is by definition a 
$(p+1)$-tuple of 
{\em multiplicity vectors (MVs)}, i.e. vectors whose components are 
non-negative integers,  
their sum (called the {\em length} of the MV) being equal to $n$. 
The MVs and PMVs with which we deal in this paper are defined by the 
multiplicities of the eigenvalues of the matrices $A_j$ or $M_j$. (We allow 
zero components for the sake of convenience.) 
Call a PMV {\em simple} (resp. {\em non-simple}) if 
the greatest common divisor of all its non-zero 
components equals 1 (resp. if not).
\end{defi}

\begin{rem}\label{tricky}
In the case of matrices $A_j$ generic eigenvalues exist only 
for simple PMVs. Indeed, if all multiplicities are divisible by 
$1<q\in {\bf Z}$, then the sum of all eigenvalues with multiplicities 
divided by $q$ equals 0 which 
is a non-genericity relation $(\gamma )$. In the case of matrices $M_j$ 
the divisibility by $q$ would imply only that the product of all  
eigenvalues with multiplicities divided by $q$ equals one of the roots of 
unity of $q$-th order, not necessarily 1, i.e. a non-genericity 
relation might or might not hold. 
However, in both cases there exist generic eigenvalues for every simple PMV.
Generic eigenvalues form a Zariski open dense subset in the set of all 
eigenvalues with a fixed simple PMV. The latter set is a linear 
space in the additive version and a non-singular variety in the 
multiplicative one, see condition (\ref{sumEVszero}).  
\end{rem}

For a PMV of length $n$ we use the notation 
$\Lambda ^n=$$(\Lambda ^n_1, \ldots , \Lambda ^n_{p+1})$ where 
$\Lambda ^n_j$ are the MVs.  
For diagonalizable matrices the MV $\Lambda ^n_j$ to 
have only one component implies $A_j$ or $M_j$ to be scalar. 

Set $\Lambda ^n_j=(m_{1,j},\ldots ,m_{k_j,j})$, $m_j=\max _im_{i,j}$. 
For a diagonalizable conjugacy class one has $r_j=n-m_j$, 
$d_j=n^2-\sum _{i=1}^{k_j}(m_{i,j})^2$. In accordance with the corresponding 
definitions for Jordan normal forms, we say that 
$\Lambda ^n$ satisfies Condition $(\alpha _n)$  (Condition $(\beta _n)$, 
Condition $(\omega _n)$) if 
$d_1+\ldots +d_{p+1}\geq 2n^2-2$ (if 
$\min _j(r_1+\ldots +r_{p+1}-r_j)\geq n$, if $r_1+\ldots +r_{p+1}\geq 2n$).

\subsection{Formulation of the basic result 
\protect\label{thebasicresult}}

For a given $(p+1)$-tuple $(J_1^n,\ldots , J_{p+1}^n)$ of Jordan normal 
forms with $n>1$ (the upper index indicates the size of the 
matrices), which satisfies condition $(\beta _n)$ 
and does not satisfy condition $(\omega _n)$ set $n_1=r_1+\ldots +r_{p+1}-n$. 
Hence, $n_1<n$ and $n-n_1\leq n-r_j$ for all $j$. Define 
the $(p+1)$-tuple of Jordan normal forms $J_j^{n_1}$ as follows: to obtain 
the Jordan normal form $J_j^{n_1}$ 
from $J_j^n$ one chooses one of the eigenvalues of $J_j^n$ with 
greatest number $n-r_j$ of Jordan blocks, then decreases  
by 1 the sizes of the $n-n_1$ {\em smallest} Jordan blocks with this 
eigenvalue and deletes the Jordan blocks of size 0. We write this 
symbolically in the form 

\[ \Psi :(J_1^n,\ldots ,J_{p+1}^n)\mapsto (J_1^{n_1},\ldots ,J_{p+1}^{n_1})\]

For a given $(p+1)$-tuple of Jordan normal forms 
$J^n=(J_1^n,\ldots ,J_{p+1}^n)$ define a sequence of $(p+1)$-tuples of 
Jordan normal forms 
$J^{n_{\nu}}$, $\nu =0$, $\ldots$, $s$ by iterating the map $\Psi$ as long 
as it is defined; we set $n_0=n$.

\begin{rem}
Notice that $n>n_1>\ldots >n_s$ (we define $n_1$ only if $J^n$ 
does not 
satisfy condition $(\omega _n)$, hence, $n_1<n$ etc.).
\end{rem}

\begin{tm}\label{generic}
For given Jordan normal forms $J_j^n$ with a simple PMV and for generic 
eigenvalues the DSP is solvable for matrices $A_j$ or $M_j$ 
if and only if the following two conditions hold:

{\em i)} The $(p+1)$-tuple of Jordan normal forms $J_j^n$ satisfies the 
inequality $(\beta _n)$;

{\em ii)} Either the $(p+1)$-tuple of Jordan normal forms $J_j^{n_s}$ 
satisfies the inequality $(\omega _{n_s})$ or $n_s=1$.
\end{tm}

\begin{rems}\label{aftergeneric}
1) The theorem holds whichever choice of eigenvalue with 
maximal number of Jordan blocks is made to define $\Psi$ (one 
choice is sufficient, it holds automatically for all other choices; the 
numbers $n_k$ are the same for all choices). 
See Remarks~\ref{psiinvariant}.

2) Condition $(\alpha _n)$ does not appear explicitly 
in the formulation of the theorem. However, it is implicitly present 
because if the 
$(p+1)$-tuple of Jordan normal forms $J_j^{n_{\nu}}$ satisfies condition 
$(\alpha _{n_{\nu}})$, then the 
$(p+1)$-tuple of Jordan normal forms $J_j^{n_{\nu +1}}$ satisfies condition 
$(\alpha _{n_{\nu +1}})$, see Corollary~\ref{alsogood1}, hence, it 
suffices to check that condition $(\alpha _{n_s})$ holds for the 
$(p+1)$-tuple of Jordan normal forms $J_j^{n_s}$. This is true -- 
if $n_s=1$, then $d(J_j^{n_s})=0$ for all $j$ and $(\alpha _{n_s})$ is an 
equality. If there holds 
$(\omega _{n_s})$, then there holds $(\alpha _{n_s})$ as well and it is 
a strict inequality, see Remark~\ref{sameindexofrigidity}.

3) In \cite{Si1} C.Simpson proved in the multiplicative version of the problem 
that for generic eigenvalues and one of the matrices $M_j$ having distinct 
eigenvalues the necessary and sufficient conditions for the solvability of 
the DSP is the 
inequalities $(\alpha _n)$ and $(\beta _n)$ to hold. The author has shown in 
\cite{Ko3} that this is true also if one of the matrices $M_j$ or $A_j$ 
has only eigenvalues of multiplicity $\leq 2$. With  
Theorem~\ref{generic} one gets rid of the condition on the multiplicities of 
the eigenvalues of one of the matrices.

4) The case when condition $(\alpha _n)$ is an equality for matrices $M_j$ 
is considered in 
detail in \cite{Ka} where it is explained how to construct such irreducible 
$(p+1)$-tuples (called {\em rigid}) of matrices $M_j$. Examples of existence 
of rigid $(p+1)$-tuples can be found in \cite{Si1}, \cite{Gl} and \cite{Ko3}. 
\end{rems} 

The theorem does not cover the case of matrices $M_j$, when the PMV of the 
eigenvalues $\sigma _{k,j}$ is non-simple but the eigenvalues are generic. 
In this case the following theorem clarifies partially the situation. 

Denote by $\Sigma _{j,l}(\sigma )$ the number of Jordan blocks of $M_j$ of 
size $l$, with eigenvalue $\sigma$, and by $d$ the greatest common divisor 
of the numbers $\Sigma _{j,l}(\sigma )$ (over all $j$ and $l$, over all 
eigenvalues $\sigma$). 
Even for non-simple PMV one has 

\begin{tm}\label{genericbis}
If $d=1$, then for generic eigenvalues the necessary and sufficient 
condition 
for the solvability of the DSP for matrices $M_j$ with 
given Jordan normal forms $J_j^n$ is the conditions {\em i)} and 
{\em ii)} from Theorem~\ref{generic} to hold.
\end{tm} 

In the case when the eigenvalues are not necessarily generic there holds the 
following  

\begin{tm}\label{nongenericevs}
1) If $d=1$ and if inequality $(\alpha _n)$ is strict, then conditions 
{\em i)} and {\em ii)} from Theorem~\ref{generic} are necessary and 
sufficient for 
the weak solvability of the DSP in the case of matrices $A_j$ for any 
eigenvalues. 

2) If $d=1$ and if inequality $(\alpha _n)$ is strict, then conditions 
{\em i)} 
and {\em ii)} from Theorem~\ref{generic} are necessary and sufficient for 
the solvability of the DSP in the case of matrices $M_j$ with generic 
eigenvalues.
\end{tm}

For matrices $A_j$ part 1) of the theorem is not true if $(\alpha _n)$ is an 
equality. Example: $p=2$, 
$n=2$ and each matrix $A_j$ is nilpotent, of rank 1. Such Jordan normal 
forms satisfy conditions {\em i)} 
and {\em ii)} from Theorem~\ref{generic}, but the triple is (up to 
conjugacy) upper-triangular and its centralizer is generated by $I$ and 
$\left( \begin{array}{cc}0&1\\0&0\end{array}\right)$.

\subsection{Plan of the paper} 

The next three sections introduce the basic ingredients used to prove 
Theorems~\ref{generic}, \ref{genericbis} and \ref{nongenericevs}. 
In Section~\ref{prelim} we describe {\em the basic technical tool} which 
is a way to obtain irreducible 
$(p+1)$-tuples of matrices by deforming $(p+1)$-tuples of matrices with 
trivial centralizers. Such a deformation allows one to keep the Jordan normal 
forms of the $p+1$ matrices the same while changing the eigenvalues. It 
allows also to change their Jordan normal forms to new ones.

In Section \ref{LRNREV} we introduce a result due to A.H.M. Levelt describing 
the structure of the solution to a regular system in a neighbourhood of a 
pole. Lemma~\ref{non-resonant} from that section is important because it is 
used further 
to transform solving the DSP in the multiplicative version into solving it 
in the additive one. This lemma gives also 
a hint why the answers to the DSP in 
both versions are the same in the cases covered by this paper 
(see Corollary~\ref{sameanswer}). 
(We should note that there are cases not covered by the present paper 
in which the formulation of the result in the 
multiplicative version is more complicated than the one in the additive 
version, see Remark~\ref{tricky}.)

Before proving Theorem \ref{generic} we prove its weakened version 
first:

\begin{tm}\label{genericweak}
For given Jordan normal forms whose PMV is simple conditions {\em i)} and 
{\em ii)} from Theorem \ref{generic} are necessary and sufficient for the 
solvability of the DSP (for matrices $A_j$ or $M_j$) for all eigenvalues from 
some Zariski open dense subset in the set 
of all generic eigenvalues with this PMV.
\end{tm}

In Section~\ref{reduction} we explain how to reduce the 
proof of Theorem~\ref{genericweak} to the case of diagonalizable matrices 
$A_j$ or $M_j$. This reduction uses the basic technical tool.

In Section~\ref{diagonal} we formulate the result 
(Theorem~\ref{basicres}) in the case of diagonalizable 
matrices $A_j$ or $M_j$. We deduce Theorem~\ref{genericweak} from 
Theorem~\ref{basicres} at the end of that section.

In Section \ref{sufficiency} we prove the sufficiency and in 
Section \ref{necessity} we prove the necessity of conditions {\em i)} and 
{\em ii)} of Theorem~\ref{generic} for the 
existence of irreducible $(p+1)$-tuples of diagonalizable matrices. This is 
the proof of Theorem~\ref{basicres}. 
In principle, when inequality $(\alpha _n)$ is 
strict, the sufficiency follows from Theorem~\ref{nongenericevs}. We 
prove the sufficiency in Section~\ref{sufficiency} to cover also the case 
when $(\alpha _n)$ is an equality.

In the case when $(\alpha _n)$ is a strict inequality Theorem \ref{generic} 
follows from Theorem \ref{genericweak} and from 
Theorem \ref{nongenericevs}. The latter in the case of generic eigenvalues 
provides the existence of irreducible $(p+1)$-tuples of matrices. 

In the case when $(\alpha _n)$ is an equality 
Theorem \ref{generic} results from  

\begin{tm}\label{genericrigid}
For given Jordan normal forms for which $(\alpha _n)$ is 
an equality, conditions {\em i)} and {\em ii)} from 
Theorem~\ref{generic} are necessary and sufficient for the solvability 
of the DSP (for matrices $A_j$ or $M_j$) for 
any generic eigenvalues.
\end{tm}

The theorem is proved in Section \ref{n_s=1}. In the proof we use 
Theorem~\ref{genericweak}. In fact, the sufficiency and the necessity being 
already proved respectively in Sections~\ref{sufficiency} and \ref{necessity} 
there remains only to be proved in Section~\ref{n_s=1} that if conditions 
{\em i)} and {\em ii)} from 
Theorem~\ref{generic} hold, then for such Jordan normal forms (admitting 
generic eigenvalues) the DSP is 
solvable for {\em all} generic eigenvalues.

Theorem \ref{nongenericevs} is proved in Section~\ref{n_s>1} after some 
preparation, i.e. 
after Section~\ref{adjac} where we discuss adjacency of nilpotent orbits. 

In Section~\ref{proofofgenericbis} we prove Theorem \ref{genericbis}. In the 
proof we use 
Theorem~\ref{genericweak}. 

Section \ref{proofofd_jr_j} contains the proofs of 
Propositions~\ref{d_jr_j} and \ref{necessarynongen}.

\section{The basic technical tool\protect\label{prelim}}

\subsection{The basic technical 
tool in the additive version\protect\label{BTTAV}}

\begin{defi}Call {\em basic technical tool} the procedure described 
below. One starts with a $(p+1)$-tuple of matrices $A_j$ (not necessarily 
irreducible) satisfying 
(\ref{A_j}) and having a {\em trivial centralizer}. Set $A_j=Q_j^{-1}G_jQ_j$, 
$G_j$ being Jordan matrices. One looks for a 
$(p+1)$-tuple of matrices $\tilde{A}_j$ of the form 

\[ \tilde{A}_j=(I+\sum _{i=1}^l\varepsilon _iX_{j,i}(\varepsilon ))^{-1}
Q_j^{-1}
(G_j+\sum _{i=1}^l\varepsilon _iV_{j,i}(\varepsilon ))Q_j
(I+\sum _{i=1}^l\varepsilon _iX_{j,i}(\varepsilon ))\] 
where 
$\varepsilon =(\varepsilon _1,\ldots ,\varepsilon _l)\in ({\bf C}^l,0)$ and 
$V_{j,i}(\varepsilon )$ are given 
matrices analytic 
in $\varepsilon$ (in each concrete application their properties will be 
specified). One has 
tr$(\sum _{j=1}^{p+1}\sum _{i=1}^l\varepsilon _iV_{j,i}(\varepsilon ))
\equiv 0$. Often but not always one chooses $l=1$ and $V_{j,1}$ such that the 
eigenvalues of the $(p+1)$-tuple of matrices $\tilde{A}_j$ are generic for 
$\varepsilon \neq 0$. One looks for $X_j$ analytic in $\varepsilon$ such that 
$\sum _{j=1}^{p+1}\tilde{A}_j\equiv 0$. 
\end{defi}

The condition $\sum _{j=1}^{p+1}\tilde{A}_j=0$ yields (in first 
approximation w.r.t. $\varepsilon$) 

\begin{equation}\label{A_jX_j} 
{\rm for~all~}i=1,\ldots ,l~{\rm one~has}~~~
\sum _{j=1}^{p+1}(Q_j^{-1}V_{j,i}(0)Q_j+[A_j,X_{j,i}(0)])=0
\end{equation}

\begin{prop}\label{[A_j,X_j]}
The centralizer of the $p$-tuple of matrices $A_j$ ($j=1,\ldots ,p$) is 
trivial if and only if the mapping 
$(sl(n,{\bf C}))^p\rightarrow sl(n,{\bf C})$, 
$(X_1,\ldots ,X_p)\mapsto \sum _{j=1}^p[A_j,X_j]$ is surjective.
\end{prop}

{\em Proof:} The mapping is not surjective if and only if the image of each 
mapping $X_j\mapsto [A_j,X_j]$ belongs to one and the same proper linear 
subspace of $sl(n,{\bf C})$. This means that there exists a matrix 
$0\neq D\in sl(n,{\bf C})$ 
such that tr$(D[A_j,X_j])=0$ for all $X_j\in sl(n,{\bf C})$ and for 
$j=1,\ldots ,p$. This is equivalent to tr$([D,A_j]X_j)=0$ for all 
$X_j\in sl(n,{\bf C})$, i.e. $[D,A_j]=0$ for $j=1,\ldots ,p$. 

The proposition is proved. $\hspace{2cm}\Box$

By Proposition \ref{[A_j,X_j]}, equation (\ref{A_jX_j}) is solvable w.r.t. 
$X_{j,i}(0)$. Hence, the equation $\sum _{j=1}^{p+1}\tilde{A}_j=0$ is 
solvable w.r.t. $X_{j,i}$ for $\varepsilon$ small enough by the implicit 
function theorem (we use the surjectivity here). If for $\varepsilon \neq 0$ 
small enough the eigenvalues of the matrices $\tilde{A}_j$ are generic, then 
their $(p+1)$-tuple is irreducible. 

\subsection{The basic technical tool in the multiplicative version
\protect\label{BTTMV}}

We explain here how the {\em basic technical tool} works in the 
multiplicative version. Given a $(p+1)$-tuple of matrices $M_j^1$ with a 
trivial centralizer and satisfying condition (\ref{M_j}), look for 
$M_j$ of the form 

\[ M_j=(I+\sum _{i=1}^l\varepsilon _iX_{j,i}(\varepsilon ))^{-1}(M_j^1+
\sum _{i=1}^l\varepsilon _iN_{j,i}(\varepsilon ))(I+
\sum _{i=1}^l\varepsilon _iX_{j,i}(\varepsilon ))\] 
where the given matrices $N_{j,i}$ depend analytically on 
$\varepsilon \in ({\bf C}^l,0)$ and one looks for matrices 
$X_{j,i}$ analytic in 
$\varepsilon$. (Like in the additive version one can set 
$M_j^1=Q_j^{-1}G_jQ_j$, $N_{j,i}=Q_j^{-1}V_{j,i}Q_j$.)
 
The matrices $M_j$ must satisfy 
equality (\ref{M_j}). In first approximation w.r.t. $\varepsilon$ this 
implies that 

\[ {\rm for~all~}i=1,\ldots ,l~{\rm one~has}~~~
\sum _{j=1}^{p+1}M_1^1\ldots M_{j-1}^1([M_j^1,X_{j,i}(0)]+
N_{j,i}(0))
M_{j+1}^1\ldots M_{p+1}^1=0\]
or

\begin{equation}\label{P_j} 
{\rm for~all~}i=1,\ldots ,l~{\rm one~has}~~~
\sum _{j=1}^{p+1}P_{j-1}([M_j^1,X_{j,i}(0)(M_j^1)^{-1}]+
N_{j,i}(0)(M_j^1)^{-1})
P_{j-1}^{-1}=0
\end{equation}
with $P_j=M_1^1\ldots M_j^1$, $P_{-1}=I$ (recall that there holds 
(\ref{M_j}), therefore $M_j^1M_{j+1}^1\ldots M_{p+1}^1=P_{j-1}^{-1}$).

Equality (\ref{M_j}) implies that $\det M_1\ldots \det M_{p+1}=1$. 
One has $\det M_j=\det M_j^1\det (I+\sum _{i=1}^l\varepsilon _i(M_j^1)^{-1}
N_{j,i})$
$=(\det M_j^1)(1+\sum _{i=1}^l\varepsilon _i$tr$((M_j^1)^{-1}N_{j,i}(0))+$
terms of order $\geq 2$ in $\varepsilon$. As 
$\det M_1^1\ldots \det M_{p+1}^1=1$, one has for all $i$ 
tr$(\sum _{j=1}^{p+1}(M_j^1)^{-1}N_{j,i}(0))$$=0$ (terms of first order w.r.t. 
$\varepsilon$ in $\det M_1\ldots \det M_{p+1}$). 

Equation (\ref{P_j}) can be written in the form 

\begin{equation}\label{S_j} 
{\rm for~all~}i=1,\ldots ,l~{\rm one~has}~~~
\sum _{j=1}^{p+1}([S_j,Z_{j,i}]+T_{j,i})=0 
\end{equation}
with $S_j=P_{j-1}M_j^1P_{j-1}^{-1}$, 
$Z_{j,i}=P_{j-1}X_{j,i}(0)(M_j^1)^{-1}P_{j-1}^{-1}$, 
$T_{j,i}=P_{j-1}N_{j,i}(0)(M_j^1)^{-1}P_{j-1}^{-1}$. The centralizers of the 
$(p+1)$-tuples of matrices $M_j^1$ and $S_j$ are the same (to be checked 
directly), i.e. they are both trivial. Hence, for all $i$ the mappings  

\[ (sl(n,{\bf C}))^{p+1}\rightarrow sl(n,{\bf C})~,~
(Z_{1,i},\ldots ,Z_{p+1,i})\mapsto \sum _{j=1}^{p+1}[S_j,Z_{j,i}]\]
are surjective (Proposition \ref{[A_j,X_j]}). Recall that for all $i$ one has  
tr$(\sum _{j=1}^{p+1}(M_j^1)^{-1}N_{j,i}(0))=0$, i.e. 
tr$(\sum _{j=1}^{p+1}T_{j,i})=0$. Hence, equation (\ref{S_j}) can be solved 
w.r.t. the unknown matrices $Z_{j,i}$ and, hence, equation (\ref{P_j}) can be 
solved w.r.t. the matrices $X_{j,i}(0)$. By the implicit function 
theorem (we use the surjectivity here), one can find $X_{j,i}$ analytic in 
$\varepsilon \in ({\bf C}^l,0)$, i.e. one can find the necessary matrices 
$M_j$. 

A first application of the basic technical tool is the following 

\begin{lm}\label{allgen}
For a given $(p+1)$-tuple of Jordan normal forms admitting generic 
eigenvalues denote by $L$ (resp. by $L'$) the set of all 
possible eigenvalues (resp. of all possible generic eigenvalues). If there 
exists a $(p+1)$-tuple 
of matrices $A_j$ (or $M_j$) with a trivial centralizer for some 
$\lambda _0\in L$, then there exist $(p+1)$-tuples 
of matrices $A_j$ (or $M_j$) with trivial centralizers for all 
$\lambda \in L$ sufficiently close to $\lambda _0$ and for all $\lambda$ 
from some Zariski open and dense subset of the connected component $L'_0$ 
of $L'$ containing $\lambda _0$. 
\end{lm}

\begin{rem}
In the multiplicative version of the DSP the set $L$ can consist of several 
(namely, $q$) connected components if 
$q>1$, see Remark~\ref{tricky}. For a given component the product of the 
eigenvalues $\sigma _{k,j}$ with multiplicities divided by $q$ equals one and 
the same root of unity of $q$-th order.
\end{rem}

{\em Proof of the lemma:}

If there exists a $(p+1)$-tuple 
of matrices $A_j$ (or $M_j$) with a trivial centralizer for some 
$\lambda _0\in L$, then there exist $(p+1)$-tuples 
of matrices $A_j$ (or $M_j$) with trivial centralizers for all 
$\lambda \in L$ sufficiently close to 
$\lambda _0$ (it suffices to apply the basic 
technical tool with diagonal 
matrices $V_j$ which are polynomials of the semi-simple parts of the 
matrices $G_j$). The set of 
all such values $\lambda$ is constructible. Hence, its intersection with 
$L'_0$ contains a Zariski open dense subset of $L'_0$. 

The lemma is proved.$\hspace{2cm}\Box$

\section{Levelt's result and non-resonant eigenvalues\protect\label{LRNREV}}
\subsection{Levelt's result\protect\label{Leveltsresult}}

In \cite{L} Levelt gives the structure of the solution to a regular 
system at a pole:

\begin{tm}\label{Lv}
In a neighbourhood of a pole the solution to the  
regular linear system 

\begin{equation}\label{bs}
\dot{X}=A(t)X 
\end{equation}
is representable in the form
\begin{equation}
\label{L*}
 X = U_j(t - a_j) (t-a_j)^{D_j} (t - a_j)^{E_j} G_j 
\end{equation}
where $U_j$ is holomorphic in a neighbourhood of the pole $a_j $, with 
$\det U_j\neq 0$ for $t\neq a_j$; 
$D_j = \mbox{diag} \, ( \varphi_{1,j} ,\ldots , \varphi_{n,j} )$,
$\varphi_{n,j} \in {\bf Z}$; $G_j\in GL(n,{\bf C})$. The matrix $E_j$ is in  
upper-triangular form and the real parts of its eigenvalues
belong to $[0,1)$ (by definition, $(t-a_j)^{E_j} = e^{E_j \ln
(t-a_j)}$). The numbers $\varphi_{k,j}$  satisfy the 
condition (\ref{L***})  formulated below.

System (\ref{bs}) is fuchsian at $a_j$ if and only if
\begin{equation}
\label{L**}
\det U_j (0) \ \neq 0
\end{equation}
We formulate the condition on $\varphi_{k,j} .$ Let $E_j$ have one and the 
same 
eigenvalue in the rows with indices $s_1<s_2<$$\ldots$$<s_q$. Then we 
have
\begin{equation}
\label{L***}
\varphi_{s_1,j} \geq \varphi_{s_2,j} \geq \ldots \geq \varphi_{s_q,j}
\end{equation}
\end{tm}

\begin{rems}\label{remLevelt}
1) Denote by $\beta _{k,j}$ the diagonal entries (i.e. the 
eigenvalues) of the matrix $E_j$. If the system is fuchsian, then the sums 
$\beta _{k,j}+\varphi _{k,j}$ are the eigenvalues $\lambda _{k,j}$ 
of the matrix-residuum $A_j$, see \cite{Bo1}, Corollary 2.1.

2) The numbers $\varphi _{k,j}$ are defined as valuations in the solution 
eigensubspace for the eigenvalue $\exp (2\pi i\beta _{k,j})$ of the 
monodromy operator, see the details in \cite{L}.

3) One can assume without loss of generality that equal eigenvalues of 
$E_j$ occupy consecutive positions on the diagonal and that the matrix 
$E_j$ is block-diagonal, with diagonal blocks of sizes equal to their 
multiplicities. The diagonal blocks themselves are upper-triangular.

4) The following proposition will be used several times in the proofs 
and is of independent interest. It deals with the case when the 
monodromy group of a fuchsian system is reducible, i.e. there is a 
proper subspace invariant for all monodromy operators.
\end{rems}

\begin{prop}\label{bolibr}
The sum of the eigenvalues $\lambda _{k,j}$ of the matrices-residua $A_j$ of 
system (\ref{Fuchs}) corresponding to an invariant 
subspace of the solution space is a non-positive integer.
\end{prop}

The proposition is proved in \cite{Bo1}, see Lemma 3.6 there.

\subsection{Non-resonant eigenvalues\protect\label{somelemmas}}

\begin{defi}
Define as {\em non-resonant} the 
eigenvalues $\lambda _{k,j}$ of the conjugacy class $c_j\in gl(n,{\bf C})$ 
if there are no non-zero integer differences between them. 
\end{defi} 

\begin{rem}
Let the PMV of the eigenvalues $\sigma _{k,j}$ of the 
monodromy operators $M_j$ of system (\ref{Fuchs}) be non-simple; 
denote by $d^*$ the greatest common divisor of its components. Set 
$\lambda _{k,j}=\beta _{k,j}+\varphi _{k,j}$ for the 
eigenvalues of the matrices-residua $A_j$ where Re$\beta _{k,j}\in [0,1)$ 
and $\varphi _{k,j}\in {\bf Z}$, see Subsection~\ref{Leveltsresult}. These 
conditions define unique numbers 
$\beta _{k,j}$. If $d^*$ does not divide the sum 
$\sum _{j=1}^{p+1}\sum _{k=1}^n\beta _{k,j}$ (this sum is always integer 
because $\varphi _{k,j}\in {\bf Z}$ and there holds (\ref{sumEVszero})), 
then the monodromy group cannot 
be realized by a fuchsian system with non-resonant eigenvalues because for 
such eigenvalues equality (\ref{sumEVszero}) would not hold.
\end{rem}

\begin{lm}\label{sameJNF}
For non-resonant eigenvalues $\lambda _{k,j}$ of the matrix-residuum $A_j$ 
of system (\ref{Fuchs}) the Jordan normal forms of the 
matrix $A_j$ and of the monodromy operator $M_j$ are the same and $M_j$ 
is conjugate to $\exp (2\pi iA_j)$. 
\end{lm}

{\em Proof:}

Use Levelt's form (\ref{L*}) of the solution to system (\ref{bs}) (presumed 
to be fuchsian at $a_j$) and 1) and 
3) from Remarks~\ref{remLevelt}. One has $A(t)=\dot{X}X^{-1}$. Hence, if the 
eigenvalues of $A_j$ are non-resonant, then to equal eigenvalues of $E_j$ 
there 
correspond equal eigenvalues of $D_j$, the matrices $D_j$ and $E_j$ commute 
and $A_j=U_j(0)(D_j+E_j)(U_j(0))^{-1}$ (to be checked directly). One has 
$M_j=G_j^{-1}\exp (2\pi iE_j)G_j$. Hence, for the Jordan normal form 
$J(M_j)$ of $M_j$ one has $J(M_j)=J(E_j)=J(D_j+E_j)=J(A_j)$. 

The lemma is proved. $\hspace{2cm}\Box$

\begin{lm}\label{non-resonant}
Every irreducible monodromy group with a simple PMV of the eigenvalues 
$\sigma _{k,j}$ can be realized by a fuchsian system with non-resonant 
eigenvalues $\lambda _{k,j}$. 
\end{lm}

The lemma follows directly from Lemma 10 from \cite{Ko4}. 

\begin{cor}\label{sameanswer}
For a given $(p+1)$-tuple of Jordan normal forms with a simple PMV the DSP 
is solvable for some generic eigenvalues for matrices 
$A_j$ if and only if it is solvable 
for some generic eigenvalues for matrices $M_j$.
\end{cor}

{\em Proof:} 

If for the given 
$(p+1)$-tuple of Jordan normal forms the DSP is solvable 
for some generic eigenvalues for matrices 
$A_j$, then 
one can choose such a $(p+1)$-tuple with not only generic, but with strongly 
generic non-resonant eigenvalues (this can be achieved by 
multiplying the given $(p+1)$-tuple of matrices 
$A_j$ by some constant $c\in {\bf C}^*$), and then use 
Lemma~\ref{sameJNF} 
(if the PMV of the eigenvalues $\lambda _{k,j}$ is simple, then  
non-resonant eigenvalues $\lambda _{k,j}$ exist). 

If for the given 
$(p+1)$-tuple of Jordan normal forms with a simple PMV the DSP is solvable  
for some generic eigenvalues for matrices 
$M_j$, then it is possible to realize such a  
monodromy group by a fuchsian system with non-resonant eigenvalues 
$\lambda _{k,j}$, see Lemma~\ref{non-resonant}. Hence, for all 
$j$ one would have $J(A_j)=J(M_j)$ (Lemma~\ref{sameJNF}).$\hspace{2cm}\Box$

\section{How to reduce the problem to the case of diagonalizable matrices ?
\protect\label{reduction}}

\subsection{Correspondence between Jordan normal forms \protect\label{diag}} 

All Jordan matrices in this subsection are presumed upper-triangular.

\begin{defi}\label{corJordforms}
Let a non-diagonal Jordan normal form $J_0=\{ b_{i,l}\}$ of size $n$ be 
given. We define its {\em associated} semi-simple Jordan normal form $J_1$ 
(also of size $n$) such 
that the quantities $r_j:=r(J_j)$ and 
$d_j:=d(J_j)$ are the same for $j=0$ and $j=1$. 

A semi-simple Jordan normal form is the same as a partition of $n$, the 
parts being the multiplicities of the eigenvalues. If $J_0=\{ b_{i,l}\}$ 
($i\in I_l$, $l\in L$), one views for each $l$ the set $\{ b_{i,l}\}$ as 
a partition of $\sum _ib_{i,l}$ and one takes for $J_1$ the disjoint sum 
of the dual partitions.

We will also say that the Jordan normal form $J_1$ {\em corresponds} 
to $J_0$ and that $J_0$ {\em corresponds} to $J_1$. Any Jordan normal form 
$J$ corresponding to $J_1$ corresponds to $J_0$ and $J_0$ corresponds to $J$.
\end{defi}

\begin{defi}\label{corJordforms1}
Denote by $G^0$, $G^1$ two Jordan matrices with Jordan normal 
forms $J_0$, $J_1$ corresponding to each other, where 
$G^1$ is diagonal and $G^0$ is block-diagonal, each diagonal block having 
a single eigenvalue, different blocks having different eigenvalues. In the 
block-decomposition defined by the multiplicities of the eigenvalues of $G^0$ 
the matrix $G^1$ has diagonal blocks with mutually different eigenvalues. 
For each diagonal block the eigenvalues of $G^1$ occupying the last but $q$ 
positions of the Jordan blocks of $G^0$ are equal (we denote them by $h_q$); 
for $q_1\neq q_2$ one has $h_{q_1}\neq h_{q_2}$.
\end{defi}

\begin{ex}
Let $J_0=\{ \{ 4,3,2\} \{ 3,1\} \}$, i.e. there are two eigenvalues to the 
first (resp. the second) of which there correspond Jordan blocks of sizes 
4, 3, 2 (resp. 3, 1). Hence, $J_1$ is defined by the MV (3,3,2,2,1,1,1). 
Indeed, the partition of $9=4+3+2$ dual to 4,3,2 is 3,3,2,1, the partition of 
$4=3+1$ dual to 3,1 is 2,1,1. When taking the direct sum of these dual 
partitions one rearranges the components of the MV so that they form a 
non-increasing sequence.
\end{ex}

\begin{rems}\label{corrr}
1) If $J_0$ is the Jordan normal form of $A$, then the multiplicities of the 
eigenvalues for $J_1$ are the numbers dim Ker$((A-\lambda I)^{k+1}-$dim 
Ker$((A-\lambda I)^k$ which are non-zero. 

2) One can show that any generic deformation of a matrix with 
Jordan normal form $J_0$ contains matrices with Jordan normal form $J_1$ and 
that these are the diagonalizable matrices from orbits of {\em least} 
dimension encountered in the deformation.

3)If a Jordan normal form is a direct sum of two Jordan normal forms with 
no eigenvalue in common, i.e. $J=J^*\oplus J^{**}$, and if the Jordan normal 
forms ${J^*}'$, ${J^*}''$ with no eigenvalue in common correspond to 
$J^*$, $J^{**}$, then ${J^*}'\oplus {J^*}''$ corresponds to $J$.
\end{rems}

\begin{prop}\label{samerank}
The quantities $r(J^{(i)})$ computed for two Jordan normal forms $J'$ and 
$J''$ corresponding to one another coincide.
\end{prop}

{\em Proof:} 

It suffices to prove the proposition for $J'=J_0$, $J''=J_1$ (see 
Definition~\ref{corJordforms}), i.e. to 
prove that $r_0=r_1$. Let $m^0$ be the greatest number of Jordan blocks of 
$J_0$ with a given eigenvalue. Hence, $r_0=n-m^0$. The 
construction of $J_1$ implies that the greatest of the multiplicities of the 
eigenvalues of $J_1$ equals $m^0$ -- this follows from the definition of a 
dual partition. Thus, $r_1=r_0$.$\hspace{2cm}\Box$

\begin{prop}\label{samedimension}
The dimensions of the orbits of two matrices with Jordan normal forms 
corresponding to one another are the same.
\end{prop}

{\em Proof:} 

$1^0$. It suffices to prove the proposition in the case when one of the 
Jordan normal forms is diagonal. Denote the two matrices by $G^0$ and $G^1$ 
where $G^i$ are defined by Definition~\ref{corJordforms1}. 
The dimension of the orbit of $G^j$, $j=0,1$, equals 
$n^2-{\rm dim}Z(G^j)$ where $Z(G^j)$ is the centralizer of $G^j$ in 
$gl(n,{\bf C})$. Block-decompose the matrices from $gl(n,{\bf C})$ 
with sizes of the diagonal blocks equal to the multiplicities of the 
eigenvalues of $G^0$. Then the off-diagonal blocks of 
$Z(G^j)$ are $0$; indeed, two diagonal blocks of $G^j$ ($j=1,2$) 
have no eigenvalue in common. This observation 
allows when computing the dimensions of the orbits to consider only the 
case when $J_0$ has only one eigenvalue.

$2^0$. Show that in this case one has  
${\rm dim}\, Z(G^0)={\rm dim}\, Z(G^1)$, hence, the dimensions of the orbits 
of $G^0$ and $G^1$ are the same. One has

\[ {\rm dim}\, Z(G^0)=b_1+3b_2+5b_3+\ldots +(2r-1)b_r \] 
where $b_1\geq b_2\geq b_3\geq \ldots \geq b_r$ are the sizes of the 
Jordan blocks of $J_0$, see~\cite{Ar}, p. 229; 

\[ {\rm dim}\, Z(G^1)=(k_1)^2+(k_2)^2+\ldots +(k_{n_1})^2 \] 
where $k_i$ are the multiplicities of the 
eigenvalues of $J_1$. 

$3^0$. The first $b_r$ of the numbers $k_j$ equal $r$, the next 
$(b_{r-1}-b_r)$ equal $r-1$, the next $(b_{r-2}-b_{r-1})$ equal $r-2$ etc. 
Thus 

\[ (k_1)^2+\ldots +(k_{n_1})^2=b_rr^2+(b_{r-1}-b_r)(r-1)^2+\ldots +
(b_1-b_2)\times 1^2= \] 

\[ =b_r[r^2-(r-1)^2]+b_{r-1}[(r-1)^2-(r-2)^2]+\ldots +b_1\times 1= \] 

\[ =(2r-1)b_r+(2r-3)b_{r-1}+\ldots +b_1={\rm dim}\, Z(G^0) \] 

The proposition is proved. $\hspace{2cm}\Box$

\begin{prop}\label{simultaneously}
Let the two Jordan normal forms ${J^n}'$ and 
${J^n}''$ correspond to one another. Choose in each of them an eigenvalue 
with maximal number of Jordan blocks. By Proposition~\ref{samerank} these 
numbers coincide. Denote them by $k'$. 
Decrease by 1 the sizes of the $k$ smallest Jordan blocks with these 
eigenvalues where $k\leq k'$. Then the two Jordan normal 
forms of size $n-k$ obtained in this way correspond to one another. 
\end{prop}

\begin{rem}\label{severalchoices}
The proposition implies in particular that if a Jordan normal 
form has several (say, $s$) eigenvalues 
with the maximal number $k'$ of Jordan blocks and if one constructs 
$s$ new Jordan 
normal forms by decreasing by 1 the sizes of the $k$ smallest blocks with a 
given one of these eigenvalues, then these $s$ Jordan normal forms 
correspond to one another.
\end{rem}

{\em Proof:}

$1^0$. It suffices to prove the proposition in the case when the Jordan 
normal form ${J^n}''$ is diagonal (if this is not so, then consider together 
with ${J^n}'$ and ${J^n}''$ the diagonal Jordan normal form ${J^n}'''$ 
corresponding to them, then prove the proposition for the couples 
${J^n}'$,${J^n}'''$ and ${J^n}''$,${J^n}'''$). 
In the case when ${J^n}''$ is 
diagonal one simply decreases by $k$ the biggest component of the MV. 

$2^0$. Assume that the $\nu$-th eigenvalue of ${J^n}'$ has 
$k'$ Jordan blocks, of sizes $b_{i,\nu }$, $i=1,\ldots ,k'$. 
Decreasing by 1 the least $k$ of the integers $b_{i,\nu }$ (considered as 
parts of 
the partition of $\sum _{i=1}^{k'}b_{i,\nu}$) results in decreasing by 
$k$ the biggest part of its 
dual partition. By definition, this biggest part equals $k'$ and it 
is (one of) the biggest component(s) of 
the MV defining ${J^n}''$.  

The proposition is proved. $\hspace{2cm}\Box$

\begin{prop}\label{singleEV}
For each diagonal Jordan normal form $J_1$ there exists a unique 
Jordan normal form $J_0$ with a single eigenvalue which corresponds 
to $J_1$. Hence, the same is true for any Jordan normal form. 
\end{prop}

{\em Proof:}

It follows from the construction of $J_1$ after $J_0$ that if the 
multiplicities of the eigenvalues of $J_1$ equal $g_1\geq \ldots \geq g_d$, 
$g_1+\ldots +g_d=n$, then $J_0$ has exactly $g_{\nu}$ Jordan blocks of 
size $\geq \nu$. This condition defines a unique Jordan normal form $J_0$ 
with a single eigenvalue. $\hspace{2cm}\Box$

Recall that the matrices $G^0$ and $G^1$ were defined in 
Definition \ref{corJordforms}.

\begin{prop}\label{sameorbit}
1) If $G^0$ is nilpotent, then the 
orbits of the matrices $\varepsilon G^1$ and 
$G^0+\varepsilon G^1$ are the same for $\varepsilon \in {\bf C}^*$.

2) If $G^0$ is not necessarily nilpotent, then the matrix 
$G^0+\varepsilon G^1$ is diagonalizable and for 
$\varepsilon \in {\bf C}^*$ small enough its Jordan normal form is $J_1$, 
its orbit is the one of $G^0_s+\varepsilon G^1$ where $G^0_s$ is the 
semisimple part of $G^0$. 
\end{prop} 

{\em Proof:}

$1^0$. Let $G^0$ be nilpotent (hence, there is just one diagonal block of 
size $n$). Conjugate the matrices $\varepsilon G^1$ and 
$G^0+\varepsilon G^1$ with a permutation matrix $Q$ 
such that after the permutation the eigenvalues $h_0$ occupy the last 
positions on the diagonal preceded by the eigenvalues $h_1$ preceded by the 
eigenvalues $h_2$ etc. 

$2^0$. If one block-decomposes a matrix with sizes of the 
diagonal blocks equal to the multiplicities of the eigenvalues $h_q$, then 
the units of the matrix $Q^{-1}G^0Q$ will be all in the blocks above the 
diagonal. Hence, the matrix $G^*=Q^{-1}(G^0+\varepsilon G^1)Q$ in this block 
decomposition is block upper-triangular and has scalar diagonal blocks with 
mutually distinct eigenvalues. 
Hence, one can conjugate this matrix with a block upper-triangular matrix 
and after the conjugation the units above the diagonal disappear and the 
resulting matrix is diagonal, with the same diagonal blocks as $G^*$. Hence, 
this is the matrix $\varepsilon G^1$.

$3^0$. If $G^0$ has one eigenvalue (not necessarily equal to 0), then the 
second statement of the proposition follows from the first one.

$4^0$. If $G^0$ is arbitrary, then one can block decompose it, the diagonal 
blocks having each one eigenvalue, the eigenvalues of different diagonal 
blocks being different, and then apply the result from $3^0$ to 
every diagonal block. (For small values of $\varepsilon \in {\bf C}^*$ two 
different diagonal blocks will have no eigenvalue in common.)

The proposition is proved. $\hspace{2cm}\Box$

\begin{prop}\label{JJ'}
Denote by $J=\{ b_{i,l}\}$ an arbitrary Jordan normal form of size $n$ and by $J'$ its 
corresponding 
Jordan normal form with a single eigenvalue. Recall that 
for each fixed $l$ one has $b_{1,l}\geq b_{2,l}\geq \ldots \geq b_{s_l,l}$. 
Then the size of the $k$-th Jordan block of $J'$ (in decreasing order) 
equals $\sum _lb_{k,l}$ (if some of the numbers participating in this sum 
are not defined, then they are presumed to equal 0).
\end{prop}

{\em Proof:}

$1^0$. It suffices to consider the case of two eigenvalues. The general case 
can be treated by induction on the number of eigenvalues (one represents a 
Jordan normal form $J_1^n$ with $k$ eigenvalues as direct sum of a Jordan 
normal form $J_2^m$ with a 
single eigenvalue and a Jordan normal form $J_3^{n-m}$ with 
$k-1$ eigenvalues; then one finds the Jordan normal form $J_4^{n-m}$ with a 
single 
eigenvalue corresponding to $J_3^{n-m}$ and finally the Jordan normal form 
$J_5^n$ with a 
single eigenvalue corresponding to $J_2^m\oplus J_4^{n-m}$; $J_5^n$ 
corresponds to $J_1^n$ and $J_2^m\oplus J_4^{n-m}$ has two eigenvalues). 

In the case of two eigenvalues 
it suffices to show that if the sizes of the blocks of $J'$ are as in the 
proposition, then to $J$ and $J'$ there corresponds one and the same 
diagonal Jordan normal form $J_d=J'_d$. 

$2^0$. Denote the two eigenvalues of $J$ by $\lambda$ and $\sigma$. 
Assume that to $\lambda$ there correspond no less Jordan blocks than to 
$\sigma$, i.e. $s_1\geq s_2$. This means that the greatest of the 
multiplicities of eigenvalues both of $J_d$ and of $J'_d$ equals $s_1$. 

Consider the last $b_{s_1,1}$ positions of every Jordan block with 
eigenvalue $\lambda$. They give rise to $b_{s_1,1}$ eigenvalues each of 
multiplicity $s_1$ in $J_d$. 

Decrease 

1) the size $n$ of the 
matrices by $s_1b_{s_1,1}$, 

2) the sizes of each of the Jordan blocks of $J$ with eigenvalue $\lambda$ 
by $b_{s_1,1}$ and 

3) the sizes of all Jordan blocks of $J'$ by $b_{s_1,1}$.

Hence, in each of the diagonal Jordan normal forms $J_d$ and $J'_d$ one loses 
$b_{s_1,1}$ eigenvalues each of 
multiplicity $s_1$. Hence, the Jordan normal forms $J_d$ and $J'_d$ coincide or 
not simultaneously before and after the reduction by $s_1b_{s_1,1}$ of 
the sizes of the matrices.

$3^0$. After a finite number of such reductions one of the two eigenvalues 
becomes of multiplicity 0; in this case there is nothing to prove.

The proposition is proved.$\hspace{2cm}\Box$\\

Denote by $J$ and $J'$ an arbitrary Jordan normal form and its corresponding 
Jordan normal form with a 
single eigenvalue. Consider a couple $\Delta$, $D'$ of Jordan matrices with 
these Jordan normal forms $D'$ being nilpotent. Suppose that the Jordan 
blocks of $\Delta$ of sizes $b_{k,l}$ for $k$ 
fixed are situated in the same rows where is situated the Jordan block 
of size $\sum _lb_{k,l}$ of $D'$ (see the previous proposition). Denote by 
$\Delta _s$ the diagonal (i.e. semi-simple) part of the matrix $\Delta$.

\begin{prop}\label{DD'} 
For all $\varepsilon \neq 0$ the matrix $\varepsilon \Delta _s+D'$ is 
conjugate to $\varepsilon \Delta$.
\end{prop}

It suffices to prove the proposition in the case when $D'$ 
has a single Jordan block of size $n$. In this case one checks directly that 
for all eigenvalues $a$ of $\Delta$ one has 
rk$(\varepsilon (\Delta _s-aI)+D')=$rk$(\varepsilon (\Delta -aI))=n-1$. 
For all other values of $a\in {\bf C}$ these ranks 
equal $n$.$\hspace{2cm}\Box$ 

\begin{rem}
Permute the diagonal entries of $\Delta _s$ so that before and after 
the permutation each entry remains in one of the rows of one and the same 
Jordan block of $D'$. Then the proposition holds again and the proof is the 
same. 
\end{rem}

\subsection{Subordinate conjugacy classes and normalized chains of eigenvalues
\protect\label{NC}}

\begin{defi}
Given two conjugacy classes $c'$, $c''$ with one and the 
same eigenvalues, of one and the same multiplicities, we say that $c''$ is 
{\em subordinate} to $c'$ if $c''$ lies in the closure of $c'$, i.e. for any 
matrix $A\in c''$ there exists a deformation $\tilde{A}(\varepsilon )$, 
$\tilde{A}(0)=A$ such that for $0\neq \varepsilon \in ({\bf C},0)$ one has 
$\tilde{A}(\varepsilon )\in c'$.
\end{defi} 

\begin{ex}
Let $n=4$ and let the eigenvalues be $a$, $a$, $b$, $b$, 
$a\neq b$. Let 
$c'$ (resp. $c''$) have one Jordan block $2\times 2$ (resp. two Jordan blocks 
$1\times 1$) with eigenvalue $a$ and both $c'$ and $c''$ have two Jordan 
blocks $1\times 1$ with eigenvalue $b$. Then $c''$ is subordinate to $c'$. 
If the conjugacy class $c'''$ has the same eigenvalues, two Jordan blocks 
$1\times 1$ with eigenvalue $a$ and one Jordan block $2\times 2$ with 
eigenvalue $b$, then neither $c'$ is subordinate to $c'''$ nor $c'''$ is 
subordinate to $c'$ (and $c''$ is subordinate to $c'''$).
\end{ex}

Notice that in the above example the Jordan normal forms of $c'$ and $c'''$ 
are the same.

\begin{defi}
Given two Jordan normal forms $J'$, $J''$, we say that 
$J''$ is {\em subordinate} to $J'$ if there exist conjugacy classes 
$c'$, $c''$ defining the Jordan normal forms $J'$, $J''$ such that $c''$ is 
subordinate to $c'$. 
\end{defi}

\begin{defi}
1) For a diagonalizable matrix represent the set of its eigenvalues as a union 
of maximal non-intersecting subsets of eigenvalues congruent modulo 
${\bf Z}$ (called further ${\bf Z}$-{\em subsets}). For each ${\bf Z}$-subset 
define its multiplicity vector where the different 
eigenvalues of the ${\bf Z}$-subset are ordered so that their real parts 
form a decreasing sequence. 
Then the eigenvalues of the matrix are said to form a {\em normalized chain} 
if for every such multiplicity vector its components form a non-decreasing 
sequence.

2) If the eigenvalues of the diagonalizable matrix $A$ form a 
normalized chain, then the multiplicity vector of each ${\bf Z}$-subset 
defines a 
diagonal Jordan normal form $J_i$. Denote by $J'_i$ its corresponding Jordan 
normal form with a single eigenvalue. Denote by $\tilde{J}(A)$ the 
Jordan normal form $\oplus _iJ'_i$ where the sum is taken over all 
${\bf Z}$-subsets; for $i_1\neq i_2$ the eigenvalues of $J'_{i_1}$ and 
$J'_{i_2}$ are different. Hence, $\tilde{J}(A)$ 
corresponds to $J(A)$, see 3) of Remarks~\ref{corrr}. If $J(A)=J^0$, 
then we set $\tilde{J}(J^0)=\tilde{J}(A)$.
\end{defi}

\begin{lm}\label{JNFofM} 
Let in system (\ref{Fuchs}) the matrix $A_1$ be  
with Jordan normal form $J_1$ and let its eigenvalues form a normalized 
chain. Then the Jordan 
normal form of the monodromy operator $M_1$ is either $J^*:=\tilde{J}(J_1)$ 
or is one subordinate to it.
\end{lm}  

{\em Proof:}

$1^0$. Consider first the case when $A_1$ is diagonal and 
$J^*$ has just one eigenvalue. Use Theorem~\ref{Lv}. If the solution to 
system (\ref{Fuchs}) is represented in form (\ref{L*}), with 
$\varphi _{1,1}\geq \ldots \geq \varphi _{n,1}$, then one has 
$E_1=\alpha I+F$ (Re$(\alpha )\in [0,1)$) where the matrix $F$ is nilpotent 
and upper-triangular. 

$2^0$. More exactly, $F$ is block upper-triangular, with zero 
diagonal blocks; the diagonal blocks are of sizes equal to the 
multiplicities of the eigenvalues of the matrix $D_1$ from (\ref{L*}). 
Indeed, the presence 
of non-zero entries in the diagonal blocks of $F$ would result in $A_1$ not 
being diagonalizable (we propose to the reader to check this oneself).

$3^0$. Denote the MV of the eigenvalues of $D_1$ (it is also the one of 
$J_1$) by 
$(l_d,\ldots ,l_1)$ where $l_d\leq \ldots \leq l_1$ (these inequalities 
follow from the definition of $J_1$ in the previous subsection -- for each 
$q$ the number $l_q$ of eigenvalues $h_q$ equals the number of Jordan blocks 
of $J^*$ of size $\geq q+1$). 

$4^0$. The rank of the matrix $(F)^{\nu}$ cannot exceed 
$\tilde{l}_{\nu}:=l_d+\ldots +l_{\nu +1}$ 
(only the first $\tilde{l}_{\nu}$ rows of $(F)^{\nu}$ can be non-zero). 
This is exactly the rank of $(N)^{\nu}$, $N$ being a nilpotent matrix with 
Jordan normal form $J^*$. Hence, the Jordan normal form of $F$ is 
either $J_0$ or is one subordinate to it. Indeed, the inequalities 
rk$(F)^{\nu}\leq$rk$(N)^{\nu}$, $\nu =1,2,\ldots$ imply that either the 
orbits of 
$F$ and $N$ coincide (if there are equalities everywhere) or that the orbit 
of $F$ lies in the closure of the one of $N$ (if at least one inequality is 
strict), see \cite{Kr}, p. 21. 

$5^0$. On the other hand, one has (up to conjugacy) 
$M_1=\exp (2\pi iE_1)$$=\exp (2\pi i\alpha )\exp (2\pi iF)$. This means 
that the Jordan normal forms of $M_1$ and $F$ coincide. Hence, the 
Jordan normal form of $M_1$ is either $J^*$ or is one subordinate to it.

$6^0$. In the general case (when $J^*$ has several eigenvalues) one uses 
3) of Remarks~\ref{remLevelt} and applies the above reasoning to each 
diagonal block of $E_j$, i.e. to each eigenvalue of the monodromy 
operator $M_j$. 

The lemma is proved.$\hspace{2cm}\Box$
 
\subsection{Reduction to the case of diagonalizable matrices $A_j$}

Denote by $J_j^0$ the Jordan normal forms of the matrices $A_j$ or $M_j$. 
Denote by $J_j^1$ their corresponding diagonal Jordan normal forms defined in 
Subsection~\ref{diag}. 

\begin{lm}\label{J^0impliesJ^1}
1) The DSP is solvable for Jordan normal forms $J_j^0$ with a simple PMV and 
for some generic eigenvalues 
if and only if it is solvable for the Jordan normal forms $J_j^1$ and 
for some generic eigenvalues.

2) If for some generic eigenvalues and given Jordan normal forms $J_j$ with 
a simple PMV the DSP is solvable, then it is solvable  
(for some generic eigenvalues) for all $(p+1)$-tuples of 
Jordan normal forms $J_j'$ where for each $j$ either $J_j$ is subordinate to 
$J_j'$ or $J_j=J_j'$.

The lemma holds for matrices $A_j$ and for matrices $M_j$.
\end{lm}

{\em Proof:} 

$1^0$. Prove the lemma first for matrices $A_j$. Denote by $G_j^i$, $i=0,1$,  
two Jordan matrices defining the same Jordan normal forms as $J_j^i$ and 
such that $A_j=Q_j^{-1}G_j^0Q_j$; we define the matrices $G_j^i$ like the 
matrices $G^i$ from Definition~\ref{corJordforms1}. The existence of 
irreducible $(p+1)$-tuples of matrices 

\[ \tilde{A}_j=(I+\varepsilon X_j(\varepsilon ))^{-1}Q_j^{-1}
(G_j^0+\varepsilon G_j^1)Q_j(I+\varepsilon X_j(\varepsilon )),\] 
satisfying (\ref{A_j}), with $\varepsilon \in ({\bf C},0)$ is proved 
using the basic technical tool, see Subsection~\ref{BTTAV}. Hence, for 
$\varepsilon \neq 0$ small enough the Jordan normal form of $\tilde{A}_j$ is 
$J_j^1$ (and $\tilde{A}_j$ is conjugate to $(G_j^0+\varepsilon G_j^1)$, see 
Proposition~\ref{sameorbit}). For these 
values of $\varepsilon$ the eigenvalues of $\tilde{A}_j$ will still be  
generic. 

Thus the existence of $(p+1)$-tuples 
with Jordan normal forms $J_j^0$ implies the existence of ones with Jordan 
normal forms $J_j^1$. 

$2^0$. By analogy one proves that the existence of irreducible 
$(p+1)$-tuples of matrices $A_j$ 
for the $(p+1)$-tuple of Jordan normal forms $J_j$ (and for some generic 
eigenvalues) implies the one for 
the $(p+1)$-tuple of Jordan normal forms $J_j'$ (and for some generic 
eigenvalues) where for each $j$ either 
$J_j'=J_j$ or $J_j$ is subordinate to $J_j'$. To this end one looks 
for the new $(p+1)$-tuple of matrices 

\[ \tilde{A}_j=(I+\varepsilon X_j(\varepsilon ))^{-1}Q_j^{-1}
(G_j+\varepsilon V_j(\varepsilon ))Q_j(I+\varepsilon X_j(\varepsilon ))\]
(where $G_j$ are Jordan matrices with Jordan normal forms $J_j$ and the 
matrices $V_j$ (holomorphic in $\varepsilon \in ({\bf C},0)$) 
are chosen such that $\tilde{A}_j$ have for $\varepsilon \neq 0$ Jordan 
normal form $J_j'$). The possibility to choose such matrices $V_j$ follows 
from the definitions of subordinate orbits and subordinate Jordan normal 
forms. This proves part 2) of the lemma for matrices $A_j$.

$3^0$. Assume that the Jordan matrices $G_j^i$ have the same meaning as in 
$1^0$. Choose strongly generic eigenvalues of $G_j^1$ such that for every $j$ 
they form a normalized chain (see the previous subsection).
Suppose that there exists a fuchsian 
system (\ref{Fuchs}) with $A_j$ conjugate to $G_j^1$ (for all $j$). Then for 
every $j$ the 
Jordan normal form ${J_j'}^0$ of the monodromy operator $M_j$ is either 
$J_j^0$ or is a Jordan normal form subordinate to $J_j^0$ 
(Lemma~\ref{JNFofM}). Such an irreducible  
monodromy group can be realized by a fuchsian system whose matrices-residua 
have Jordan normal forms ${J_j'}^0$ 
(by Lemma~\ref{non-resonant} -- such matrices-residua 
correspond to a non-resonant choice of the eigenvalues $\lambda _{k,j}$). 
By 2) of the present lemma, there exist such $(p+1)$-tuples of 
matrices-residua also for Jordan normal forms $J_j^0$. This 
proves part 1) of the lemma for matrices $A_j$.

$4^0$. Having proved the lemma for matrices $A_j$, one knows from 
Corollary~\ref{sameanswer} that it is true for matrices $M_j$ as well. 

The lemma is proved. $\hspace{2cm}\Box$

\section{The basic theorem for diagonalizable matrices
\protect\label{diagonal}}

\begin{defi} 
A simple PMV is called {\em good} if the DSP is solvable for some 
generic eigenvalues and for the $(p+1)$-tuple of diagonal 
Jordan normal forms defined by the PMV. For $n=1$ the only possible PMV 
is also defined as good.
\end{defi}

For the PMV $\Lambda ^n=(\Lambda ^n_1,\ldots ,\Lambda ^n_{p+1})$ 
(where $\Lambda ^n_j=(m_{1,j},\ldots ,m_{k_j,j})$,
$m_{1,j}+\ldots +m_{k_j,j}=n$) we presume that the following condition holds:

\[ m_{1,j}\geq \ldots \geq m_{k_j,j}~~~~~~~~~~~~(*_n)\]
Hence, $r_j=n-m_{1,j}$.

\begin{lm}\label{***}
A simple PMV satisfying condition $(\omega _n)$ is good.
\end{lm}

{\em Proof:} 

$1^0$. We use Definition~\ref{nicerepr} and 
Theorem~\ref{nilpunip}. For any $(p+1)$-tuple of diagonal Jordan normal 
forms $J_j^1$ one can 
find the $(p+1)$-tuple of corresponding Jordan normal forms with a single 
eigenvalue $J_j^0$, see Proposition~\ref{singleEV}. If the PMV of the 
Jordan normal forms $J_j^1$ is simple, then the Jordan normal 
forms $J_j^0$ do not correspond to any of the four exceptional cases cited 
in Theorem~\ref{nilpunip}. If condition $(\omega _n)$ holds 
for the Jordan normal forms $J_j^1$, 
then it holds for the Jordan normal forms $J_j^0$ as well 
(Proposition~\ref{samerank}). 

$2^0$. Hence, there exist nice $(p+1)$-tuples of nilpotent matrices $A_j$ or 
of unipotent matrices $M_j$ 
with Jordan normal forms $J_j^0$, see Definition~\ref{nicerepr} and 
Theorem~\ref{nilpunip}. The existence of 
irreducible $(p+1)$-tuples of matrices with Jordan normal forms $J_j^1$ is 
deduced from the existence for $J_j^0$ by means of the basic technical 
tool, see Subsections~\ref{BTTAV} and \ref{BTTMV}, by complete analogy with 
$1^0$ of the proof of Lemma~\ref{J^0impliesJ^1}.

The lemma is proved. $\hspace{2cm}\Box$

For a given simple PMV $\Lambda ^n$ define the numbers
$n=n_0>n_1>\ldots >n_s$ like this was done before Theorem~\ref{generic} by 
means of the map $\Psi$ 
(in our particular case of diagonalizable matrices we operate over  
PMVs instead of Jordan normal forms): if 
$\Lambda ^n$ satisfies
condition $(\omega _n)$ or if it does not satisfy condition $(\beta _n)$ or 
if $n=1$, then set $s=0$. If not, then set $n_1=r_1+\ldots
+r_{p+1}-n$. Hence, $n_1<n$ (otherwise $\Lambda ^n$
satisfies condition $(\omega _n)$).

Define the PMV $\Lambda ^{n_1}=(\Lambda _1^{n_1},\ldots ,\Lambda
_{p+1}^{n_1})$. Set $\Lambda
_j^{n_1,0}=$$(m_{1,j}-n+n_1,m_{2,j},\ldots ,m_{k_j,j})$ (recall that 
$m_{1,j}-n+n_1=n_1-r_1\geq 0$ because there hold conditions $(*_n)$ and 
$(\beta _n)$). For 
each $j$ rearrange the components of $\Lambda _j^{n_1,0}$ to obtain 
condition $(*_{n_1})$ -- this gives the MVs $\Lambda _j^{n_1}$.

Suppose that the PMVs $\Lambda ^{n_i}$ are constructed for $i=0,\ldots ,k$. 
If $\Lambda ^{n_k}$ satisfies condition $(\omega _{n_k})$ or if it does not 
satisfy condition $(\beta _{n_k})$ or if $n_k=1$, then set $s=k$.
If not, then define $n_{k+1}$ and $\Lambda ^{n_{k+1}}$ after $n_k$ and
$\Lambda ^{n_k}$ in the same way as $n_1$ and $\Lambda ^{n_1}$ were
defined after $n$ and $\Lambda ^n$ etc. In the end we have either
$n_s=1$ or $\Lambda ^{n_s}$ satisfies condition $(\omega _{n_s})$ or it 
does not satisfy condition $(\beta _{n_s})$.

\begin{rems}\label{psiinvariant}
1) Given a $(p+1)$-tuple of arbitrary Jordan normal forms $J_j^n$, 
construct the PMV $\Lambda ^n$ of the $(p+1)$-tuple of diagonal Jordan normal 
forms corresponding to them. Hence, the quantities $r_j$ and $n_1$ defined 
for both $(p+1)$-tuples coincide (this can be deduced from 
Proposition~\ref{samerank}). By 
Proposition~\ref{simultaneously}, the PMV $\Lambda ^{n_1}$ defines the 
$(p+1)$-tuple of diagonal Jordan normal forms corresponding to the Jordan 
normal forms $J_j^{n_1}$. In the same way one shows that the numbers 
$s$ and $n_1$, $\ldots$, $n_s$ are the same when defined for the 
$(p+1)$-tuple of Jordan normal forms $J_j^n$ and when defined for the PMV 
$\Lambda ^n$ and that for all $\nu$ the PMVs $\Lambda ^{n_{\nu}}$ define the 
$(p+1)$-tuples of diagonal Jordan normal forms corresponding to 
$J_j^{n_{\nu}}$.

2) For diagonal Jordan normal forms the PMVs $\Lambda ^{n_{\nu}}$ do not 
depend on the choice of eigenvalue in the map $\Psi$. This together with 
Remark~\ref{severalchoices} explains why 
Theorem~\ref{generic} is true whichever choice of eigenvalue is 
made in $\Psi$ -- for two such choices the two respective 
Jordan normal forms $J_j^{n_\nu}$ will correspond to one another for all 
$j$ and $\nu$ (and the sizes $n_{\nu}$ will be the same) because they 
correspond to one and the same PMVs $\Lambda ^{n_{\nu}}$. Hence, the PMV 
$\Lambda ^{n_s}$ and the $(p+1)$-tuple of Jordan normal forms $J_j^{n_s}$ 
satisfy or not condition $(\omega _n)$ (resp. $(\beta _n)$) simultaneously.
\end{rems}

\begin{lm}\label{simple}
If the PMV $\Lambda ^{n_{\nu}}$ is simple, then the PMV 
$\Lambda ^{n_{\nu +1}}$ is also simple.
\end{lm}

{\em Proof:} 

We prove the lemma for $\nu =0$, for arbitrary $\nu$ it is 
proved by analogy. Suppose that $\Lambda ^{n_1}$ is non-simple. Then for
every $j$ the greatest common divisor $l$ of its components divides
$m_{2,j}$, $\ldots$, $m_{k_j,j}$ and $m_{1,j}-n+n_1$, hence, it
divides $n_1$ (the length of $\Lambda ^{n_1}$). But 
$n_1=r_1+\ldots +r_{p+1}-n$ and $l$ 
divides $r_j$ (because $r_j=m_{2,j}+\ldots +m_{k_j,j}$); hence, $l$
divides $n$ and $m_{1,j}$ (because $m_{1,j}=n-r_j$). This means that  
$\Lambda ^n$ is non-simple -- a contradiction.$\hspace{2cm}\Box$

\begin{tm}\label{basicres}
A simple PMV $\Lambda ^n$ is good if and only if it satisfies condition 
$(\beta _n)$ and either the PMV $\Lambda ^{n_s}$ defined 
above satisfies condition 
$(\omega _{n_s})$ or one has $n_s=1$. The theorem is true both in the 
additive and in the multiplicative version of the DSP.
\end{tm}

\begin{defi}
For a given $(p+1)$-tuple of Jordan normal forms $J_j^n$ with $d_j=d(J_j^n)$ 
call the quantity $\kappa :=2n^2-d_1-\ldots -d_{p+1}$ {\em index of rigidity} 
of the $(p+1)$-tuple. (This notion was introduced by N.Katz in \cite{Ka}.) 
For a PMV define its {\em index of rigidity} as the one of the 
$(p+1)$-tuple of diagonal Jordan normal forms defined by it. 
\end{defi}

\begin{lm}\label{alsogood}
The PMVs $\Lambda ^{n_{\nu}}$ and $\Lambda^{n_{\nu +1}}$ have the same 
index of rigidity. In particular, they satisfy or not the respective 
conditions $(\alpha _{n_{\nu}})$ and $(\alpha _{n_{\nu +1}})$ simultaneously.
\end{lm}

{\em Proof:} We prove the lemma for $\nu =0$, for arbitrary $\nu$ the proof 
is analogous. Set $d_j=d(\Lambda ^n_j)$, $d_j^1=d(\Lambda ^{n_1}_j)$. One has 
$d_j^1=d_j-2(n-n_1)r_j$ (by direct computation) and

\[ \sum _{j=1}^{p+1}d_j^1=\sum _{j=1}^{p+1}d_j-2(n-n_1)\sum
_{j=1}^{p+1}r_j=\sum _{j=1}^{p+1}d_j-2(n-n_1)(n+n_1)=
\sum _{j=1}^{p+1}d_j-2n^2+2(n_1)^2 \]
which shows that the index of rigidity remains the same.$\hspace{2cm}\Box$

\begin{cor}\label{alsogood1}
If the $(p+1)$-tuple of Jordan normal forms $J_j^n$ (not necessarily 
diagonal) satisfies the equality $\sum _{j=1}^{p+1}d_j=2n^2-2+\chi$, 
$\chi \geq 0$, then for the quantities $d_j^{\nu}=d(J_j^{n_{\nu}})$ (where 
the Jordan normal forms $J_j^{n_{\nu}}$ are defined before 
Theorem~\ref{generic}) 
one has $d^{\nu}:=\sum _{j=1}^{p+1}d_j^{\nu}=2(n_{\nu})^2-2+\chi$. 
\end{cor}

Indeed, one can define the PMV $\Lambda ^n$ of the eigenvalues of the 
diagonal Jordan 
normal forms corresponding to $J_j^n$ and then the PMVs $\Lambda ^{n_{\nu}}$. 
For all $\nu$ the PMVs $\Lambda ^{n_{\nu}}$ 
define diagonal Jordan normal forms corresponding to $J_j^{n_{\nu}}$, see 
Remarks~\ref{psiinvariant}. For 
the $(p+1)$-tuples $\Lambda _j^{n_{\nu}}$ and $J_j^{n_{\nu}}$ the quantity 
$d^{\nu}$ is the same (Proposition~\ref{samedimension}).$\hspace{2cm}\Box$ 

\begin{rem}\label{sameindexofrigidity}
The conditions 
$n_s>1$ and $(\alpha _n)$ being a strict inequality are equivalent. 
Indeed, if $n_s>1$, then condition $(\omega _{n_s})$ holds for the 
$(p+1)$-tuple of Jordan normal forms $J_j^{n_s}$. By Lemma 3 from 
\cite{Ko1}, inequality $(\alpha _{n_s})$ holds for the Jordan normal forms 
$J_j^{n_s}$ and is strict. Corollary~\ref{alsogood1} allows to conclude 
that condition $(\alpha _n)$ is a strict inequality (i.e. $\chi >0$).
\end{rem}

\begin{cor}\label{n1}
If the PMV $\Lambda ^n$ is simple and good, then so are the PMVs 
$\Lambda ^{n_{\nu}}$, $\nu =1$, $\ldots$, $s$.
\end{cor}

The corollary follows from the definition of the PMVs $\Lambda ^{n_{\nu}}$,  
from Lemmas~\ref{alsogood} and \ref{simple} and from Theorem~\ref{basicres}.

{\em Proof of Theorem \ref{genericweak}:}

Proposition \ref{samerank} and Remarks~\ref{psiinvariant} show 
that the Jordan normal forms 
$J_j^n$ and $J_j^{n_s}$ satisfy conditions {\em i)} and {\em ii)} of 
Theorem~\ref{generic} if and only if the PMVs $\Lambda ^n$ and 
$\Lambda ^{n_s}$ satisfy the conditions of Theorem~\ref{basicres} (where 
$\Lambda ^n$ defines the diagonal Jordan normal forms corresponding to 
the $(p+1)$-tuple of Jordan normal forms $J_j^n$). 

By Lemma~\ref{J^0impliesJ^1}, there exist for some generic eigenvalues 
matrices $A_j$ or $M_j$ satisfying (\ref{A_j}) or (\ref{M_j}) with Jordan 
normal forms $J_j^n$ if and only if this is 
the case of diagonalizable matrices defined by the PMV $\Lambda ^n$. Thus 
Theorem~\ref{genericweak} results from 
Theorem~\ref{basicres}.$\hspace{2cm}\Box$

\section{Proof of the sufficiency in Theorem \protect\ref{basicres}
\protect\label{sufficiency}}
\subsection{Proof of the theorem itself}

The lemmas from this subsection are proved in the next ones. We prove the 
sufficiency in the case of matrices $A_j$, for matrices $M_j$ it follows 
then from Corollary~\ref{sameanswer}. 

Induction on $n$. For $n=1$ and 2 the theorem is checked
straightforwardly. If $\Lambda ^n$ satisfies condition $(\omega _n)$, then
$\Lambda ^n$ is good, see Lemma~\ref{***}. If not, then $\Lambda ^{n_1}$ 
satisfies the conditions of the theorem with $n$ replaced by $n_1$ (this 
follows from the definition of the PMVs $\Lambda ^{n_{\nu}}$ before 
Remarks~\ref{psiinvariant}). By 
inductive assumption, there exist (for generic eigenvalues) irreducible
$(p+1)$-tuples of diagonalizable $n_1\times n_1$-matrices $B_j$
(satisfying (\ref{A_j})) with PMV equal to $\Lambda ^{n_1}$. We assume
that the eigenvalue $\lambda _1$ of $B_1$ of multiplicity $m_{1,1}-n+n_1$ 
(when this multiplicity is not 0) equals 1 and that for $j>1$ the 
eigenvalue $\lambda _j$ of 
$B_j$ of multiplicity $m_{1,j}-n+n_1$ equals 0. This can be achieved by 
replacing the matrices $B_j$ by $B_j-\lambda _jI$, $j>1$, and 
$B_1$ by $B_1+(\lambda _2+\ldots +\lambda _{p+1})I$, and by multiplying 
all matrices by $c\in {\bf C}^*$. 

For the sake of convenience we make a circular permutation of the components 
of the MVs $\Lambda ^n_j$ and $\Lambda ^{n_1}_j$ putting their first 
components (i.e. $m_{1,j}$ and $m_{1,j}-n+n_1$) in last position.

Define the PMV $\tilde{\Lambda }^n$ as follows: for $j>1$ set
$\tilde{\Lambda }^n_j=\Lambda ^n_j$; set $\tilde{\Lambda
}_1^n=$$(m_{2,1}$,$\ldots$,$m_{k_1,1}$,$m_{1,1}-n+n_1$,$n-n_1)$. Define
the diagonalizable $n\times n$-matrices $A_j^0$ with PMV 
$\tilde{\Lambda }^n$ as $A_j^0=\left(
\begin{array}{cc}B_j&0\\0&0\end{array}\right)$, with $B_j$ as above. 
(The multiplicity of 0 as eigenvalue of $A_j^0$ equals $m_{1,j}$ for 
$j>1$ and $n-n_1$ for $j=1$.) Construct a $(p+1)$-tuple of matrices 
$A_j^1=\left( \begin{array}{cc}B_j&B_jY_j\\0&0\end{array}\right)$, 
$Y_j$ being $n_1\times (n-n_1)$, such that the monodromy operator $M_1^1$ at 
$a_1$ of the fuchsian system

\begin{equation}\label{B_jY_j}
\dot{X}=(\sum _{j=1}^{p+1}A_j^1/(t-a_j))X
\end{equation}
is diagonalizable (see Lemma \ref{M_1diag} below); we set $Y_1=0$. Notice 
that for each $j$ the 
matrix $A_j^1$ is conjugate to the matrix $A_j^0$. We assume that the only 
couple of eigenvalues of some matrix $A_j^1$ whose difference is a non-zero 
integer are the eigenvalues 0 and 1 of $A_1^1$. This is not restrictive, see 
Lemma~\ref{allgen}.

We also assume $A_1^0=A_1^1$ to be diagonal (hence, $B_1$ as well) and  
the eigenvalues of the matrices $B_j$ to be generic. 

\begin{lm}\label{M_1diag}
The operator $M_1^1$ is diagonalizable if and only if the following
conditions hold:

\begin{equation}
\sum _{j=2}^{p+1}(A_j^1/(a_j-a_1))_{\kappa ,\nu}=0~~~~,~~~~
\kappa =r_1+1,\ldots ,n_1~;~ \nu =n_1+1,\ldots ,n
\end{equation}
(double subscripts indicate matrix entries).
\end{lm}

\begin{rems}
1) The lemma is vacuous if $r_1=n_1$ when there is no condition 
to verify and $M_1^1$ is automatically diagonalizable.

2) If at least one of the matrices $Y_j$, $j>1$, is
non-zero, then the $(p+1)$-tuple of matrices $A_j^1$ is not conjugate
to the $(p+1)$-tuple of matrices $A_j^0$. Indeed, if this were the
case, then the conjugation should be carried out by a matrix 
commuting with $A_1^0=A_1^1$, i.e. block-diagonal, with diagonal blocks 
of sizes $n_1\times n_1$ and $(n-n_1)\times (n-n_1)$. Such a conjugation 
cannot annihilate the blocks $B_jY_j$.
\end{rems}

\begin{defi}
We say that the columns of the $(p+1)$-tuple of 
$q\times r$-matrices  
$C_j$ are linearly independent if for no $r$-tuple of constants 
$\beta _i\in {\bf C}$ (not all of them being 0) one has $\sum
_{i=1}^r\beta _iC_{j,i}=0$ for $j=1,\ldots ,p+1$ where $C_{j,i}$
is the $i$-th column of the matrix $C_j$. In the same way one defines 
independence of rows.
\end{defi}

\begin{defi}
Denote by $\tilde{{\bf C}}^v$ the linear space of 
$p$-tuples of vectors $T_j\in {\bf C}^{n_1}$, $j=2,\ldots ,p+1$, where 
\begin{equation}\label{T_2=0}
T_2+\ldots +T_{p+1}=0
\end{equation}
and $T_j=B_jU_j$ for some $U_j\in {\bf C}^{n_1}$. Denote by 
$\tilde{{\bf C}}^w\subset \tilde{{\bf C}}^v$ its subspace satisfying the 
condition 
\begin{equation}\label{alphaT_2=0}
(\alpha _2T_2+\ldots +\alpha _{p+1}T_{p+1})|_{\kappa }=0~,~\alpha _j=
1/(a_j-a_1)~,~\kappa=r_1+1,\ldots ,n_1
\end{equation}
The notation $|_{\kappa }$ means the $\kappa$-th 
coordinate of the vector, see Lemma \ref{M_1diag}.
\end{defi}

It is clear that $v\stackrel{{\rm def}}{=}\dim \tilde{{\bf C}}^v
\geq r_2+\ldots +r_{p+1}-n_1=n-r_1$ (the image of the linear operator 
$\xi _j:{\bf C}^{n_1}\rightarrow {\bf C}^{n_1}$, $\xi _j:(.)\mapsto B_j(.)$ 
is of dimension $r_j$ and equation (\ref{T_2=0}) is equivalent to $\leq n_1$ 
linearly independent equations). In the same way one deduces the inequality 
$w\stackrel{{\rm def}}{=}\dim \tilde{{\bf C}}^w\geq v-(n_1-r_1)\geq n-n_1$.

\begin{lm}\label{vw}
One has $v=n-r_1$ and $w=n-n_1$.
\end{lm}

\begin{lm}\label{fourconditions}
There exists a $(p+1)$-tuple of matrices $Y_j$ such that

1) $Y_1=0$ and for $j>1$ $Y_j$ belongs to the image of the
linear operator $\tau _j:(.)\mapsto B_j(.)$ acting on the space of 
$n_1\times (n-n_1)$-matrices;

2) $\sum _{j=2}^{p+1}B_jY_j=0$;

3) the monodromy operator $M_1^1$ at $a_1$ of the fuchsian system
(\ref{B_jY_j}) is diagonalizable;

4) the columns of the $p$-tuple of matrices $B_jY_j$, $j=2,\ldots ,p+1$, are
linearly independent; they are a basis of the space $\tilde{{\bf C}}^w$ 
defined above. 
\end{lm}

\begin{lm}\label{trivcentral}
1) The centralizer of the $(p+1)$-tuple of matrices $A_j^1$ satisfying 1)
-- 4) of the previous lemma is trivial. 

2) The centralizer of the monodromy group of system (\ref{B_jY_j}) (in which 
the $(p+1)$-tuple of matrices $A_j^1$ satisfies 1) -- 4) of 
Lemma~\ref{fourconditions}) is trivial. 
\end{lm}

Lemma~\ref{trivcentral} is necessary for the proof of the following lemma 
from which follows the proof of the sufficiency.
 
\begin{lm}\label{existM_j}
Denote by ${\cal L}$ the set of eigenvalues of the 
$(p+1)$-tuple of matrices $M_j^1$. 
There exist (for generic eigenvalues close to ${\cal L}$) irreducible 
$(p+1)$-tuples of diagonalizable matrices $M_j$ satisfying (\ref{M_j}), with 
PMV equal to $\Lambda ^n$.
\end{lm}
 
By Lemma \ref{allgen} there exist such $(p+1)$-tuples for all 
eigenvalues from a Zariski open dense subset of the set of all generic 
eigenvalues with multiplicities defined by $\Lambda ^n$. Thus we have proved 
that in the multiplicative 
version the simple PMV $\Lambda ^n$ is good. It is good in the additive 
one as well due to Corollary~\ref{sameanswer}.

The sufficiency is proved.

\subsection{Proof of Lemma \protect\ref{M_1diag}
\protect\label{proofsoflemmassuff}}

The fuchsian 
system (\ref{B_jY_j}) represented by its Laurent 
series at $a_1$ looks like this:

\[ \dot{X}=[A_1^1/(t-a_1)+B+o(1)]X~~,~~
B=-(\sum _{j=2}^{p+1}A_j^1/(a_j-a_1))\]
One can assume that $A_1^1=$diag$(\lambda _{1,1},\ldots ,\lambda _{n,1})$ 
where $\lambda _{r_1+1,1}=\ldots =\lambda _{n_1,1}=1$, 
$\lambda _{n_1+1,1}=\ldots =\lambda _{n,1}=0$. The local (at $a_1$) change 
of variables

\[ X\mapsto {\rm diag}(1,\ldots ,1,(t-a_1)^{-1},\ldots ,(t-a_1)^{-1})X\]
($n-n_1$ times $(t-a_1)^{-1}$) brings the system to the form 

\[ \dot{X}=[A_1^*/(t-a_1)+O(1)]X\]
where the matrix $A_1^*$ 

1) is upper-triangular; 

2) has no non-zero integer 
differences between its eigenvalues; 

3) has an 
eigenvalue 1 of multiplicity $m_{1,1}$ occupying the last $m_{1,1}$ 
positions on its diagonal;  

4) its right lower $m_{1,1}\times m_{1,1}$-block 
equals $\left( \begin{array}{cc}I&\Delta \\0&I\end{array}\right)$; here 
$\Delta$ is the restriction of the matrix $B$ to the last $n-n_1$ columns 
intersected with the rows with indices $r_1+1$, $\ldots$, $n_1$.

Hence, the eigenvalues of $A_1^*$ are non-resonant and the monodromy operator 
$M_1^1$ is conjugate to $\exp (2\pi iA_1^*)$, see Lemma~\ref{sameJNF}. 

Hence, $A_1^*$ and $M_1^1$ are diagonalizable if and only if $\Delta =0$. 
This proves the lemma.

\subsection{Proof of Lemma \protect\ref{vw}}

$1^0$. Multiply the matrices $B_j$ by $c\in {\bf C}^*$ so that the new 
matrices $B_j$ have strongly generic eigenvalues (the lemma is true or not 
simultaneously for the old and for the new matrices). It suffices to prove 
the second equality which would imply that both 
inequalities $v\geq n-r_1$ and $w\geq v-(n_1-r_1)$ are equalities. The 
equality is 
true exactly if conditions (\ref{T_2=0}) and (\ref{alphaT_2=0}) together are 
linearly independent. We consider them as a system of linear equations 
with unknown variables the entries of the vectors $U_j\in {\bf C}^{n_1}$ 
where $T_j=B_jU_j$. Their linear dependence is equivalent to the statement:

{\em there exist vector-rows $V,W\in {\bf C}^{n_1}$, $(V,W)\neq (0,0)$, such 
that} 

\begin{equation}\label{alphavw} 
(V+\alpha _jW)B_j=0~{\rm for~}j=2,\ldots ,p+1~{\rm and~}WB_1=W
\end{equation} 
Indeed, if (\ref{T_2=0}) and (\ref{alphaT_2=0}) together are not linearly 
independent, then some non-trivial linear combination of theirs is of the 
form $0=0$. This linear combination is of the form 

\[ V(\sum _{j=2}^{p+1}B_jU_j)+W(\sum _{j=2}^{p+1}\alpha _jB_jU_j)=
\sum _{j=2}^{p+1}(V+\alpha _jW)B_jU_j=0\] 
Its left hand-side must be identically 0 in the entries of $U_j$, 
i.e. $(V+\alpha _jW)B_j=0$ for $j=2,\ldots ,p+1$. The condition $WB_1=W$ 
follows from $\kappa =r_1+1, \ldots ,n_1$, see (\ref{alphaT_2=0}); recall 
that $B_1$ is diagonal and that its last eigenvalue equal to 1 occupies the 
positions with indices $r_1+1$, $\ldots$, $n_1$, therefore $W$ is left 
eigenvector of $B_1$ corresponding to the eigenvalue 1. 

$2^0$. Consider the fuchsian system

\begin{equation}\label{B} 
\dot{X}=A(t)X~,~A(t)=\left( \sum _{j=1}^{p+1}\left( 
\begin{array}{cc}B_j&0\\0&0\end{array}\right) /(t-a_j)\right) X
\end{equation}
of dimension $n_1+1$. Perform the change $X\mapsto R(t)X$, $R(t)=\left( 
\begin{array}{cc}I&0\\V+W/(t-a_1)&1\end{array}\right)$. The matrix $A(t)$ 
changes to $-R^{-1}\dot{R}+R^{-1}A(t)R$. One can check directly that $A(t)$ 
does not change under the above change of variables (i.e. 
$-\dot{R}+A(t)R=RA(t)$) if and only if conditions (\ref{alphavw}) hold 
(after the change the system is fuchsian at $a_j$ for $j>1$ and the residuum 
equals 
$\left( \begin{array}{cc}B_j&0\\-(V+W/(a_j-a_1))B_j&0\end{array}\right)$; 
its polar part at $a_1$ equals 
$\left( \begin{array}{cc}0&0\\-WB_1+W&0\end{array}\right) /(t-a_1)^2+$
$\left( \begin{array}{cc}B_1&0\\-VB_1+\sum _{j=2}^{p+1}WB_j/(a_j-a_1)&0
\end{array}\right) /(t-a_1)$). 

$3^0$. The solution to system (\ref{B}) with initial data $X|_{t=a_0}=I$ 
changes from $X$ to $R(t)X$ and this must be again a solution to system 
(\ref{B}) (because the system does not change). 
Hence, $R(t)X=XD$ for some $D\in GL(n_1+1,{\bf C})$. The solution $X$ is 
block-diagonal (with blocks $n_1\times n_1$ and $1\times 1$) for all values 
of $t$ due to the block-diagonal form of the 
system and, hence, the one of the monodromy group as well. 

The first $n_1$ 
coordinates of the last column of the matrix $R(t)X$ are identically zero 
and its restriction to $H$ (the left upper $n_1\times n_1$-block) are 
identically equal 
to the ones of $X$. This together with the linear independence of the 
columns of $X|_H$ implies the form of the matrix $D$: 
$D=\left( \begin{array}{cc}I&0\\C&g\end{array}\right)$. 

The conditions $X|_{t=a_0}=I$ and $R(t)X=XD$ imply $g=1$. 

$4^0$. The analytic continuations of $R(t)X$ and $XD$ coincide, therefore 
for every monodromy operator $M'_j$ of the system one must have 
$R(t)XM'_j$=$XM'_jD$. But one has $R(t)XM'_j=XDM'_j$, i.e. $[M'_j,D]=0$ for 
every monodromy operator.

The monodromy operators are block-diagonal: 
$M'_j=\left( \begin{array}{cc}M''_j&0\\0&1\end{array}\right)$ and the group 
${\cal G}\subset GL(n_1,{\bf C})$ generated by the operators $M''_j$ is 
irreducible (this follows from the strong 
genericity of the eigenvalues of the matrices $B_j$). 

The condition $[M'_j,D]=0$ implies $(M''_j-I)C=0$ for all $j$. This together 
with the irreducibility of the group ${\cal G}$ yields $C=0$. But then 
$R(t)=I$, i.e. $V=W=0$ which proves the lemma.

\subsection{Proof of Lemma \protect\ref{fourconditions}} 

The space $\Theta$ of matrices $Y_j$ satisfying 1) is of dimension
$(r_2+\ldots +r_{p+1})(n-n_1)=(n+n_1-r_1)(n-n_1)$ (for $j\geq 2$ the 
dimension of the image of $\tau _j$ is $(n-n_1)r_j$). Its subspace $\Phi$
defined by 2) is of codimension $(n-n_1)n_1$ in $\Theta$, hence, of 
dimension $(n-r_1)(n-n_1)$. This follows from $\dim \tilde{{\bf C}}^v=n-r_1$, 
see Lemma~\ref{vw}, because one has 
$\Phi =\tilde{{\bf C}}^v\times \ldots \times \tilde{{\bf C}}^v$ ($n-n_1$ 
times).

The subspace $\Xi$ of $\Phi$ defined by 
condition 3) is of codimension $(n_1-r_1)(n-n_1)$ in $\Phi$ (see 
Lemmas~\ref{M_1diag} and \ref{vw} -- $\dim \tilde{{\bf C}}^w=n-n_1$), i.e. 
of dimension $(n-n_1)^2$. 

This dimension is $n-n_1$ times the dimension of the
space $\tilde{{\bf C}}^w$ of vector-columns $Y_j$ of length $n_1$ (instead 
of $n_1\times (n-n_1)$-matrices) which satisfy 1) -- 3) of the conditions of 
the lemma. Indeed, one has 
$\Xi=\tilde{{\bf C}}^w\times \ldots \times \tilde{{\bf C}}^w$ ($n-n_1$ 
times).

By Lemma \ref{vw}, dim$\tilde{{\bf C}}^w=n-n_1$, i.e. one can choose exactly 
$n-n_1$ $(p+1)$-tuples of 
vector-columns satisfying conditions 1) -- 3) of the lemma which are
linearly independent. The exactitude implies that they are a basis 
of the space $\tilde{{\bf C}}^w$. Hence, the choice of matrices $Y_j$ 
satisfying 1) -- 4) is also possible. 

The lemma is proved.

\subsection{Proof of Lemma \protect\ref{trivcentral}} 

$1^0$. Prove 1). A matrix $Z$ 
commuting with $A_1^1$ must be of the form 
$Z=\left( \begin{array}{cc}Z'&0\\0&Z''\end{array}\right)$, $Z''$ 
being $(n-n_1)\times (n-n_1)$. 

One must have $Z'=\alpha I$, $\alpha \in {\bf C}$ due to Schur's lemma 
because the $(p+1)$-tuple of matrices $B_j$ is irreducible, one has 
$[A_j^1,Z]=0$ for all $j$ and, hence, $[B_j,Z']=0$. 

Hence, for all $j$ one has $\alpha B_jY_j=B_jY_jZ''$. The linear 
independence of the columns of the $(p+1)$-tuple of matrices $B_jY_j$ 
implies $Z''=\alpha I$. Part 1) of the lemma is proved.

$2^0$. Prove 2). Let $X|_{t=a_0}=I$, $a_0\neq a_j$, $j=1$, $\ldots$, $p+1$. 
One can conjugate the monodromy operators (defined for these initial data) 
to the same form as the one of the matrices-residua: 
$M_j^1=\left( \begin{array}{cc}N_j&N_jW_j\\0&I\end{array}\right)$, 
$N_j$ being $n_1\times n_1$, with $W_1=0$. If it were known that the columns 
of the $(p+1)$-tuple of matrices $N_jW_j$ are independent, then part 2) of 
the lemma could be proved like part 1). So suppose that this is not the 
case. 

A conjugation with a matrix 
$\tilde{D}=\left( \begin{array}{cc}I&0\\0&D\end{array}\right)$, 
$D\in GL(n-n_1,{\bf C})$, brings the matrices $M_j^1$ to the form 
$M_j^1=\left( \begin{array}{ccc}N_j&N_jW_j'&0\\0&I&0\\0&0&1\end{array}\right)$ 
where $I$ is $(n-n_1-1)\times (n-n_1-1)$. This conjugation is tantamount to 
the change of the basis of the solution space: $X\mapsto X\tilde{D}$. 

$3^0$. Perform the change of the dependent variables 
$\eta :X\mapsto \left( \begin{array}{cc}I&0\\0&D^{-1}\end{array}\right) X$. 
This changes system (\ref{B_jY_j}) but preserves its block upper-triangular 
form, the form of its monodromy operators and the size of the blocks $B_j$. 

%

Hence, after the change for any value of $t$ the solution $X$ is of the form 
$X=\left( \begin{array}{cc}X'&X''\\0&I\end{array}\right)$ because the 
derivative of any entry of the last $n-n_1$ rows is 0 (recall that the 
last $n-n_1$ eigenvalues of the matrices $A_j^1$ before and after the change 
$\eta$ are 0). Moreover, one has $X|_{t=a_0}=I$.

$4^0$. The form of the monodromy operators implies that each entry of the 
last column $X^n$ of $X$ is a meromorphic (i.e. univalued) function on 
${\bf C}P^1$. Moreover, the last $n-n_1$ entries of $X^n$ equal identically 
0, $\ldots$, 0, 1. 

Hence, there 
exists a change of variables $X\mapsto V(t)X$ (with 
$V(t)=\left( \begin{array}{cc}I&\tilde{V}(t)\\0&1\end{array}\right)$,  
the matrix-function $\tilde{V}$ being meromorphic on 
${\bf C}P^1$, its last $n-n_1-1$ entries being identically 0)  
after which the new 
matrix-solution $VX$ is of the form 
$\left( \begin{array}{ccc}X'&X'''&0\\0&I&0\\0&0&1\end{array}\right)$
($I$ being $(n-n_1-1)\times (n-n_1-1)$). Show that system 
(\ref{B_jY_j}) becomes after this change fuchsian again and block-diagonal. 

$5^0$. Indeed, under the 
change $X\mapsto V(t)X$ the linear system $\dot{X}=A(t)X$ undergoes the gauge 
transformation $A(t)\rightarrow C(t)=-V^{-1}(t)\dot{V}(t)+V^{-1}(t)A(t)V(t)$. 
Hence, the left $(n-1)$ columns of the matrix $A(t)$ from system 
(\ref{B_jY_j}) do not change at all (we use 
the fact that the last row of $A(t)$ equals $(0,\ldots ,0,0)$). The last 
column of the new matrix $A(t)=\dot{(VX)}(VX)^{-1}$ is identically 0, see 
the form of $VX$. Hence, the poles of $C(t)$ are of first order and its 
matrices-residua are of the form 
$A_j^2=\left( \begin{array}{ccc}B_j&G_j&0\\0&0&0\\0&0&0\end{array}\right)$, 
$G_j$ being $n_1\times (n-n_1-1)$.

$6^0$. We show that the matrix $V$ is constant, see $7^0$. This implies that 
the $(p+1)$-tuple of matrices-residua $A_j^2$ is conjugate to the 
$(p+1)$-tuple of matrices-residua $A_j^1$ and, hence, the columns of the 
$(p+1)$-tuple of matrices $B_jY_j$ are dependent -- a contradiction. This 
contradiction proves part 2) the lemma.

$7^0$. The change $X\mapsto V(t)X$ from $5^0$ preserves up to conjugacy 
the matrices-residua $A_j^1$ for $j>1$. Represent system (\ref{B_jY_j}) and 
the matrix $V$ in the neighbourhood of $a_j$ by their Laurent series: 

\[ \dot{X}=[A_j^1/(t-a_j)+O(1)]X~,~V=V^*/(t-a_j)^k+o(1/(t-a_j)^k)~,~k\geq 0.
\]

One has (see $5^0$) $C(t)=(A_j^1)'/(t-a_j)+O(1)$ in the neighbourhood of 
$a_j$ where $(A_j^1)'\in gl(n,{\bf C})$ is conjugate to $A_j^1$. 
The equation $VC=-\dot{V}+AV$ implies $V^*(A_j^1)'=kV^*+A_j^1V^*$ (these are 
the coefficients before $1/(t-a_j)^{k+1}$). 

If $k>0$, 
this equation yields $V^*=0$. Indeed, the eigenvalues of the linear operator 
$(.)\mapsto -(.)(A_j^1)'+k(.)+A_j^1(.)$ acting on $gl(n,{\bf C})$ equal 
$\lambda '=-\lambda _{\nu ,j}+k+\lambda _{\mu ,j}$. The absence of non-zero 
integer differences between the eigenvalues of $A_j$ for $j>1$ implies that 
$\lambda '\neq 0$, hence, $V^*=0$, i.e. $V$ has no pole at $a_j$ for $j>1$. 

The form of the last column of the solution $X$ at $a_1$ and the one of $VX$ 
imply that $V(a_1)=I$, i.e. $V$ has 
no pole at $a_1$ either, hence, no poles on ${\bf C}P^1$, i.e. $V$ is 
constant.

The lemma is proved.

\subsection{Proof of Lemma \protect\ref{existM_j}} 

Apply the basic technical tool in the multiplicative version, see 
Subsection~\ref{BTTMV}. 
To prove the lemma it suffices to choose for each $j$ a matrix $N_j$ which 
is a suitable 
polynomial of $M_j^1$. The $(p+1)$-tuple of matrices $M_j$ is with trivial 
centralizer, but can be reducible. Choose $N_j$ such that for 
$\varepsilon \neq 0$ the eigenvalues of the matrices $M_j$ to be generic. 
Hence, the $(p+1)$-tuple of matrices $M_j$ will be irreducible for 
$\varepsilon \neq 0$. 

The lemma is proved.

\section{Proof of the necessity in Theorem \protect\ref{basicres}
\protect\label{necessity}}
\subsection{Proof of the theorem itself\protect\label{schemeoftheproof}}

$1^0$. In this section we consider system (\ref{Fuchs}) with generic but not 
strongly generic eigenvalues, with diagonalizable matrices $A_j$
whose PMV $\Lambda ^n$ is simple and good. Without loss of generality we 
assume that for $j=2$, $\ldots$, $p+1$ one of the 
eigenvalues of greatest multiplicity of $A_j$ is 0 and 
for $j=1$ one of them equals 1 (the last  
condition is obtained by multiplying the 
residua by $c\in {\bf C}^*$). Hence, the corresponding eigenvalues 
$\sigma _{k,j}$ of the 
matrices $M_j$ equal 1, i.e. they satisfy at least one non-genericity 
relation (denoted by $(\gamma ^0)$). None of the other eigenvalues 
$\lambda _{k,j}$ is integer. 

$2^0$. We assume that for all $j$ 
the eigenvalues of $A_j$ are non-resonant. We assume 
also that 

A) either $(\gamma ^0)$ is the only non-genericity relation that the 
eigenvalues $\sigma _{k,j}$ 
satisfy or the greatest common divisor $l$ of the multiplicities of the 
non-integer eigenvalues of all matrices $A_j$ is $>1$; if 
$l=1$, then it is 
possible to choose the eigenvalues $\lambda _{k,j}$ so that the eigenvalues 
$\sigma _{k,j}$ satisfy only the non-genericity relation $(\gamma ^0)$ and 
no other. If $l>1$, 
then one can divide by $l$ the multiplicities of 
the eigenvalues $\sigma _{k,j}$ which 
are not 1 -- their product (which is a priori a root of unity of order $l$, 
see (\ref{sumEVszero})) 
might turn out to be a non-primitive such root. This could give rise to 
another 
non-genericity relation $(\gamma ^1)$. In this case one can choose the 
eigenvalues $\lambda _{k,j}$ so that every non-genericity relation 
satisfied by the eigenvalues $\sigma _{k,j}$ should be 
a linear combination of $(\gamma ^0)$ and $(\gamma ^1)$;

B) neither $n=1$, nor the PMV $\Lambda ^n$ satisfies condition $(\omega _n)$ 
(in which cases there is nothing to prove).

$3^0$. Assumption B) above implies that the $(p+1)$-tuple of matrices $M_j$ 
must be reducible -- part 2) of Proposition~\ref{necessarynongen} does not 
hold 
(recall that 1 is eigenvalue of greatest multiplicity for all $j$; hence, 
rk$(M_j-I)=r_j$; if one sets $b_j=1$, then the necessary 
condition for existence of irreducible $(p+1)$-tuples coincides with 
condition $(\omega _n)$ which does not hold).

$4^0$. 

\begin{lm}\label{reprtrivcentr}
The monodromy group of a fuchsian system with generic non-resonant 
eigenvalues of the matrices-residua 
is with trivial centralizer. In particular, the monodromy group of system 
(\ref{Fuchs}) with eigenvalues defined as above is with trivial centralizer. 
\end{lm}

All lemmas from this subsection are proved in the next ones.

\begin{lm}\label{theMGisreducible}
The monodromy group of system (\ref{Fuchs}) with eigenvalues defined as 
above can be conjugated to the form $\left( \begin{array}{cc}\Phi &\ast \\
0&I\end{array}\right)$ where $\Phi$ is $n_1\times n_1$.
\end{lm}

\begin{rem}
Notice that the subrepresentation $\Phi$ can be reducible.
\end{rem} 

\begin{lm}\label{subreprtrivcentr}
The centralizer $Z(\Phi )$ of the subrepresentation $\Phi$ is trivial.
\end{lm}

$5^0$. The subrepresentation $\Phi$ being of dimension 
$n_1<n$, one can use induction on $n$ to prove the necessity. For $n=1$ and 
$2$ the 
necessity is evident. The PMV of the matrices $M_j'$ defining $\Phi$ equals 
$\Lambda ^{n_1}$. It follows from 
Lemma~\ref{subreprtrivcentr} that for generic eigenvalues close to the ones 
of the matrices $M_j'$ defining $\Phi$ there exist irreducible 
$(p+1)$-tuples of 
diagonalizable matrices $\tilde{M}_j'\in GL(n_1,{\bf C})$ with PMV 
$\Lambda ^{n_1}$ and satisfying (\ref{M_j}) 
(this can be proved by complete analogy with Lemma~\ref{existM_j}, by using 
the basic technical tool in the multiplicative version; recall 
that the triviality of the centralizer was essential in the proof of 
Lemma~\ref{existM_j} and was assured by Lemma~\ref{trivcentral}).

Hence, if $\Lambda ^n$ is good, then $\Lambda ^{n_1}$ is good. The necessity 
of $(\beta _n)$ was 
proved in Proposition~\ref{d_jr_j} and condition $(\omega _n)$ 
does not hold by assumption. Finally, the PMV $\Lambda ^{n_s}$ is the same for 
$\Lambda ^n$ and for $\Lambda ^{n_1}$ (this follows from the definition of 
the PMVs $\Lambda ^{n_{\nu}}$ before Remarks~\ref{psiinvariant} -- the PMV 
$\Lambda ^{n_s}$ is the last of this chain of PMVs). If 
$\Lambda ^{n_1}$ is good, then either $\Lambda ^{n_s}$ satisfies condition 
$(\omega _{n_s})$ or one has $n_s=1$. Hence, if the 
PMV $\Lambda ^n$ is good, then it satisfies the conditions of 
Theorem~\ref{basicres}, i.e. they are necessary. 

The necessity holds in both versions (additive and multiplicative), 
see Corollary~\ref{sameanswer}.

The necessity is proved.

\subsection{Proof of Lemma \protect\ref{reprtrivcentr}} 

$1^0$. Suppose the lemma not to be true. Then the centralizer 
either contains a diagonalizable matrix $D$ with at least two distinct 
eigenvalues or it contains a nilpotent matrix $N\neq 0$. 
(Indeed, let $A$ and $S$ be respectively a Jordan matrix and its 
semisimple part. If $[X,A]=0$, then $[X,S]=0$.) In the first case we can 
assume that $D$ has exactly two eigenvalues which can be achieved by 
considering instead of $D$ some suitable polynomial of it. In the second 
case without restriction one can assume that $N^2=0$ (by considering $N^k$ 
instead of $N$ for some $k\in {\bf N}$). 

$2^0$. In the first case one conjugates $M_j$ and $D$ to the form 
$M_j=\left( \begin{array}{cc}M_j^1&0\\0&M_j^2\end{array}\right)$, 
$D=\left( \begin{array}{cc}\alpha I&0\\0&\beta I\end{array}\right)$, 
$\alpha \neq \beta$. The form of $M_j$ follows from $[M_j,D]=0$. Apply 
Proposition~\ref{bolibr} twice -- once to each of the subspaces defined by 
the $(p+1)$-tuples of blocks $M_j^1$ and $M_j^2$. One sees that the sums 
$\lambda ^1$, $\lambda ^2$ of 
the eigenvalues $\lambda _{k,j}$ 
corresponding to $M_j^1$, $M_j^2$ must be $\leq 0$. On the other hand, there 
holds $\lambda ^1+\lambda ^2=0$, see (\ref{sumEVszero}), hence, 
$\lambda ^1=\lambda ^2=0$. This non-genericity relation 
contradicts the genericity of the eigenvalues. We used the fact that the 
eigenvalues are non-resonant -- knowing the eigenvalues $\sigma _{k,j}$ 
of the blocks $M_j^i$, we know the corresponding eigenvalues $\lambda _{k,j}$ 
as well (the absence of non-zero 
integer differences (for fixed $j$) between the eigenvalues $\lambda _{k,j}$ 
implies that to equal eigenvalues of $M_j$ there correspond equal 
eigenvalues of $A_j$).

$3^0$. In the second case the matrices $M_j$ and the matrix $N$ can be 
conjugated respectively to the form

\[ M_j=\left( \begin{array}{ccc}M_j'&\ast &\ast \\
0&M_j''&\ast \\0&0&M_j'\end{array}\right) ~,~
N=\left( \begin{array}{ccc}0&0&I\\0&0&0\\0&0&0\end{array}\right) \]
where $I$ is $w\times w$, $w\leq n/2$; if $w=n/2$, then the blocks of the 
middle columns and of the middle rows are absent. The form of $M_j$ 
follows from $[M_j,N]=0$. 

$4^0$. By Proposition \ref{bolibr}, the sum $\lambda '$ of the 
eigenvalues $\lambda _{k,j}$ corresponding to the block 
$\left( \begin{array}{cc}M_j'&\ast \\
0&M_j''\end{array}\right)$ must be non-positive and the sum $\lambda ''$ of 
the eigenvalues $\lambda _{k,j}$ corresponding to the block 
$\left( \begin{array}{cc}M_j''&\ast \\
0&M_j'\end{array}\right)$ must be non-negative (because there holds 
(\ref{sumEVszero}) and the sum of the eigenvalues $\lambda _{k,j}$ 
corresponding to the upper blocks $M_j'$ is $\leq 0$ by 
Proposition~\ref{bolibr}). One must have 
$\lambda '=\lambda ''$ (use like in $2^0$ the fact that the eigenvalues are 
non-resonant). It follows from (\ref{sumEVszero}) that 
$\lambda '=\lambda ''=0$. This non-genericity 
condition contradicts the genericity of the eigenvalues $\lambda _{k,j}$. 

The lemma is proved.

\subsection{Proof of Lemma \protect\ref{theMGisreducible}}

$1^0$. The monodromy group can 
be conjugated to a block upper-triangular form. The diagonal blocks define 
either irreducible or one-dimensional representations. The eigenvalues of 
each diagonal block $1\times 1$ 
satisfy the non-genericity relation $(\gamma ^0)$. 

$2^0$. Recall that the integer $l$ was defined in $2^0$ of 
Subsection~\ref{schemeoftheproof}. 
The block in the right  
lower corner must be of size 1. Indeed, if $l=1$, then by 
Proposition~\ref{bolibr} the left upper block cannot be of size 1 (because 
the corresponding sum of eigenvalues $\lambda _{k,j}$ equals $1>0$). Hence, 
it must be the only block of size $>1$ and the matrices $M_j$ 
look like this: 
$M_j=\left( \begin{array}{cc}M_j'&L_j\\0&I\end{array}\right)$. 

The block $M'$ must be of size $\leq n_1$. 
Indeed, if its size is $>n_1$ (i.e. this is the only diagonal block of 
size $>1$), then the columns of the $(p+1)$-tuples of 
matrices $L_j$ are not linearly independent (this is proved by complete 
analogy with the proof of 4) of Lemma~\ref{fourconditions}). 

This proves the lemma in the case $l=1$.

$3^0$. Let $l>1$. In the absence of second non-genericity relation 
$(\gamma ^1)$ (it was defined in $2^0$ of Subsection~\ref{schemeoftheproof}) 
the proof is finished like in $1^0$ -- $2^0$. So suppose that $(\gamma ^1)$ 
holds. The diagonal blocks can be of two types. The first are 
of size 1, the eigenvalues satisfying the non-genericity relation 
$(\gamma ^0)$. 

Describe the second type of diagonal blocks. Their sizes are $>1$ and can 
be different. Define the {\em unitary set} of eigenvalues: for each $j$ 
divide by $l$ the multiplicities of all eigenvalues 
$\sigma _{k,j}$ of the ones that are $\neq 1$. A 
block $F$ of the second type contains $h$ times the unitary set, 
$1\leq h\leq l$, and a certain number of eigenvalues equal to 1. 
(To different matrices $M_j$ there correspond, 
in general, different numbers of eigenvalues from the unitary set; therefore 
one must, in general, add some number of eigenvalues 1 for some values of 
$j$ to make the number of 
eigenvalues of the restrictions of the matrices $M_j$ to $F$  
equal; one then could eventually add one and the same number of eigenvalues 
equal to 1 to all matrices $M_j|_F$.) 
 
The eigenvalues of each block of the second type satisfy a corollary of the 
non-genericity relations $(\gamma ^1)$ and $(\gamma ^0)$.

$4^0$. Denote by $\kappa (F)$ the 
ratio "number of eigenvalues $\sigma _{k,j}$ equal to 1"/"number of 
eigenvalues $\sigma _{k,j}$ not equal to 1" (eigenvalues of the 
restriction of the monodromy group to $F$), and by $\kappa _0$ the same 
ratio computed for the entire matrices $M_j$ (in both ratios one takes into 
account the eigenvalues of all matrices $M_j$). Then one must have 
$\kappa (F)<\kappa _0$.

Indeed, one cannot have $\kappa (F)\geq \kappa _0$ 
because $\Lambda ^n$ does not satisfy condition $(\omega _n)$, hence, 
the restriction of the monodromy group to $F$ wouldn't satisfy this 
condition either. In the presence of the non-genericity relation 
$(\gamma ^0)$ this implies a contradiction with 
Proposition~\ref{necessarynongen}
(like in $3^0$ of Subsection~\ref{schemeoftheproof}). 

But then the sum $\tilde{\lambda}$ of the eigenvalues $\lambda _{k,j}$ 
corresponding to the 
eigenvalues $\sigma _{k,j}$ from $F$ will be negative. If the block $F$ is 
to be in the right lower corner, then the sum $\tilde{\lambda}$ must be 
positive (Proposition~\ref{bolibr} and (\ref{sumEVszero}) -- the sum of the 
eigenvalues of the union of all other diagonal blocks must be $\leq 0$ and 
it cannot be 0 because the eigenvalues 
$\lambda _{k,j}$ are generic). Hence, the right lower block is of size 1. 

$5^0$. Denote by $\Pi$ the left upper $(n-1)\times (n-1)$-block. Conjugate 
it to make all non-zero rows of the restriction of 
the $(p+1)$-tuple $\tilde{M}$ of matrices $M_j-I$ to $\Pi$ linearly 
independent. After the conjugation some of the rows of the restriction of 
$\tilde{M}$ to $\Pi$ might be 0. In this case conjugate the matrices $M_j$ 
by one and the same permutation matrix which places the zero rows of 
$M_j-I$ in the last (say, $m$) positions (recall that the last row of 
$M_j-I$ is 0, see $4^0$, so $m\geq 1$). Notice that if the restriction 
to $\Pi$ of a row of  
$M_j-I$ is zero, then its last (i.e. $n$-th) position is 0 as well, 
otherwise $M_j$ is not diagonalizable. 
 
$6^0$. Show that $m\geq n-n_1$ (and this will be the end of the proof of the 
lemma). One has $M_j=\left( \begin{array}{cc}G_j&R_j\\0&I\end{array}\right)$, 
$I\in GL(m,{\bf C})$. 

Denote by $\tilde{G}$ the representation defined by the matrices $G_j$. We 
regard the columns of 
the $(p+1)$-tuple of matrices $R_j$ as elements of the space 
${\bf C}'''(\tilde{G})$ defined as follows.

Each column of the $(p+1)$-tuple of matrices $R_j$ 
belongs to a linear space ${\bf C}'(\tilde{G})$ of dimension 
$\theta =r_1+\ldots +r_{p+1}$ 
which is the sum of the dimensions of the images of the linear operators 
$(.)\mapsto (G_j-I)(.)$ acting on ${\bf C}^{n-m}$ (every column of $R_j$ 
belongs to the image of this operator, otherwise $M_j$ will not be 
diagonalizable). Equality (\ref{M_j}) is equivalent to $n-m$ 
linear equations which the entries of the column must satisfy (for the 
block $R$ this equality implies 
$\sum _{j=1}^{p+1}G_1\ldots G_{j-1}R_j=0$; we prove in $7^0$ that these 
$n-m$ linear equations are linearly independent). Hence, this equality 
defines a subspace ${\bf C}''(\tilde{G})$ of ${\bf C}'(\tilde{G})$ 
of dimension $\theta -(n-m)$. 

One then factorizes ${\bf C}''(\tilde{G})$ by the space of $(p+1)$-tuples  
of blocks $(G_j-I)V$, $V\in {\bf C}^{n-m}$. These blocks are obtained as 
$R$-blocks when the $(p+1)$-tuple of matrices 
$\left( \begin{array}{cc}G_j&0\\0&1\end{array}\right)$ is conjugated by 
the matrix $V^*=\left( \begin{array}{cc}I&V\\0&1\end{array}\right)$. 
This factorization gives the space ${\bf C}'''(\tilde{G})$.

$7^0$. {\em The space ${\bf C}''(\tilde{G})$ is of codimension $n-m$ in 
${\bf C}'(\tilde{G})$.} 

One has to show that the $n-m$ linear relations defining 
${\bf C}''(\tilde{G})$ are linearly independent. If they are not, then the 
images of all linear operators 
$(.)\mapsto (G_j-I)(.)$ (acting on ${\bf C}^{n-m}$) must be contained in a 
proper subspace of ${\bf C}^{n-m}$ (say, the one defined by the first 
$n-m-1$ vectors of its standard basis). This means that all entries of the 
last rows of the matrices $G_j-I$ are 0. The matrices $M_j$ being 
diagonalizable, this implies that the entire $(n-m)$-th rows of $M_j-I$ are 
0. This contradicts the condition the first $n-m$ rows of the 
restriction to $\Pi$ of the $(p+1)$-tuple of matrices $M_j-I$ to be 
linearly independent, see $5^0$.

{\em The space ${\bf C}'''(\tilde{G})$ is of codimension 
$n-m$ in ${\bf C}''(\tilde{G})$, i.e. of dimension $\theta -2(n-m)$.}

Indeed, each column of $V$ 
belongs to ${\bf C}^{n-m}$ and the intersection ${\cal I}$ of the kernels of 
the operators $(.)\mapsto (G_j-I)(.)$ (acting on ${\bf C}^{n-m}$) is 
$\{ 0\}$, otherwise the matrices $M_j$ 
would have a non-trivial common centralizer. Indeed, if 
${\cal I}\neq \{ 0\}$, then after a change of the 
basis of ${\bf C}^{n-m}$ one can assume that a non-zero vector from 
${\cal I}$ equals 
$^t(1,0,\ldots ,0)$. Hence, the matrices $G_j$ are of the form 
$\left( \begin{array}{cc}1&\ast \\0&G_j^*\end{array}\right)$, 
$G_j^*\in GL(n-m-1,{\bf C})$, and one checks 
directly that $[M_j,E_{1,n}]=0$ for $E_{1,n}=\{ \delta _{i-1,n-j}\}$. 

$8^0$. The columns of the $(p+1)$-tuple of matrices $R_j$ (regarded 
as elements of ${\bf C}'''(\tilde{G})$) must be linearly independent, 
otherwise the monodromy group can be conjugated by a matrix 
$\left( \begin{array}{cc}I&\ast \\0&P\end{array}\right)$, 
$P\in GL(m,{\bf C})$, to a block-diagonal 
form, the right lower block (of size 1) 
for each $j$ being equal to 1 which means that 
the monodromy group is a direct sum and, hence, its centralizer is 
non-trivial -- a contradiction with Lemma~\ref{reprtrivcentr}. 

This means that $\dim {\bf C}'''(\tilde{G})\geq m$, i.e.  

\[ \theta -2(n-m)=r_1+\ldots +r_{p+1}-2(n-m)\geq m\]
which is equivalent to $m\geq n-n_1$; recall that 
$n_1=r_1+\ldots +r_{p+1}-n$. In the case of equality (and only in it) the 
columns of the $(p+1)$-tuple of matrices $R_j$ are a basis of the space 
${\bf C}'''(\tilde{G})$.

The lemma is proved.

\subsection{Proof of Lemma \ref{subreprtrivcentr}} 

$1^0$. If the lemma is not true, then $Z(\Phi )$ either 
contains a diagonalizable matrix $D$ with exactly two distinct eigenvalues 
or it contains a nilpotent matrix $N\neq 0$, $N^2=0$, 
see $1^0$ of the proof of Lemma~\ref{reprtrivcentr}.

$2^0$. In the first case one can conjugate the monodromy group to the form 
$\left( \begin{array}{cc}G_j&R_j\\0&I\end{array}\right)$ with 
$G_j=\left( \begin{array}{cc}M_j'&0\\0&M_j''\end{array}\right)$ 
where the sizes of $M_j'$, $M_j''$ equal the multiplicities of the two 
eigenvalues of $D$. One has 
$D=\left( \begin{array}{cc}\alpha I&0\\0&\beta I\end{array}\right)$, 
$\alpha \neq \beta$, $D$ and $G_j$ are $n_1\times n_1$. 

$3^0$. Denote by $\tilde{M}^{(i)}$ the $(p+1)$-tuple of matrices $M_j^{(i)}$, 
$i=1,2$. 
The space ${\bf C}'''(\tilde{G})$ was defined in $6^0$ -- $7^0$ of the proof 
of Lemma~\ref{theMGisreducible}. Set $k'=\dim {\bf C}'''(\tilde{M}')$, 
$k''=\dim {\bf C}'''(\tilde{M}'')$. Hence, $k'+k''=n-n_1$. Indeed, 
it follows from the definition of ${\bf C}'''(\tilde{G})$ that 

\[ {\bf C}'''(\tilde{G})=
{\bf C}'''(\tilde{M}')\oplus {\bf C}'''(\tilde{M}'')\]

$4^0$. Hence, there exists a conjugation of the monodromy 
group with a matrix $Q^*=\left( \begin{array}{cc}I&0\\
0&Q\end{array}\right)$, $Q\in GL(n-n_1,{\bf C})$ after which it is of the 
form $M_j=\left( \begin{array}{cccc}M_j'&0&F_j'&0\\0&M_j''&0&F_j''\\
0&0&I&0\\0&0&0&I\end{array}\right)$. Here the columns of the 
$(p+1)$-tuples of matrices $F'$, $F''$ are bases of the spaces 
${\bf C}'''(\tilde{M}')$, ${\bf C}'''(\tilde{M}'')$. (One first obtains the 
form 
$\left( \begin{array}{cc}F_j'&0\\\ast &F_j''\end{array}\right)$ of the block 
$R_j$ by such a conjugation; after this by conjugation with another matrix  
of the same form as $Q^*$ one makes the block $*$ equal to 0.)

This shows that the monodromy group 
is a direct sum (one has to perform a self-evident permutation of the 
blocks and columns to make the matrices block-diagonal). Hence, its 
centralizer is non-trivial which contradicts Lemma~\ref{reprtrivcentr}.

$5^0$. If there exists a matrix $0\neq N\in Z(\Phi )$, $N^2=0$, then one 
can conjugate $N$ to the form 
$N=\left( \begin{array}{ccc}0&0&I\\0&0&0\\0&0&0\end{array}\right)$ (or 
$N=\left( \begin{array}{cc}0&I\\0&0\end{array}\right)$) with $I$ being 
$q\times q$, $q\leq n_1/2$; the second case corresponds to $q=n_1/2$. 

$6^0$. Hence, 
$M_j=\left( \begin{array}{cccc}M_j'&R_j&T_j&L_j\\0&M_j''&S_j&H_j\\
0&0&M_j'&P_j\\0&0&0&I\end{array}\right)$ where $M_j'$ is $q\times q$ and if 
$q=n_1/2$, then the blocks $R_j$, $M_j''$, $S_j$ and $H_j$ are absent. 
In a similar way one brings the blocks 
$W_j=\left( \begin{array}{c}L_j\\H_j\\P_j\end{array}\right)$ to the form 
$\left( \begin{array}{ccc}F_j'&\ast &\ast \\0&F_j''&\ast \\0&0&F_j'
\end{array}\right)$ with the same meaning of $M_j^{(i)}$ and $F_j^{(i)}$ as 
in $3^0$ (such a form of $W_j$ can be achieved by conjugation with a matrix 
$Q^*$, see $4^0$).

$7^0$. By permuting the rows and columns of $M_j$ (which results from a 
conjugation) one brings $M_j$ to the form indicated below, with the matrix 
$Z$ belonging to the centralizer of the monodromy group which again 
contradicts Lemma~\ref{reprtrivcentr}.

\[ M_j=\left( \begin{array}{cccccc}M_j'&F_j'&R_j&\ast &T_j&\ast \\
0&I&0&0&0&0\\0&0&M_j''&F_j''&S_j&\ast \\
0&0&0&I&0&0\\0&0&0&0&M_j'&F_j'\\0&0&0&0&0&I\end{array}\right) ~,~
Z=\left( \begin{array}{cccccc}0&0&0&0&I&0\\
0&0&0&0&0&I\\0&0&0&0&0&0\\0&0&0&0&0&0\\0&0&0&0&0&0\\0&0&0&0&0&0
\end{array}\right) \]

The lemma is proved. 

\section{On adjacency of nilpotent orbits\protect\label{adjac}}

Denote by $J_1$ and $J_2$ two nilpotent orbits (i.e. conjugacy classes). Let 
$D_1\in J_1$, $D_2\in J_2$. Denote by $\rho _i$ and $\theta _i$ the ranks of 
the matrices $(D_1)^i$ and $(D_2)^i$. It is known that the orbit $J_1$ 
belongs to the closure of the orbit $J_2$ if and only if 
one has (for all $i\in {\bf N}$) $\rho _i\leq \theta _i$ 
(see \cite{Kr}, p. 21). 

In this section we prove a more concrete statement, see 
Theorem~\ref{adjacency1}, from which we deduce Corollary~\ref{adjacency2}. 
The latter is used in the proof of Theorem~\ref{nongenericevs}.

It is evident that if $J_1$ belongs to the closure of $J_2$, then 
one has for all $i\in {\bf N}$ $\rho _i\leq \theta _i$. 
To prove the implication in the 
other direction we use the following {\em operation $(s,l)$}, defined 
for $s\geq l$, $s,l\in {\bf N}^*$: suppose that the nilpotent orbit $J$ has 
two Jordan blocks, of sizes $s$ and $l$, $s\geq l$. We say that the nilpotent 
orbit $J'$ (of the same size as $J$) is obtained from $J$ with the help of 
the operation $(s,l)$ if $J'$ has all Jordan blocks of the same sizes as 
$J$ except these two which 
are replaced by two blocks of sizes $s+1$ and $l-1$. If $l-1=0$, then only 
one block of size $s+1$ replaces the two blocks of sizes $s$ and 1. 

\begin{prop}\label{sl}
The orbit $J$ lies in the closure of the orbit $J'$.
\end{prop} 

{\em Proof:} 

$1^0$. Assume first that the orbit $J$ is of size $s+l$. Consider the 
matrix 
$U(\varepsilon )=\left( \begin{array}{cr}S&\varepsilon R\\0&L
\end{array}\right)$ 
where $S$ and $L$ are 
upper-triangular nilpotent Jordan blocks of sizes $s$ and $l$ and 
$\varepsilon \in {\bf C}$. The block $R$ contains a single unit in its 
lowest row and last column, its other entries are zeros. 

$2^0$. For $\varepsilon =0$ one has 
$U\in J$. For $\varepsilon \neq 0$ the matrix $U$ belongs to $J'$, i.e. it 
is conjugate to a nilpotent Jordan matrix with two  
Jordan blocks, of sizes $l+1$ and $s-1$. Indeed, the number and sizes of 
Jordan blocks of a nilpotent orbit $\Omega$ are defined by the ranks of the 
matrices $A^i$, $A\in \Omega$. These ranks computed for 
$U|_{\varepsilon \neq 0}$ and for $V\in J'$ coincide 
(to be checked directly). 

$3^0$. It is obvious that the matrix $U|_{\varepsilon =0}$ lies in the 
closure of 
the orbit of any of the matrices $U|_{\varepsilon =\varepsilon _0}$, 
$\varepsilon _0\neq 0$ (which is 
$J'$, i.e. one and the same for all $\varepsilon _0\neq 0$). This follows from 
the inclusion of $U(0)$ in the family $U(\varepsilon )$.  

$4^0$. If the size of the matrices is $>s+l$, then the proposition is proved 
by 
analogy (one sets $J=J_0\oplus J_*$, $J'=J_0'\oplus J_*$ where $J_0$, $J_0'$ 
are nilpotent orbits of size $l+s$, with blocks of sizes $s$, $l$ and 
$s+1$, $l-1$ and $J_*$ is some nilpotent orbit). 

The proposition is proved. $\hspace{2cm}\Box$

\begin{tm}\label{adjacency1}
The two orbits $J_1$ and $J_2$ can be connected by a chain of 
intermediate orbits such that each orbit of the chain is obtained from the 
previous one by some operation $(s,l)$ and, hence, each previous orbit lies 
in the closure of the next one.
\end{tm}

{\em Proof of Theorem \ref{adjacency1}:} 

$1^0$. Assume that for the two nilpotent orbits $J_1$, $J_2$ there holds 
$\rho _i\leq \theta _i$ for all $i$, with strict inequality for at least 
one $i$. If each of the orbits contains a Jordan block of size $k$, then 
one can decrease the size of the orbits by $k$ by excluding the two 
equal blocks from consideration. 
So assume that the two orbits have no such couple of blocks; in particular, 
that for the blocks of greatest size $H_1\in J_1$ and $H_2\in J_2$ one has 
$h_1=$size$(H_1)<h_2=$size$(H_2)$. Indeed, $h_1>h_2$ would imply that 
$\rho _{h_1-1}>\theta _{h_1-1}$. 

$2^0$. Denote the chain of nilpotent orbits joining $J_1$ with $J_2$ by 
($J_1$, $J_3$,$\ldots$, $J_{\nu }$,$J_2$). 
Denote by $h$ the size of the second largest Jordan block of $J_1$. Then the 
operation $(h_1,h)$ applied to $J_1$ preserves the quantities $\rho _i$ 
for $i\leq h_1-1$ and increases $\rho _{h_1}$ by 1 (it changes from 0 to 1).

$3^0$. Define $J_3$ as obtained from $J_1$ by the 
operation $(h_1,h)$. If $J_3$ coincides with 
$J_2$, then the construction of the chain is finished. If not, then we  
construct $J_4$ after $J_3$ in the same way as $J_3$ 
was constructed 
after $J_1$. Namely,  
if the greatest of the sizes of the Jordan blocks of $J_3$, i.e. $h_1+1$, 
equals $h_2$, then one can exclude the greatest blocks 
of $J_2$ and $J_3$ from consideration 
and continue in the same way with orbits of smaller size. This means that 
the block of size $h_2$ will be present in all orbits $J_3$, $J_4$, $\ldots$, 
$J_{\nu}$, $J_2$. If $h_1+1<h_2$, then 
one can repeat what was done in $1^0$ -- $2^0$ with $J_3$ on the place of 
$J_1$ etc. 

$4^0$. After finitely many such steps one will have $\rho _i=\theta _i$ for 
all $i$, i.e. one obtains the orbit $J_2$. Each orbit of 
the chain is obtained from the previous one by some 
operation $(s,l)$. Each previous orbit lies in the closure of the next one, 
see Proposition~\ref{sl}, so $J_1$ lies in the closure of $J_2$. 

The theorem is proved.$\hspace{2cm}\Box$\\

If the nilpotent orbit $J_1$ belongs to the closure of the nilpotent orbit 
$J_2$ (both of size $g$), then in general $J_1$ might have more Jordan blocks 
than $J_2$, i.e. rk$(J_1)<$rk$(J_2)$. In this case we assume that $J_2$ has 
rk$(J_2)-$rk$(J_1)$ Jordan blocks of size 0, so that both orbits have the 
same number of Jordan blocks. When the numbers of Jordan blocks are defined 
in this way, one can add one and the same number of Jordan blocks of size 
0 to $J_1$ and $J_2$. In what follows we assume that the number of Jordan 
blocks of size 0 is known.

\begin{cor}\label{adjacency2}
Increase by 1 the sizes of the $k$ smallest blocks 
of $J_1$ and of the $k$ smallest blocks of $J_2$ -- this defines two 
nilpotent orbits $J_1'$, $J_2'$, both of size $g+k$. Then $J_1'$ lies in the 
closure of $J_2'$.
\end{cor} 

\begin{rem}\label{onlypartof}
It might happen that $k$ is such that one has to increase 
by 1 only part of the Jordan blocks of a given size. Example: there are 
4 Jordan blocks in $J_1$, of sizes 2, 2, 1, 1 and $k=3$. Then the new 
sizes are 3, 2, 2, 2, i.e. only one of the two blocks of size 2 becomes of 
size 3.
\end{rem}  

{\em Proof:}

$1^0$. Consider first the case when $J_2$ is obtained from $J_1$ by an 
operation $(s,l)$. 

Denote by $m$ the greatest of the sizes of the $k$ Jordan blocks of $J_1$ 
to be increased by 1. 

{\em Case 1)} One has $m\neq l$ and $m\neq s$. Hence, $J_2'$ is 
obtained from $J_1'$ by an operation 

$(s,l)$ if $m<l$, 

$(s,l+1)$ if $l<m<s$,

$(s+1,l+1)$ if $s<m$.

{\em Case 2)} One has $m=l<s$. Hence, $J_2'$ is 
obtained from $J_1'$ by an operation $(s,l+1)$. 

{\em Case 3)} One has $m=l=s$. In this case  
either $J_2'=J_1'$ or $J_2'$ is 
obtained from $J_1'$ by an operation $(l+1,l+1)$. The first (resp. the second) 
possibility 
takes place when not all Jordan blocks of $J_1$ of size $l$ have to be 
chosen as smallest blocks and their sizes increased by 1 (resp. 
when all of them have to be chosen as such).

{\em Case 4)} One has $m=s>l$. Hence, $J_2'$ is 
obtained from $J_1'$ by an operation $(s,l+1)$ or $(s+1,l+1)$. The first 
(resp. the second) possibility takes place when not all (resp. when all) 
Jordan blocks 
of $J_1$ of size $s$ have to be chosen as smallest blocks and their sizes 
increased by 1. 

Hence, in all these cases $J_1'$ lies in the closure of $J_2'$ or coincides 
with it.

$2^0$. In the general case one uses Theorem \ref{adjacency1} and applies 
$1^0$ to each couple of consecutive orbits from the chain connecting 
$J_1$ and $J_2$.

The corollary is proved. $\hspace{2cm}\Box$

\section{Proof of Theorem \protect\ref{nongenericevs}\protect\label{n_s>1}}
\subsection{The basic lemma and its corollaries}

It is clear that conditions {\em i)} and {\em ii)} from 
Theorem~\ref{generic} 
are necessary for the existence of $(p+1)$-tuples of matrices with trivial 
centralizers and with $d=1$ -- the basic technical tool 
allows one to deform such a $(p+1)$-tuple ${\cal A}_0$ (resp. ${\cal M}_0$)
with a trivial centralizer 
into a nearby irreducible one ${\cal A}_1$ (resp. ${\cal M}_1$) 
of matrices from the corresponding 
diagonal Jordan normal forms and with generic eigenvalues. The condition 
$d=1$ implies that the PMV of the eigenvalues of the $(p+1)$-tuple 
${\cal A}_1$ (resp. ${\cal M}_1$) is simple, hence, generic eigenvalues with 
such PMVs exist. 
Conditions 
{\em i)} and {\em ii)} from Theorem~\ref{generic} must hold for 
${\cal A}_1$ (resp. for ${\cal M}_1$), hence, they hold for 
${\cal A}_0$ (resp. for ${\cal M}_0$) as well, see Remarks~\ref{psiinvariant}. 

Therefore we prove only the sufficiency of conditions {\em i)} and {\em ii)} 
for the existence of $(p+1)$-tuples of matrices with trivial 
centralizers. 

Recall that the integers $n_i$, $i=0,1,\ldots ,s$ were defined before 
Theorem~\ref{generic}   
and that the conditions 
$n_s>1$ and $(\alpha _n)$ being a strict inequality are equivalent, see 
Corollary~\ref{alsogood1} and Remark~\ref{sameindexofrigidity}.

\begin{lm}\label{trivialcentr}
If for the Jordan normal forms $J_j^n$ each with a single eigenvalue, 
with $d=1$ 
and with $n_s>1$ conditions {\em i)} and {\em ii)} from Theorem~\ref{generic} 
hold, then the DSP is weakly solvable for nilpotent 
matrices $A_j^0$ with Jordan normal forms $J_j^n$.
\end{lm}

\begin{rem}
The lemma is true also in the case when $d>1$ and the index of rigidity 
$\kappa$ of 
the $(p+1)$-tuple of Jordan normal forms is strictly negative. The proof is 
the same with the exception of $1^0$ of it where the four exceptional cases 
are eliminated due to $d>1$ and $\kappa <0$.
\end{rem}

The lemma is proved in the next subsection. It implies the following 
corollary which finishes the proof of the theorem. 

\begin{cor}\label{suff}
1) For any $(p+1)$-tuple of Jordan normal forms $J_{j;0}^n$  
with $d=1$, with $n_s>1$ and satisfying conditions {\em i)} and {\em ii)} 
from Theorem~\ref{generic} the DSP is weakly solvable for matrices $A_j$. 
In particular, for any generic eigenvalues and such Jordan normal forms 
it is solvable for matrices $A_j$.

2) For any $(p+1)$-tuple of Jordan normal forms $J_{j;0}^n$ as in 1) and 
for any generic eigenvalues whose product is 1   
the DSP is solvable for matrices $M_j$.
\end{cor}

{\em Proof:}

$1^0$. Set $A_j^0=Q_j^{-1}D'_jQ_j$ where $A_j^0$ are nilpotent and their 
$(p+1)$-tuple is with trivial centralizer. Look for matrices $A_j$ of the form 
$A_j=(I+\varepsilon X_j(\varepsilon ))^{-1}Q_j^{-1}(D'_j+\varepsilon D_j)
Q_j(I+\varepsilon X_j(\varepsilon ))$. Here the matrices $D_j'$ and $D_j$ 
have the same meaning as $D'$ and $\Delta _s$ from Proposition~\ref{DD'}. 

$2^0$. The basic technical tool provides the existence of such matrices 
$A_j$ for $\varepsilon$ small enough. Hence, the matrix $A_j$ is 
conjugate to $\varepsilon D_j$ (Proposition~\ref{DD'}). One can multiply 
such a $(p+1)$-tuple by $1/\varepsilon$ and, hence, in the new $(p+1)$-tuple 
the matrix $\varepsilon ^{-1}A_j$ will be conjugate to $D_j$. As the 
matrices $D_j$ have Jordan normal forms $J_j^n$ and can have any 
eigenvalues, this proves the corollary for matrices $A_j$.

$3^0$. Prove part 2) using part 1) already proved. One needs to 
consider the monodromy operators $M_j$ of the fuchsian system (\ref{Fuchs}).  
The matrices $A_j$ are chosen from the diagonal 
Jordan normal forms corresponding to 
$J(M_j)$. The eigenvalues of each matrix $A_j$ form a normalized 
chain. The condition $d=1$ implies that the PMV of the eigenvalues 
of the matrices $A_j$ is simple and one can find 
strongly generic eigenvalues $\lambda _{k,j}$. The eigenvalues of each 
matrix $A_j$ form a normalized chain. Hence, 
for each $j$ one will have $J(M_j)=J_{j;0}^n$ or $J(M_j)$ will be 
subordinate to $J_{j;0}^n$, see Lemma~\ref{JNFofM}.

The strong genericity of the eigenvalues implies that the monodromy group 
is irreducible. It follows from Lemma~\ref{J^0impliesJ^1} that one can 
construct an irreducible monodromy group with $J(M_j)=J_{j;0}^n$ for all $j$. 

The corollary is proved. $\hspace{2cm}\Box$

\subsection{Proof of Lemma \protect\ref{trivialcentr}}

{\bf The case $s=0$.}

$1^0$. If condition $(\omega _n)$ holds, and if $s=0$ (hence, $n=n_s>1$), 
then there exist nice $(p+1)$-tuples of nilpotent matrices $A_j$ 
from these nilpotent orbits, see Theorem~\ref{nilpunip}. Indeed, one is 
never in one of the four 
exceptional cases cited in  Theorem~\ref{nilpunip} due to 
$d=1$. This proves the lemma in the case $s=0$.

The condition $d=1$ implies (for any possible value of $s$, not only 
for $s=0$) that the Jordan normal forms $J_j^{n_s}$ never correspond 
to one of these four exceptional cases. Indeed, in these four cases 
the PMV of the diagonal Jordan normal forms corresponding to $J_j^{n_s}$ 
is non-simple. By Lemma~\ref{simple}, the PMV of 
the diagonal Jordan normal forms corresponding to $J_j$ is non-simple. 
This implies that $d>1$ -- a contradiction.

{\bf The case $s=1$.}

$2^0$. If $s=1$ and if $(\omega _{n_1})$ holds for the Jordan normal 
forms $J_j^{n_1}$, then 
there exists a $(p+1)$-tuple 
of nilpotent matrices $A_j^0$ (satisfying (\ref{A_j})) blocked as follows:
$A_j^0=\left( \begin{array}{cc}A_j'&B_j\\0&0\end{array}\right)$ where 
$A_j'$ are $n_1\times n_1$, the 
$(p+1)$-tuple of conjugacy classes $C(A_j')$ defines the Jordan normal forms 
$J_j^{n_1}$ 
and satisfies condition 
$(\omega _{n_1})$. The $(p+1)$-tuple of matrices $A_j'$ is presumed 
nice. Show that one can choose the blocks $B_j$ such that 

1) for every $j$ the matrix $A_j^0$ is from the closure of the 
necessary conjugacy class;

2) the centralizer ${\cal Z}$ of the $(p+1)$-tuple of matrices 
$A_j^0$ is trivial.

$3^0$. Set

\[ \Delta =\{ (A_1'T,\ldots ,A_{p+1}'T)~~|~~T\in {\bf C}^{n_1}\} ~.\]
Hence, dim$\Delta =n_1$. Indeed, the $(p+1)$-tuple of 
matrices $A_j'$ being nice implies that the intersection of the kernels of the 
linear operators $T\mapsto A_j'T$ is $\{ 0\}$. 

$4^0$. Every column of the block $B_j$ belongs to a linear space ${\cal S}_j$. 
It can be described as follows. There exists $m\in {\bf N}$ such that all 
Jordan blocks 
of size $\geq m$ of $A_j'$ are Jordan blocks of $A_j^0$ of the same size 
and all Jordan blocks of size $<m$ have their size increased by 1 in $A_j^0$ 
(and, of course, $A_j^0$ may contain Jordan blocks of size 1 which are 
not present in $A_j'$). The matrices $A_j^0$ and $A_j'$ satisfy the 
conditions rk$(A_j^0)^i=$rk$(A_j')^i$ for $i\geq m$.

$5^0$. Denote by ${\cal Q}_j$ the space  
Im$(A_j')$. Hence, dim${\cal Q}_j=$rk$A_j'$. One can choose as 
columns of $B_j$ $r_j-$rk$A_j'$ vector-columns from Ker$(A_j')^{m-1}$  
linearly independent modulo ${\cal Q}_j$ so that the Jordan normal form of 
$A_j^0$ be the necessary one. Indeed, one can 
conjugate $A_j^0$ so that the block $A_j'$ be in upper-triangular Jordan 
normal form; the conjugation can be performed by a block-diagonal matrix, 
with diagonal blocks of sizes $n_1$, 1,$\ldots$,1. 
 
Conjugate after this the matrix $A_j^0$ by a block upper-triangular matrix, 
the diagonal blocks (of sizes $n_1$ and $n-n_1$) being equal to $I$, so that 
the block $B_j$ contain only zeros in the rows where the units of $A_j'$ are. 

After this conjugate $A_j^0$ by a block-diagonal matrix, with diagonal 
blocks of sizes $n_1$ and $n-n_1$, the first of them being equal to $I$, so 
that in the columns of the block $B$ there be exactly one non-zero entry 
equal to 1 in each of the rows described below (call them {\em marked}) and 
in different columns. If 
$A_j'$ has a Jordan block in the rows with indices $l$, $l+1$, $\ldots$, 
$l+s$ whose size has to be increased by 1 when it is considered as a block 
of $A_j^0$, then the $(l+s)$-th row is marked. 

(Notice that all these conjugations 
preserve the size $n_1\times (n-n_1)$ of the block $B_j$.)

One checks directly that the Jordan normal form of the matrix $A_j^0$ is the 
necessary one (which is easy to do in the present form of $A_j^0$). The space 
${\cal S}_j$ is the preimage (before the above conjugations) of the space 
spanned by the vector-columns having non-zero entries only in the marked rows.

From now on we presume that for every $j$ the $r_j-$rk$A_j'$ vector-columns 
from ${\cal S}_j$ (linearly independent modulo the subspace ${\cal Q}_j$) 
are fixed.

$6^0$. Denote by $\Omega$ the space of $(p+1)$-tuples of 
columns of the blocks $B_j$ modulo the space $\Delta$. Hence,  
dim$\Omega \geq r_1+\ldots +r_{p+1}-2n_1=n-n_1$. One subtracts $n_1$  
twice because the sum of 
the $p+1$ columns (of length $n_1$) must be 0 and to factor out $\Delta$. 
(We do not need to discuss the question when the inequality is strict and when 
it is an equality.) 

Choose $n-n_1$ $(p+1)$-tuples of columns of 
the block $B$ which are linearly independent modulo $\Delta$ and whose sum 
is 0. Check that there holds condition 2) from $2^0$. 

Let a matrix $X\in {\cal Z}$ equal 
$\left( \begin{array}{cc}Y&Z\\T&U\end{array}\right)$, $Y$ being 
$n_1\times n_1$. The commutation 
relations yield 

a) $TA_j'=0$ for all $j$, (hence, $TA=0$ for every matrix from the matrix 
algebra ${\cal A}$ generated by the matrices $A_j'$; their $(p+1)$-tuple 
being nice, one has $T=0$) 

b) $[Y,A_j']=0$ for all $j$, (hence, $Y=\alpha I$ because the $(p+1)$-tuple 
of matrices $A_j'$ is nice)

c) $A_j'Z+B_jU-YB_j=0$, i.e. $A_j'Z+B_jU-\alpha B_j=0$; as the columns of 
the $(p+1)$-tuple of blocks 
$B_j$ are independent modulo the space $\Delta$, i.e.  modulo 
columns of the form $A_j'Z$, one must have $Z=0$ and $U=\alpha I$. Hence, 
the centralizer ${\cal Z}$ is trivial.

$7^0$. Condition 1) from $2^0$ holds, see the construction of the spaces 
${\cal S}_j$ in $4^0$ and $5^0$. We admit that for some values of $j$ the 
conjugacy class of $A_j^0$ might 
be not the necessary one but from its closure for the following reason --  
when one constructs $n-n_1$ $(p+1)$-tuples of columns of 
the block $B$ which are linearly independent modulo $\Delta$ and whose sum 
is 0, one does not know whether for each $j$ the 
vector-columns of $B_j$ of this $(p+1)$-tuple 
span the {\em whole} space ${\cal S}_j/{\cal Q}_j$.  

{\bf The case of arbitrary $s$.}

$8^0$. For arbitrary $s$ one constructs the $(p+1)$-tuple of 
nilpotent matrices $A_j^0$ in a similar way: namely, block-decompose any 
$n\times n$-matrix, the diagonal blocks being of sizes $n_s$, $n_{s-1}-n_s$, 
$\ldots$, $n_0-n_1=n-n_1$. Call {\em basic minor of size $n_{s-k+1}$} 
of a given matrix $A$ (denoted by $A|_{L_k}$) the square submatrix which is 
the restriction of $A$ to the first 
$n_{s-k+1}$ rows and $n_{s-k+1}$ columns, $k=1,\ldots ,s+1$. Denote by 
$H_{\mu ,\nu}$ the blocks of a matrix from $gl(n,{\bf C})$ in this 
block-decomposition, 
$\mu ,\nu =1,\ldots ,s+1$ (we enumerate them in the usual way, from above to 
below and from left to right). 

Denote by ${\cal Z}_k$ the centralizer of the 
$(p+1)$-tuple of matrices $A_j^0|_{L_k}$. By abuse of language we denote by 
$J_j^{n_{s-k+1}}$ both the Jordan normal form with a single eigenvalue 
and the nilpotent orbit defining such a Jordan normal form.

$9^0$. The basic minors of the matrices $A_j^0$  
of sizes $n_s$ and $n_{s-1}$ are constructed  
like in the cases $s=0$ and $s=1$, see $1^0$ -- $7^0$. After this the basic 
minors of sizes $n_{s-k}$ for 
$k<s+1$ are constructed like the one of size $n_{s-1}$ but one 
has to take into account the possibility $A_j^0|_{L_k}$ to be from an 
orbit subordinate to the required one. 

Denote by $v$ and $w$ the ranks of $A_j^0|_{L_k}$ and of the required orbit 
$J_j^{n_{s-k+1}}$ of $A_j^0|_{L_k}$ (i.e. $w=r(J_j^{n_{s-k+1}})$). Hence, 
$v\leq w$. Define 
the space ${\cal S}_{j,k}\subset {\bf C}^{n_{s-k+1}}$ (the analog of 
${\cal S}_j$, see 
$4^0$ -- $5^0$). Suppose that $A_j^0|_{L_k}$ is in upper-triangular 
Jordan normal form (the conjugation of $A_j^0$ 
can be carried out by a block-diagonal matrix, with diagonal blocks of 
sizes $n_{s-k+1}$, 1, $\ldots$, 1). 

$10^0$. Define as {\em marked} the rows which are last rows of the $b$ 
smallest blocks of $A_j^0|_{L_k}$ where 

\[ b=r(J_j^{n_{s-k}})-v.\] 
If $A_j^0|_{L_k}$ belongs to the closure of the orbit $J_j^{n_{s-k+1}}$ but 
not to itself, then one might have to choose among the blocks of a given size 
the size of which ones to be increased by 1, see Remark~\ref{onlypartof}. 
The space ${\cal S}_{j,k}^0$ is defined as spanned by all 
vector-columns in ${\bf C}^{n_{s-k+1}}$ which have units in the rows where 
$A_j^0|_{L_k}$ has a unit or a marked row. The space ${\cal S}_{j,k}$ is the 
preimage of ${\cal S}_{j,k}^0$ under the conjugation of $A_j^0|_{L_k}$.  

Hence, dim${\cal S}_{j,k}=r(J_j^{n_{s-k}})$. 
 
Denote by $D_{k+1}$ the union of the blocks $H_{k,k+1}$, 
$H_{k-1,k+1}$, $\ldots$, $H_{1,k+1}$. Denote by $\Delta _k$ 
the space of $(p+1)$-tuples of the form 

\[ (A_1^0|_{L_k}T,\ldots ,A_{p+1}^0|_{L_k}T)~~,~~
T\in {\bf C}^{n_{s-k+1}}\]
We make two inductive assumptions:
  
1) for all $j$ the orbit of $A_j^0|_{L_k}$ to be either 
$J_j^{n_{s-k+1}}$ or one subordinate to it; 

2) the columns of the $(p+1)$-tuple of matrices $A_j^0|_{L_k}$ to be 
linearly independent. (This means that dim$(\Delta _k)=n_{s-k+1}$. The 
assumption is true for $k=0$ because the algebra ${\cal A}$ is nice, 
and for $k=1$ by the construction of $\Omega$, see $6^0$.) 

$11^0$. {\em If the columns of $A_j^0|_{D_{k+1}}$ belong to the 
space ${\cal S}_{j,k}$, then the orbit of $A_j^0|_{L_{k+1}}$ will be either 
$J_j^{n_{s-k}}$ or one subordinate to $J_j^{n_{s-k}}$. This follows 
from Corollary~\ref{adjacency2}.} 

Indeed, two cases are possible: 

{\em Case 1)} One has $n_{s-k}-n_{s-k+1}\geq b$. 

If the columns of $A_j^0|_{D_{k+1}}$ belong to the 
space ${\cal S}_{j,k}$, then the orbit of 
$A_j^0|_{L_{k+1}}$ is either $J_0$ (obtained from 
the one of $A_j^0|_{L_k}$ by increasing by 1 the sizes of the 
$n_{s-k}-n_{s-k+1}$ smallest Jordan blocks of $A_j^0|_{L_k}$) or belongs to 
the closure of $J_0$. Notice that 
if $n_{s-k}-n_{s-k+1}>b$, then we add $n_{s-k}-n_{s-k+1}-b$ Jordan blocks of 
size 0 to the orbit of $A_j^0|_{L_k}$.  

On the other hand, $J_j^{n_{s-k}}$ is obtained from $J_j^{n_{s-k+1}}$ by 
increasing by 1 the sizes of the $n_{s-k}-n_{s-k+1}$ smallest blocks 
(one adds  
$y:=n_{s-k}-n_{s-k+1}-(r(J_j^{n_{s-k}})-r(J_j^{n_{s-k+1}}))$ 
Jordan blocks of size 0 to $J_j^{n_{s-k+1}}$ if $y>0$; recall 
that one always has $y\geq 0$ -- this follows from the definition of the 
Jordan normal forms $J_j^{n_i}$). 

Hence, $J(A_j^0|_{L_{k+1}})$ is either $J_j^{n_{s-k}}$ or is subordinate to it 
(this follows from assumption 1) from $10^0$ and from 
Corollary~\ref{adjacency2}).)

{\em Case 2)} One has $n_{s-k}-n_{s-k+1}<b$.

If the columns of $A_j^0|_{D_{k+1}}$ belong to 
the space ${\cal S}_{j,k}$, then again the orbit of $A_j^0|_{L_{k+1}}$ 
belongs to the closure of $J_j^{n_{s-k}}$. 

Indeed, increase the sizes of 
$J_j^{n_{s-k+1}}$ and of $A_j^0|_{L_k}$ not by $n_{s-k}-n_{s-k+1}$ but by 
$b$, and increase by 1 the sizes of their $b$ smallest Jordan blocks (one 
adds the necessary number of Jordan blocks of size 0 to $J_j^{n_{s-k+1}}$). 

Denote the matrices thus obtained respectively by $L$ and $P$. The last 
$b-n_{s-k}+n_{s-k+1}$ Jordan blocks of $L$ are of size 1. By 
Corollary~\ref{adjacency2}, $P$ belongs to the closure of the orbit of $L$ 
(because $A_j^0|_{L_k}$ belongs to the closure of $J_j^{n_{s-k+1}}$). 
Hence, for all $i$ one has rk$((P)^i)\leq$rk$((L)^i)$. 
When one reduces the 
sizes of $L$ and $P$ by deleting their last $b-n_{s-k}+n_{s-k+1}$ columns 
and rows,  
the quantities rk$((L)^i)$ do not change because we delete Jordan blocks of 
size 1 while the quantities rk$(P^i)$ decrease or remain the same (recall 
that all entries of the last $b$ rows of the matrices $P$ and $L$ are 0). 

Hence, for all $i$ one has 
rk$((J_j^{n_{s-k}})^i)\geq$rk$((A_j^0|_{L_{k+1}})^i)$ and 
the matrix $A_j^0|_{L_{k+1}}$ belongs to the 
closure of the nilpotent orbit $J_j^{n_{s-k}}$.

$12^0$. {\em The columns of the $(p+1)$-tuple of blocks 
$A_j^0|_{D_{k+1}}$ 
(whose sum is 0) can be chosen linearly independent modulo the 
space $\Delta _k$.}
 
This is proved by complete analogy with the case $s=1$, by estimating the 
dimension of 
the linear space to which these $(p+1)$-tuples of columns belong. Namely, 
the dimension of the space of such $(p+1)$-tuples of columns is 

\begin{equation}\label{dim>=}
 \geq z:=r(J_1^{n_{s-k}})+\ldots  
+r(J_{p+1}^{n_{s-k}})-{\rm dim}(\Delta _k)-n_{s-k+1} 
=n_{s-k}-n_{s-k+1}
\end{equation} 
because dim$(\Delta _k)=n_{s-k+1}$ and one has 
$r(J_1^{n_{s-k}})+$ $\ldots$ 
$+r(J_{p+1}^{n_{s-k}})=n_{s-k}+n_{s-k+1}$. Subtracting $n_{s-k+1}$ in 
(\ref{dim>=}) corresponds to imposing the $n_{s-k+1}$ conditions the sum of 
the $(p+1)$ vector-columns to be 0; we do not prove that these conditions 
are independent, therefore we claim only that 
the dimension is $\geq z$, not necessarily equal to $z$. 

Hence, the $n_{s-k}$ $(p+1)$-tuples of columns of the $(p+1)$-tuple of 
matrices $A_j^0|_{L_{k+1}}$ are linearly independent. For $k=s+1$ this 
implies the linear independence of the $(p+1)$-tuples of columns of the 
$(p+1)$-tuple of matrices $A_j^0$. 

$13^0$. {\em The centralizer ${\cal Z}_{s+1}$ is 
trivial.}

This is proved like it was done for $s=1$. Namely, give the decomposition 
of a matrix from ${\cal Z}_{s+1}$ in blocks $H_{\mu ,\nu}$. 

Consider the blocks $H_{k,1}$ of a 
matrix $[X,A_j^0]$, $X\in {\cal Z}_{s+1}$. One obtains consecutively 
$X|_{H_{s+1,1}}=\ldots =X|_{H_{2,1}}=0$, $X|_{H_{1,1}}=\alpha I$. 
Indeed, the $(p+1)$-tuple of matrices $A_j^0|_{H_{1,1}}$ is nice and 
the equality $[X,A_j^0]|_{H_{s+1,1}}=0$ is equivalent to 
$(XA_j^0)|_{H_{s+1,1}}=0$ (because $A_j^0|_{H_{s+1,k}}=0$ for all $k$, $j$); 
the latter implies $X|_{H_{s+1,1}}=0$ like in a) from $6^0$. Then similarly 
one 
deduces that $X|_{H_{s,1}}=0$ (making use of $X|_{H_{s+1,1}}=0$) and that 
$X|_{H_{s-1,1}}=\ldots =X|_{H_{2,1}}=0$. The 
equality $X|_{H_{1,1}}=\alpha I$ 
follows from the $(p+1)$-tuple of matrices $A_j^0|_{H_{1,1}}$ being nice 
like in b) from $6^0$. 

Assume that $\alpha =0$. 

$14^0$. Then consider $[X,A_j^0]|_{H_{k,2}}$ for $k=1,\ldots ,s+1$. 
Hence, $X|_{H_{k,2}}=0$ for $k=1,\ldots ,s+1$, 
otherwise the columns of the $(p+1)$-tuple of matrices $A_j^0$ will be 
linearly dependent. Notice that $[X,A_j^0]|_{H_{k,2}}=0$ is equivalent to 
$(A_j^0X)|_{H_{k,2}}=0$ because one has $(XA_j^0)|_{H_{k,2}}=0$ after $13^0$. 

Then consider in the same way $[X,A_j^0]|_{H_{k,3}}$, $[X,A_j^0]|_{H_{k,4}}$ 
etc. One obtains in a similar way that the restrictions of $X$ to 
$H_{k,3}$, $H_{k,4}$ etc. are 0, i.e. $X=0$. 
Without the assumption $\alpha =0$ this 
would mean that $X=\alpha I$. Hence, the centralizer is trivial. 

$15^0$. If for some $j$ the orbit of the matrix $A_j^0$   
is subordinate to the necessary one, i.e. to $J_j^n$, then one can apply the 
basic technical tool in the additive version and deform the $(p+1)$-tuple 
of matrices $A_j^0$ into one satisfying condition 2) from $2^0$ and every 
matrix $A_j^0$ being from the necessary orbit $J_j^n$. 

The lemma is proved.$\hspace{2cm}\Box$ 

\section{Proof of Theorem \protect\ref{genericrigid}\protect\label{n_s=1}}

In this section we consider only the case of equality in condition 
$(\alpha _n)$ (i.e. $n_s=1$). The necessity of conditions {\em i)} and 
{\em ii)} from Theorem~\ref{generic} for the existence of irreducible 
$(p+1)$-tuples of matrices $A_j$ or $M_j$ with generic 
eigenvalues was proved in Section~\ref{necessity}. Their sufficiency for 
the existence of such $(p+1)$-tuples for {\em some} generic eigenvalues was 
proved in Section~\ref{sufficiency}. There remains to be proved the 
existence of such $(p+1)$-tuples (for fixed Jordan normal forms $J_j^n$) 
for {\em all} generic eigenvalues. 

\subsection{Definitions and notation}

\begin{defi}
Two $(p+1)$-tuples of conjugacy classes 
$c_j$ each with non-resonant eigenvalues and with $\sum _j${\rm Tr}$c_j=0$ 
are said to be {\em similar} if they are 
obtained from one another by adding to equal eigenvalues of a given 
conjugacy class equal integers; the sum of all added integers (taking into 
account the multiplicities) is 0.
\end{defi} 

\begin{defi}
A $(p+1)$-tuple of conjugacy classes ($c_j$ or $C_j$) is 
{\em good} (resp. is {\em bad}) if the DSP is solvable for 
matrices $A_j\in c_j$ or for matrices $M_j\in C_j$ (resp. if not). 
A $(p+1)$-tuple of Jordan normal forms is 
{\em good} (resp. is {\em bad}) 
if there exists a good $(p+1)$-tuple of conjugacy classes defining the 
corresponding Jordan normal forms, with generic eigenvalues (resp. if not). 
If the $(p+1)$-tuple of Jordan normal forms 
is fixed, then the $(p+1)$-tuple of conjugacy classes is completely defined 
by the eigenvalues and we say that the eigenvalues are good or bad if 
the $(p+1)$-tuple of conjugacy classes is such. 
\end{defi} 

In the case of matrices $A_j$ and for a fixed $(p+1)$-tuple of Jordan normal 
forms denote by ${\bf C}^s$ the space of 
eigenvalues. (If the class $c_j$ has $s_j$ distinct eigenvalues, then 
$s=s_1+\ldots +s_{p+1}$). Denote by ${\bf C}'\simeq {\bf C}^{s-1}$ the 
subspace of ${\bf C}^s$ 
defined by the condition the sum of all eigenvalues (taking the 
multiplicities into account) to be 0. The PMV of the eigenvalues is simple, 
otherwise there exist no generic eigenvalues at all.

In accordance with the above 
definition, we call the points from ${\bf C}'$ {\em good} or {\em bad} if 
they define 
good or bad $(p+1)$-tuples of conjugacy classes. A point from ${\bf C}'$ is 
called {\em (strongly) generic} (resp. {\em non-resonant}) if it defines 
(strongly) generic (resp. non-resonant) eigenvalues.

Denote by $\Omega$ the set of good points of ${\bf C}'$. By 
Lemma \ref{allgen} (in which $L$ coincides with ${\bf C}'$; recall that 
$L$ is connected in the case of matrices $A_j$), the set 
$\Omega$ contains a Zariski open dense subset of ${\bf C}'$. The set 
$\Omega$ is constructible and invariant under multiplication by ${\bf C}^*$.

\subsection{Proof of Theorem \protect\ref{genericrigid} in the 
multiplicative version}



$1^0$. If the $(p+1)$-tuple of Jordan normal forms $J_j^n$ is good, then it 
is impossible to have 
the following situation: there exists a bad non-resonant strongly generic 
point 
$P\in {\bf C}'$ and every point from the set $\Theta _P$ consisting of $P$ and 
of all points in ${\bf C}'$ defining $(p+1)$-tuples of conjugacy classes 
similar to the ones defined by $P$ is also bad. Indeed, the 
constructibility of $\Omega$ implies that if the above situation takes place, 
then $\Omega$ cannot contain a Zariski open dense subset of the set 
of all generic eigenvalues of ${\bf C}'$; hence, $\Omega$ must be empty. 

$2^0$. Hence, every set $\Theta _P$ defined like above contains a good 
strongly 
generic non-resonant point. Every strongly generic non-resonant point from 
${\bf C}'$ defines a $(p+1)$-tuple of conjugacy classes in $GL(n,{\bf C})$ 
via the rule: if $c_j$ are the conjugacy classes in $gl(n,{\bf C})$ defined 
by the point and if $A_j\in c_j$, then $C_j$ are the conjugacy classes of 
the matrices $\exp (2\pi iA_j)$. 

$3^0$. All points from $\Theta _P$ define one 
and the same $(p+1)$-tuple of conjugacy classes $C_j\in GL(n,{\bf C})$, with 
generic eigenvalues. 
Hence, this $(p+1)$-tuple of conjugacy classes is good (indeed, the 
monodromy group of a fuchsian system with matrices-residua $A_j\in c_j$ 
is irreducible for any choice of the positions of the poles and one has 
$M_j\in C_j$ with $c_j$ as in $2^0$). On the other hand, 
for every $(p+1)$-tuple of conjugacy classes in $GL(n,{\bf C})$ with generic 
eigenvalues one can find a set $\Theta _P$ as above which  
defines this $(p+1)$-tuple of conjugacy classes. 

$4^0$. Hence, if a given 
$(p+1)$-tuple of Jordan normal forms is good, then for all possible generic 
eigenvalues 
there exist $(p+1)$-tuples of matrices $M_j$ with this $(p+1)$-tuple of 
Jordan normal forms.
This proves Theorem~\ref{genericrigid} in the multiplicative version.

\subsection{Proof of Theorem \protect\ref{genericrigid} in the additive 
version}

$1^0$. Consider a good generic point $D$ from ${\bf C}'$ such that all 
eigenvalues 
of all conjugacy classes $c_j$ are integer (hence, it is not strongly 
generic). All such points from ${\bf C}'$ 
cannot be bad because the constructibility of $\Omega$ would imply that 
$\Omega$ is empty.

\begin{lm}\label{L-Dan}
The monodromy operator of a fuchsian system the eigenvalues of whose 
matrices-residua define the point $D$ is upper-triangular up to conjugacy.
\end{lm}

The lemma is proved in the next subsection. 

The lemma implies that one can choose an initial value of the solution 
$X$ such that the monodromy group is upper-triangular, with matrices $M_j$ 
arbitrarily close to $I$ in some matrix norm. Hence, for matrices-residua 
close to the given ones defined by the point $D$ the monodromy operators 
$M_j$ will be all close to $I$. 
If these matrices-residua are with the same Jordan normal forms as the ones 
defined by $D$, then the 
monodromy group is defined by a point from a set $\tilde{D}\subset {\bf C}'$ 
containing a neighbourhood of $D$ in ${\bf C}'$.  
These points are also good -- to prove it one has to apply the basic 
technical tool in the additive version.

$2^0$. In \cite{L-D} I.A. Lappo-Danilevskii proves the following result:

{\em The monodromy operators of a fuchsian system are expressed as power 
series of its matrices-residua. These series are convergent if the residua 
are small enough and for such residua the map "residua" $\mapsto$ "monodromy 
operators" is a diffeomorphism of a neighbourhood ${\cal N}_0$ of 0 to a 
neighbourhood ${\cal N}_1$ 
of $I$. The initial data are $X|_{t=\infty}=I$ assuming 
that there is no pole at $\infty$.}

$3^0$. We give another formulation of the above result. Identify the space of 
$(p+1)$-tuples of matrices $A_j$ whose sum is 0 with $(gl(n,{\bf C}))^p$ 
(one defines only the first $p$ of them). Denote by ${\cal S}$ the unit 
sphere in $(gl(n,{\bf C}))^p$ when regarded as ${\bf R}^{2pn^2}$. Introduce 
coordinates in $(gl(n,{\bf C}))^p$ which are the union of some coordinates 
on ${\cal S}$ and $h\geq 0$. Consider the fuchsian system 

\[ \dot{X}=(\sum _{j=1}^{p+1}hA_j/(t-a_j))X\]
(without a pole at $\infty$, with $(A_1,\ldots ,A_{p+1})\in {\cal S}$) and 
its solution $X$ satisfying the condition $X(\infty )=I$. 

\begin{lm}\label{concreteL-D}
For $h$ small enough and $(A_1,\ldots ,A_{p+1})\in {\cal S}$ one has 
$M_j=I+h2\pi iA_j+o(h)$. The estimation is uniform in 
$(A_1,\ldots ,A_{p+1})\in {\cal S}$. 
\end{lm}
 
The lemma is proved in Subsection \ref{proofofconcreteL-D}. It implies that 
for $h$ small enough the map "residua" $\mapsto$ "monodromy 
operators" is a diffeomorphism of ${\cal N}_0$ to ${\cal N}_1$. 

$4^0$. The monodromy group ${\cal M}$ of a fuchsian system ($F_1$) with 
conjugacy classes of its matrices-residua corresponding to every strongly 
generic non-resonant 
point from $\tilde{D}$ admits a conjugation after which it will belong 
to ${\cal N}_1$, see $2^0$. 

By $3^0$, there exists a fuchsian system ($F_2$) with 
eigenvalues close to 0 with the same monodromy group. Its eigenvalues 
are shifted w.r.t. the ones of ($F_1$) by integers (and the $(p+1)$-tuples 
of conjugacy classes of the matrices-residua of $(F_1$) and ($F_2$) are 
similar). These integers are opposite to the eigenvalues of the 
$(p+1)$-tuple of conjugacy classes defined by the point $D$. 

$5^0$. Every strongly generic non-resonant point in ${\bf C}'$ which is close 
to $D$ is good and by 
the shift of eigenvalues defined in $4^0$ it defines a good strongly generic 
non-resonant point close to 0. The shift leaves the set of strongly generic 
non-resonant points invariant. Hence, all strongly generic non-resonant points 
close to 0 are good. 

This means that all generic points close to 0 are good. Indeed, a generic 
point being close 
to 0 means that it is strongly generic and non-resonant. 

$6^0$. By $5^0$, the set $\Omega$ contains the intersection of some 
neighbourhood of $0\in {\bf C}'$ with the set of generic points of 
${\bf C}'$. The set 
$\Omega$ being invariant under multiplication by ${\bf C}^*$ (it is defined 
by linear homogeneous inequalitites), it 
must contain the set of all generic points of ${\bf C}'$.

This proves Theorem \ref{genericrigid} in the additive version.

\subsection{Proof of Lemma \protect\ref{L-Dan}}

$1^0$. Suppose that the monodromy group is not triangularizable by 
conjugation. Then it can be conjugated to a block upper-triangular form, 
its restriction to at least one diagonal block $P$ of size $m>1$ being 
irreducible. Hence, the matrix algebra generated by the restriction of the 
monodromy matrices $M_j$ to the block $P$ is $gl(m,{\bf C})$ (the Burnside 
theorem) and, hence, there exists a polynomial without a constant term 
$s(M_1-I,\ldots ,M_{p+1}-I)$ in the matrices $M_j-I$ which is a matrix with 
at least two distinct eigenvalues. 

$2^0$. All points from ${\bf C}'$ sufficiently close to $D$ are good as well 
(this follows easily from the basic technical tool in the additive version). 
Hence, the same polynomial $s$ evaluated for $M_j$ corresponding to points 
close to $D$ is still a matrix with at least two distinct eigenvalues which 
will be close to two of the eigenvalues of $s$ evaluated at $D$ (denoted  
by $a$, $b$). Suppose that $a\neq 0$. Hence, for all points from 
${\bf C}'$ close to $D$ (denote their set by $D^*$) the polynomial $s$ 
evaluated at them has an eigenvalue $\lambda$ with $|\lambda |\geq |a|/2$.  

$3^0$. {\em There exists a constant $\delta >0$ such that $||s||\geq \delta$ 
for all points from $D^*$.} (Notice that this estimation is based only on the 
presence of distinct eigenvalues; the monodromy group is defined only up to 
conjugacy and the above estimation is valid for any of the possible 
definitions of the monodromy group.)

Suppose that such a constant $\delta$ does not exist. Then for each 
$\varepsilon >0$ there exist 
points from $D^*$ arbitrarily close to $D$ such that 
$||s||<\varepsilon$. One can choose as matrix norm the sum of the absolute 
values of all entries of a matrix. 

There holds the following lemma (well-known to specialists in 
numerical methods -- the lemma of diagonal domination):

\begin{lm}
If the module of any diagonal entry of a matrix from 
$gl(n,{\bf C})$ is greater than the sum of the modules of 
the non-diagonal entries of the same row, then the matrix is non-degenerate.
\end{lm}

Hence, if the modules of all entries of a matrix $A$ from $gl(n,{\bf C})$ 
are smaller than $|a|/4n$, then the matrix cannot have an eigenvalue 
$\lambda$ with $|\lambda |\geq |a|/2$ (because the matrix $A-\lambda I$ will 
satisfy the conditions of the lemma). This implies the existence of $\delta$ 
as above.

$4^0$. Fix a generic point $Q\in \Omega$ and consider the points $hQ$, 
$h\in [0,1]$. Denote their set by $\gamma$. For $h\neq 0$ they are generic 
and belong to $\Omega$; for $h\neq 0$ small enough 
they are strongly generic and non-resonant. Let the $(p+1)$-tuple 
$(A_1^*,\ldots ,A_{p+1}^*)$ be with eigenvalues defined by $Q$.  
Consider the fuchsian system with fixed poles and with $(p+1)$-tuple of 
matrices-residua $hA_j^*$. If one defines its monodromy group by fixing 
the initial point and initial value of $X$ one and the same for all $h$, 
then for $h$ small enough 
its monodromy operators $M_j(h)$ will be arbitrarily close to $I$. This 
follows from the continuous dependence of the solution on the parameter $h$ 
(when the solution is considered on any simply connected domain not 
containing a pole of the system) -- for $h=0$ one has $X=$const, $M_j=I$. 
Hence, the norm of the polynomial $s$ computed for these monodromy groups 
will be arbitrarily small when $h\rightarrow 0$. 

For $h\neq 0$ small enough the strong 
genericity of the eigenvalues of the matrices $hA_j^*$ implies that the 
monodromy groups of the systems are irreducible. They are rigid because 
$(\alpha _n)$ is an equality. Hence, they are unique up to conjugacy.

$5^0$. Denote by $\gamma '\subset {\bf C}'$ the segment $\gamma$ translated 
so that the 
point corresponding to $h=0$ be at $D$. For $h$ small enough the points from 
$\gamma '$ belong to $\Omega$. The translation of $\gamma$ into $\gamma '$ 
means that for every fixed value $h_0\neq 0$ of $h$ its corresponding 
$(p+1)$-tuple of conjugacy classes defined by the point  
$R=\gamma |_{h=h_0}$ will be replaced by a similar $(p+1)$-tuple defined by 
the point  
$R'=\gamma '|_{h=h_0}$. Hence, up to conjugacy, the monodromy groups of the 
fuchsian systems with matrices-residua whose eigenvalues are defined by the 
points $R$ and $R'$ coincide (when $h_0$ is small enough). Indeed, they 
correspond to similar $(p+1)$-tuples of conjugacy classes. 

This however is a contradiction with $2^0$ -- by choosing $h$ small enough 
the norm of $s$ computed for $R$ can become arbitrarily small which is 
impossible to happen for $R'$, see $3^0$.

The lemma is proved.$\hspace{2cm}\Box$

\subsection{Proof of Lemma \protect\ref{concreteL-D} 
\protect\label{proofofconcreteL-D}}

$1^0$. To compute the monodromy operators $M_j$ we fix the contours of 
integration. They begin at $\infty$, go along arcs $\eta _j$ to some 
points $b_j$ close to $a_j$, 
go around $a_j$ counterclockwise along the circumferences $\zeta _j$ passing 
through $b_j$ and centered at $a_j$, and 
then go back to $\infty$ along $\eta _j$. The points $b_j$ and the arcs 
$\eta _j$ are chosen such 
that there is no other pole of the system except $a_j$ on the closed discs 
$\Xi _j$ (where $\partial \Xi _j=\zeta _j$) and no pole at all on $\eta _j$.

$2^0$. For each $j$ the value at $b_j$ of the analytic continuation of the 
solution to 
the system with initial data $X(\infty )=I$ along the segment $\eta _j$ 
equals $I+Rh+o(h)$. This estimation is uniform in 
$(A_1,\ldots ,A_{p+1})\in {\cal S}$. Indeed, this follows from the smooth  
dependence of the solution $X$ on $h$ (there are no singularities of the 
system on $\eta _j$; for $h=0$ one has $X\equiv I$). 

$3^0$. Denote by $K_j$ the operators of local monodromy defined with 
initial data $X(b_j)=I$ and mapping this value onto the value of the 
analytic continuation of $X$ along $\zeta _j$. For $h$ small enough one has 
$K_j=I+h2\pi iA_j+o(h)$. 

Indeed, for $h$ small enough and for any 
$(A_1,\ldots ,A_{p+1})\in {\cal S}$ no two eigenvalues of any of the 
matrices $A_j$ differ by a non-zero integer. Hence, the solution to the 
system in some neighbourhood of $a_j$ can be represented in the form 

\[ X=(I+P(t,h))\exp (hA_j{\rm ln}(t-a_j))G_j\]
where $G_j\in GL(n,{\bf C})$ and $P$ is a Taylor series in $(t-a_j)$ whose 
terms are expressed through the entries of the matrices $A_j$, see \cite{Wa}. 
One has $P=O(h)$, $P(a_j,h)=0$. If $X(b_j)=I$, then 

\begin{equation}\label{G_j=} 
G_j=\exp (-hA_j{\rm ln}(b_j-a_j))(I+P(b_j,h))^{-1}
\end{equation} 
The matrix-function $(I+P(t,h))$ is holomorphically invertible on 
$\Xi _j\backslash a_j$ because $G_j$, $X$ and $\exp (hA_j{\rm ln}(t-a_j))$ 
are such (the latter two are multivalued). As $P(a_j,h)=0$, it is 
holomorphically invertible inside $\Xi _j$, i.e. at $a_j$ as well. 

$4^0$. The monodromy of the matrix-function 
$\exp (hA_j{\rm ln}(t-a_j))$ around 
$a_j$ equals $\exp (h2\pi iA_j)=I+h2\pi iA_j+o(h)$. One has 

\[ K_j=(G_j)^{-1}\exp (h2\pi iA_j)G_j=(G_j)^{-1}(I+h2\pi iA_j+o(h))G_j\] 
where $G_j$ is defined by (\ref{G_j=}). The factor 
$\exp (-hA_j{\rm ln}(b_j-a_j))$ commutes with $\exp (h2\pi iA_j)$. The factor 
$(I+P(b_j,h))^{-1}$ equals $I+Qh+o(h)$. Hence, 

\[ K_j=(I+Qh+o(h))^{-1}(I+h2\pi iA_j+o(h))(I+Qh+o(h))=I+h2\pi iA_j+o(h).\]

$5^0$. The operator $M_j$ is conjugate to $K_j$ and the conjugation is 
carried out by the matrix $I+Rh+o(h)$ from $2^0$. 
Hence, $M_j=I+h2\pi iA_j+o(h)$.

The lemma is proved.$\hspace{2cm}\Box$

\section{Proof of Theorem \protect\ref{genericbis}
\protect\label{proofofgenericbis}}

$1^0$. Prove the necessity. Denote by $J_j^n$ the Jordan normal forms of the 
matrices $M_j$ and 
by $J_j^{n,1}$ their corresponding diagonal Jordan normal forms, see 
Subsection~\ref{diag}. The PMV of the $(p+1)$-tuple of Jordan normal forms 
$J_j^{n,1}$ is simple; this follows from $d=1$ and from the construction of 
$J_1$ after $J_0$ in Subsection~\ref{diag}. 

$2^0$. Apply 
the basic technical tool in the multiplicative version. Look 
for matrices $\tilde{M}_j$ of the form 

\[ \tilde{M}_j(\varepsilon )=(I+\varepsilon X_j(\varepsilon ))^{-1}
Q_j^{-1}(G_j^0+\varepsilon G_j^1)Q_j(I+\varepsilon X_j(\varepsilon ))\]
where $M_j=Q_j^{-1}G_j^0Q_j=\tilde{M}_j(0)$, $G_j^0$ being Jordan matrices 
and $G_j^1$ 
being diagonal matrices.  The matrices $G^0_j$ and $G^1_j$ are chosen with the 
properties respectively 
of $G^0$ and $G^1$ from part 2) of Proposition~\ref{sameorbit}. 
For $\varepsilon$ small enough the eigenvalues 
of the matrices $\tilde{M}_j$ are generic and their Jordan normal forms 
(when $\varepsilon \neq 0$) are $J_j^{n,1}$, see Proposition~\ref{sameorbit}. 

$3^0$. Theorem \ref{basicres} gives the necessary and sufficient conditions 
for the existence of irreducible $(p+1)$-tuples of matrices 
$\tilde{M}_j\in J_j^{n,1}$ with generic eigenvalues. 
It follows from the theorem that conditions {\em i)} and {\em ii)} of 
Theorem~\ref{generic} hold for the $(p+1)$-tuple of Jordan normal forms 
$J_j^n$, see Remarks~\ref{psiinvariant}. Hence, 
conditions {\em i)} and {\em ii)} from Theorem~\ref{generic} are necessary 
for the existence of 
irreducible $(p+1)$-tuples of matrices $M_j$. 

This proves the necessity.

$4^0$. Prove the sufficiency. If the conditions of Theorem~\ref{basicres} 
hold for the diagonal Jordan normal forms $J_j^{n,1}$ (recall that their PMV 
is simple), then there exist 
irreducible $(p+1)$-tuples of matrices $A_j^2\in J_j^{n,1}$ satisfying 
(\ref{A_j}). One can choose their eigenvalues $\lambda _{k,j}$ (presumed to 
be generic and to form for each $j$ a normalized chain, see 
Subsection~\ref{NC}) such 
that 
$\exp (2\pi i\lambda _{k,j})=\sigma _{k,j}$.  
Hence, for each $j$ the Jordan normal 
form of $M_j$ will be either $J_j^n$ or some Jordan normal form 
subordinate to it, see Lemma~\ref{JNFofM}. By part 2) of  
Lemma~\ref{J^0impliesJ^1}, there exist irreducible $(p+1)$-tuples of 
matrices $M_j\in J_j^n$, with generic eigenvalues. 

The theorem is proved.$\hspace{2cm}\Box$

\section{Proofs of Propositions \protect\ref{d_jr_j} 
and \protect\ref{necessarynongen}\protect\label{proofofd_jr_j}}

${\bf 1^0}$. We prove $(\alpha _n)$ in $1^0$ -- $4^0$ and $(\beta _n)$ in 
$5^0$. To find the dimension of the representation defined by the 
matrices $A_j$ in $(sl(n,{\bf C}))^p$  one has first to restrict oneself 
to the cartesian product $c_1\times \ldots \times c_p$ 
of the orbits 
of $A_1$, $\ldots$, $A_p$ considered as varieties in each of the $p$ copies 
of $sl(n,{\bf C})$.

${\bf 2^0}$. The algebraic variety ${\cal V}$ defined in $sl(n,{\bf C})^p$ 
by the orbits of 
$A_1$,$\ldots$,$A_{p+1}$ 
is the projection in $c_1\times \ldots \times c_p$ of the intersection of 
the two 
varieties in $c_1\times \ldots \times c_p\times sl(n,{\bf C})$ (the matrix 
$A_{p+1}$ 
corresponds to the factor $sl(n,{\bf C})$): the graph of the 
mapping 

\[ (c_1\times \ldots \times c_p)\ni (A_1,\ldots ,A_p)
\mapsto A_{p+1}=-A_1-\ldots -A_p\in sl(n,{\bf C})\]
and $c_1\times \ldots \times c_p\times c_{p+1}$. This intersection is 
transversal which implies the smoothness of the 
variety ${\cal V}$ (this can be proved by analogy with 1) of Theorem 2.2 
from \cite{Ko5}). Thus 

\[ {\rm dim}\, {\cal V}=(\sum _{j=1}^p{\rm dim}\, c_j)-[(n^2-1)-
{\rm dim}\, c_{p+1}]\]
(here $(n^2-1)-{\rm dim}\, c_{p+1}=$codim$_{sl(n,{\bf C})}c_{p+1}$). 
Hence, dim${\cal V}=\sum _{j=1}^{p+1}$dim$\, c_j-n^2+1$. 

${\bf 3^0}$. In order to obtain the dimension of the representation defined 
by  $(A_1,\ldots ,A_{p+1})$ one has to factor out the possibility to 
conjugate the $(p+1)$-tuple $(A_1,\ldots ,A_{p+1})$ with matrices from 
$SL(n,{\bf C})$. No such non-scalar matrix commutes with all the 
matrices $(A_1,\ldots ,A_{p+1})$ due to the irreducibility of the 
$(p+1)$-tuple and to Schur's lemma. Thus the dimension of the 
representation equals 

\[ {\rm dim}\, {\cal V}-n^2+1=\sum _{j=1}^{p+1}{\rm dim}\, c_j-2n^2+2\]
(where $-n^2+1=-{\rm dim}\, SL(n,{\bf C})$). 
When this dimension is negative, 
then the representation does not exist.
In the case of $(SL(n,{\bf C}))^p$ the mapping    
$(A_1,\ldots ,A_p)$ $\mapsto A_{p+1}=-A_1-\ldots -A_p$
from ${\bf 2^0}$ has to be replaced by the mapping 
\[(M_1,\ldots ,M_p)\mapsto M_{p+1}=(M_1\ldots M_p)^{-1}~.\] 

This proves $(\alpha _n)$. 

$4^0$. Prove $(\beta _n)$. If 

\[ r_1+\ldots +\hat{r}_j+\ldots r_{p+1}<n~,\]
then for suitable eigenvalues $\lambda _i$ 

\[ {\rm codim} \cap _{i\neq j}{\rm Ker}(A_i-\lambda _iI)<n\] 
and this non-trivial subspace contradicts the irreducibility assumption. 

The proposition is proved.$\hspace{2cm}\Box$

{\em Proof of Proposition \ref{necessarynongen}:}

$1^0$. The $(p+1)$-tuples of matrices $A_j$ and $A_j-b_jI$ being 
simultaneously irreducible we prove the proposition for $b_j=0$, 
$j=1,\ldots ,p+1$. Consider the maps 

\[ \begin{array}{lcl}
\sigma :(X_1,\ldots ,X_{p+1})\mapsto (A_1X_1,\ldots ,A_{p+1}X_{p+1})&& 
{\rm and}\\
\tau :(X_1,\ldots ,X_{p+1})\mapsto 
(A_1X_1+\ldots +A_{p+1}X_{p+1})&,&X_j\in {\bf C}^n\end{array}\] 
Prove that dim(Im$\, \sigma$)$ \geq 2n$ 
(which amounts to proving part 1) of the lemma). 

Denote by $S_1$ the subspace of ${\bf C}^{(p+1)n}$ where 
$S_1=\{ (X,\ldots ,X), X\in {\bf C}^n\}$, dim$S_1=n$; one has 
$\tau (S_1)=\{ 0\}$ (because $A_1+\ldots +A_{p+1}=0$). 
The irreducibility of $(A_1,\ldots ,A_{p+1})$ implies that 

A) the image of $\tau$ is the whole space ${\bf C}^n$; hence, there exists 
a subspace $S_2\subset {\bf C}^{(p+1)n}\simeq \{ (X_1,\ldots ,X_{p+1})\}$, 
dim$\, S_2=n$, such that $\tau (S_2)={\bf C}^n$;

B) $S_1\cap \, $Ker$\, \sigma =\{ 0\}$.

Hence, 

a) $S_1\cap S_2=\{ 0\}$ (otherwise a vector from $S_2$ belongs to 
Ker$\, \tau$ which contradicts $\tau (S_2)={\bf C}^n$ and dim$\, S_2=n$), 
dim$(S_1\oplus S_2)=2n$ and 

b) the image of $S_1\oplus S_2$ under $\sigma$ is of dimension $2n$. Indeed, 
the image of $S_1$ is of dimension $n$ and belongs to Ker$\, \tau$ whereas 
the one of $S_2$ (also of dimension $n$) is transversal to Ker$\, \tau$. 
This proves 1).

$2^0$. Prove 2). If the monodromy group 
generated by the matrices $M_j$ is 
irreducible, then so is the algebra generated by the matrices 
$M_j$, $j=1,\ldots ,p+1$; hence, by 
the matrices 
$N_j=(b_{j+1}\ldots b_{p+1})M_1\ldots M_{j-1}(M_j-b_jI)$ as well (recall 
that $b_1\ldots b_{p+1}=1$).

One has $N_1+\ldots +N_{p+1}=0$ and rk$N_j=$rk$(M_j-b_jI)$. Hence, 
part 2) follows from part 1) by setting $A_j=N_j$. 

The proposition is proved.$\hspace{2cm}\Box$

\end{document}